\documentclass[12pt,reqno]{amsart}
\usepackage{amsfonts}
\usepackage{amssymb}
\usepackage{amsthm}
\usepackage{upref}
\usepackage{enumerate}

\usepackage{amscd}

\makeatletter 

\def\LaTeX{\leavevmode L\raise.42ex
    \hbox{\kern-.3em\size{\sf@size}{0pt}\selectfont A}\kern-.15em\TeX}
 \makeatother
 
  \newcommand{\arccosh}{\operatorname{arccosh}}

\DeclareMathOperator{\clos}{clos}

\numberwithin{equation}{section}

\newtheorem{lemma}{Lemma}[section]
\newtheorem{theorem}[lemma]{Theorem} 
\newtheorem{corollary}[lemma]{Corollary}
\newtheorem{proposition}[lemma]{Proposition}

\theoremstyle{definition}

\newtheorem{definition}[lemma]{Definition}
\newtheorem{example}[lemma]{Example}

\newtheorem{remark}[lemma]{Remark}

\renewcommand{\det}{\operatorname{Det}}

\newcommand{\tr}{\operatorname{Tr}}

\newcommand{\diag}{\operatorname{diag}}

 \newcommand{\supp}{\operatorname{supp}}
  \newcommand{\e}{\eqref}

\newcommand{\q}{\quad}

\newcommand{\ov}{\overline}
\newcommand{\wt}{\widetilde}
\newcommand{\z}{\zeta}
\newcommand{\ti}{\tilde}

\newcommand{\la}{\langle}
\newcommand{\ra}{\rangle}

\renewcommand{\d}{\delta}

   \newcommand{\sgn}{\operatorname{sgn}}
  \newcommand{\ran}{\operatorname{Ran}}

\renewcommand\Im{\operatorname{Im}}
\renewcommand\Re{\operatorname{Re}}

\newenvironment{pf}{\begin{proof}}{\end{proof}}

\def\qqq{\mathrel{\subset\mkern-15mu\lower.38ex\hbox{${\scriptscriptstyle\rightarrow}$}}}

\let\goth\mathfrak
\let\cal\mathcal

\let\Bbb\mathbb
 
\begin{document}
\title 
[Asymptotic  behavior of orthogonal polynomials]
{Spectral analysis of Jacobi operators and 
\\
asymptotic  behavior of orthogonal polynomials} 
\author{ D. R. Yafaev  }
\address{   Univ  Rennes, CNRS, IRMAR-UMR 6625, F-35000
    Rennes, France, SPGU, Univ. Nab. 7/9, Saint Petersburg, 199034 Russia, and    Sirius Math.  Center, Olympiysky av. 1, Sochi, 354349 Russia}
\email{yafaev@univ-rennes1.fr}


 

\thanks {Supported by  project   Russian Science Foundation   17-11-01126 and  Sirius Univ. of Science and Technology (project `Spectral and Functional Inequalities of Math. Phys. and their Appl.')}

\begin{abstract}
We find and discuss  asymptotic formulas for orthonormal polynomials $P_{n}(z)$  with  recurrence coefficients $a_{n}, b_{n}$. 
Our main goal  is to consider the case where off-diagonal elements  $a_{n}\to\infty$ as $n\to\infty$. Formulas obtained are essentially different for relatively small and large diagonal elements $b_{n}$. 

Our analysis is intimately linked with
 spectral theory  of    Jacobi operators $J$ with  coefficients $a_{n}, b_{n}$ and a study of the corresponding second order difference equations. 
    We introduce the Jost solutions $f_{n}(z)$, $n\geq -1$,  of such equations by a condition for $n\to\infty$ and suggest an Ansatz for  them  playing the role of the semiclassical Liouville-Green Ansatz for   solutions of the   Schr\"odinger equation.  
  This allows us to study the spectral structure of Jacobi operators and their eigenfunctions $P_{n}(z)$ by traditional methods of spectral theory developed for differential equations. In particular, we express all coefficients in asymptotic formulas for $P_{n}(z)$ as $n \to\infty$ in terms of the Wronskian of the solutions $  P_{n} (z) $ and $ f_{n} (z)$.   The formulas obtained for    $P_{n}(z)$  generalize the  asymptotic formulas for the classical  Hermite   polynomials where  $a_{n}=\sqrt{(n+1)/2}$ and $b_{n}=0$.
  
  The spectral structure of    Jacobi operators $J$ depends crucially on a rate of growth of the off-diagonal elements  $a_{n} $ as $n\to\infty$.  If the Carleman condition is satisfied, which,  roughly  speaking, means that $a_{n}= O (n)$, and the diagonal elements  $b_{n} $ are small compared to $a_{n} $, then $J$ has the absolutely continuous spectrum covering the whole real axis. We obtain an expression for the corresponding spectral measure  in terms of the boundary values $| f_{-1}(\lambda\pm i0)|$ of the Jost solutions.  On the contrary, if 
  the Carleman condition is violated, then the spectrum of $J$ is discrete.
  
  We also review the case of stabilizing recurrence coefficients when $a_{n}$ tend to a positive constant and $b_{n}\to 0$ as $n\to\infty$. It turns out that the cases of stabilizing and increasing recurrence coefficients  can be treated in an essentially same way.
   \end{abstract}

 \maketitle

 \tableofcontents
 

\section{Overview}



\subsection{Two definitions of orthogonal   polynomials }

There  are two a priori different definitions of
orthogonal, or more precisely orthonormal,  polynomials  $P_{n } (z)$. The first definition proceeds from two sets of coefficients  $a_n >0$, $b_n=\bar{b}_n $ and polynomials  $P_{n } (z) $ are determined
by  a recurrence relation
   \begin{equation}
 a_{n-1} P_{n-1} (z) +b_{n} P_{n } (z) + a_{n} P_{n+1} (z)= z P_n (z), \q n\in{\Bbb Z}_{+}=\{0,1,\ldots\}, \q z\in{\Bbb C}, 
\label{eq:RR}\end{equation}
complemented by   boundary conditions 
 \begin{equation}
 P_{-1 } (z) =0, \q P_0 (z) =1.  
 \label{eq:RR1}\end{equation}
  Finding  $P_{n } (z)$ for $n=1,2, \ldots$ successively from \e{eq:RR}, we see that $P_{n } (z)$ is a polynomial 
with real coefficients of degree $n$:
  \begin{equation}
 P_{n } (z)= k_{n}(z^n+r_{n}z^{n-1}\cdots ).
\label{eq:RRx}\end{equation}
Substituting this expression into \e{eq:RR} and comparing the coefficients at $z^{n+1}$, we see that necessarily
   \begin{equation}
k_{n}/k_{n+1}= a_{n}\q\mbox{whence}\q k_{n}= (a_{0}a_{1}\cdots a_{n-1})^{-1}.
\label{eq:RRx1}\end{equation}
Similarly, comparing the coefficients at $z^{n}$, we see that  
   \begin{equation}
 r_{n}-r_{n+1} = b_{n}  \q\mbox{whence}\q r_{n}= -\sum_{m=0}^{n-1} b_{m}.
\label{eq:RRx2}\end{equation}
Relations \e{eq:RRx1} and \e{eq:RRx2}  allow one to recover the recurrence coefficients $a_{n}$, $b_{n}$  from the polynomials $P_{n} (z)$ satisfying  \e{eq:RR}, \e{eq:RR1}.

Another possibility is to define orthonormal polynomials via a measure $d\rho(\lambda)$ on ${\Bbb R}$. It is supposed that all moments
   \begin{equation}
s_{n}: = \int_{-\infty}^\infty \lambda^n d\rho(\lambda)< \infty, \q s_{0}=1,
\label{eq:mo}\end{equation}
  and that the support of $d\rho(\lambda)$ is infinite. One defines  polynomials  $ P_{n } (z)$ by the Gram-Schmidt orthonormalization of the monomials $1,\lambda,\ldots, \lambda^n,\ldots$
in the space $L^2 ({\Bbb R};d\rho)$  with the scalar product denoted $\la \cdot, \cdot \ra$. Then
      \begin{equation}
\la P_{n}, P_{m}\ra =\int_{-\infty}^\infty P_{n}(\lambda) P_{m}(\lambda) d\rho(\lambda) =\d_{n,m};
\label{eq:mo1}\end{equation}
as usual, $\d_{n,n}=1$ and $\d_{n,m}=0$ for $n\neq m$. Such polynomials are defined up to a sign. We accept that the coefficients  $k_{n}$ in \e{eq:RRx} are positive. We emphasize that all orthogonal polynomials considered in this paper are normalized. The measure $d\rho(\lambda)$ is often called the orthogonality measure for the polynomials   $P_{0}, P_{1}, \ldots, P_{n },\ldots$.

  The following statement shows that this definition of orthonormal polynomials implies the first one.  
  
      \begin{proposition}\label{Fav}
 The   polynomials $P_{n}(\lambda)$  defined by equalities \e{eq:mo1} satisfy  recurrence relation \e{eq:RR} where  
 \begin{equation}
 a_{n}=\int_{-\infty}^\infty \lambda P_{n } (\lambda)P_{n+1 } (\lambda) d\rho(\lambda)> 0
\label{eq:RP}\end{equation}
and
 \begin{equation}
 b_{n}=\int_{-\infty}^\infty \lambda P_{n }^2 (\lambda) d\rho(\lambda).
\label{eq:RP1}\end{equation}
    \end{proposition}

    Proposition~\ref{Fav}  as well as Propositions~\ref{adj}, \ref{simple}, \ref{FAV} stated in the next subsection are checked in Appendix~A.

     There are three specific classes of polynomials, Jacobi, Hermite and Laguerre, named classical. To  a some extent, all modern studies of general orthogonal polynomials can be considered as far reaching generalizations of the results well known for the classical polynomials. Therefore we   discuss these important special cases in Appendix~B.


 \subsection{Jacobi operators }

Now we proceed from recurrence relations \e{eq:RR}, \e{eq:RR1} for the polynomials   $P_{0}, P_{1}, \ldots, P_{n },\ldots$ and construct their orthogonality measure. To that end, we introduce a 
 Jacobi matrix 
  \begin{equation*}
{\cal J} = 
\begin{pmatrix}
 b_{0}&a_{0}& 0&0&0&\cdots \\
 a_{0}&b_{1}&a_{1}&0&0&\cdots \\
  0&a_{1}&b_{2}&a_{2}&0&\cdots \\
  0&0&a_{2}&b_{3}&a_{3}&\cdots \\
  \vdots&\vdots&\vdots&\ddots&\ddots&\ddots
\end{pmatrix}   .
\end{equation*}
If    $u= (u_{0}, u_{1}, \ldots)^\top=: (u_{n})$ is a column,
  then 
\begin{equation}
( {\cal J} u) _{0} =  b_{0} u_{0}+ a_{0} u_{1}, \q  
( {\cal J} u) _{n} = a_{n-1} u_{n-1}+b_{n} u_{n}+ a_{n} u_{n+1} \q \mbox{for}\q n\geq 1  ,
 \label{eq:ZP+1}\end{equation}
and equation \e{eq:RR} with the boundary condition $P_{-1 } (z) =0$ is
 equivalent to the equation ${\cal J} P( z )= z P(z)$ for the vector $P(z)=(P_{0}(z), P_1(z),\ldots)^{\top}$.  Thus,  $P(z)$ is an ``eigenvector" of the matrix $\cal J$ corresponding to an ``eigenvalue" $z$. 
  

Let us now consider Jacobi operators  defined   by formula \e{eq:ZP+1} in 
 the space $\ell^2 ({\Bbb Z}_{+})$. 
A minimal Jacobi operator $J_{\rm min}$ is  defined   by the equality $J_{\rm min} u= {\cal J} u$ on a set $\cal D \subset \ell^2 ({\Bbb Z}_{+})$ of vectors $u=  (u_{n}) $ with only a  finite number of non-zero components  $u_{n}$. Note that $J_{\rm min}:  {\cal D}\to  {\cal D}$.  A maximal  operator $J_{\rm max}$ is given by the  same formula $J_{\rm max} u= {\cal J} u$ on the set $ {\cal D} (J_{\rm max})$  of all vectors $u \in \ell^2 ({\Bbb Z}_{+})$ such that ${\cal J}u \in \ell^2 ({\Bbb Z}_{+})$. 
Evidently, the operator $J_{\rm min} $ is
symmetric  in the space $\ell^2 ({\Bbb Z}_{+})$. We also have the following assertion.  


   \begin{proposition}\label{adj}
    The adjoint operator    $J^*_{\rm min} = J_{\rm max}$.
 \end{proposition} 
 
It is easy to see that $J_{\rm min} $ extends to a bounded operator  on $\ell^2 ({\Bbb Z}_{+})$ if and only if the sequences $ (a_{n} )\in \ell^\infty ({\Bbb Z}_{+})$ and $ (b_{n} )\in \ell^\infty ({\Bbb Z}_{+})$. In this paper we  are particularly interested in the case $a_{n}\to\infty$ when Jacobi  operators are unbounded.

 Since the operator $J_{\rm min}  $  commutes with the complex conjugation,  its deficiency indices are equal so that $J_{\rm min}  $ has self-adjoint extensions. Actually, the deficiency indices of $J_{\rm min}  $ are   either  $(0,0)$ or  $(1,1)$. In the first case the operator $J_{\rm min}  $ is essentially  self-adjoint, that is,  its closure $\clos J_{\rm min}=J_{\rm max}$. 

    A link of orthonormal polynomials with Jacobi operators  relies on the  following observation.
 
    \begin{proposition}\label{simple}
   Let $e_{0}, e_{1},\ldots, e_{n},\ldots$ be the canonical basis in the space $\ell^2({\Bbb Z}_{+})$.  Then
   \begin{equation}
 e_{n} =P_{n} (J_{\rm min})e_{0}.
\label{eq:Phi3}\end{equation} 
      \end{proposition}
      
 Let $J$ be an arbitrary self-adjoint extension of the operator  $J_{\rm min}$, and let  $E_{J}(\lambda)$    be its spectral family.
      According to Proposition~\ref{simple}  the spectrum of the operator $J$ is simple with $e_{0} = (1,0,0,\ldots)^{\top}$ being a generating vector. Therefore  it is natural to define   the   spectral measure of $J$ by the equality 
      \begin{equation}
d\rho_{J}(\lambda)=d\la E_{J}(\lambda)e_{0}, e_{0}\ra .
\label{eq:UD1}\end{equation}

Let us now state a version of Favard's  theorem.

    \begin{proposition}\label{FAV}
      Let      polynomials $P_{n}(\lambda)$ satisfy   recurrence relation \e{eq:RR} and boundary conditions \e{eq:RR1}, and let $J$ be an arbitrary    self-adjoint extension  of the   operator $J_{\rm min}$ given  by equalities  \e{eq:ZP+1}  on the set ${\cal D}\subset \ell^2 ({\Bbb Z}_{+})$. Define the measure $d\rho_{J}$ by formula \e{eq:UD1} and set
        \begin{equation}
(\Phi e_{n})(\lambda)= P_{n} (\lambda).
\label{eq:Phi1}\end{equation}
Then the operator
  \begin{equation}
\Phi : \ell^2({\Bbb Z}_{+})\to L^2({\Bbb R}; d\rho_{J}) 
\label{eq:Phi}\end{equation}
is unitary and enjoys the intertwining property
   \begin{equation}
(\Phi J f) (\lambda)= \lambda(\Phi  f) (\lambda), \q \forall f\in {\cal D} (J) .
\label{eq:Phi2}\end{equation} 
In particular, the spectrum of the operator $J $ is simple.
    \end{proposition}
    
     \begin{corollary}\label{FAVs}
     
   \begin{enumerate}[{\rm(i)}]
 
       The following properties are true:
       
 \item The  polynomials  $P_{n}(\lambda)$ are orthogonal and normalized  in the spaces $L^2 ({\Bbb R};d\rho_{J})$, that is, relations  \e{eq:mo1} with $d\rho(\lambda)=d\rho_{J}(\lambda)$  are satisfied.
 
 \item    The set of all polynomials is dense in  the spaces $L^2({\Bbb R}; d\rho_{J}) $.
    \end{enumerate}
        \end{corollary}
    
    Conversely, suppose that  a measure $d\rho(\lambda)$ satisfying assumptions \e{eq:mo} is given.
    One first constructs  orthonormal polynomials $P_{n} (\lambda)$
    satisfying equalities \e{eq:mo1} and then   defines the recurrence coefficients $a_{n}$, $b_{n}$ by formulas \e{eq:RP} and \e{eq:RP1}. Let $J_{\rm min}$ be
     the minimal Jacobi operator   with these coefficients.  Taking its arbitrary self-adjoint extension $J$, one defines its spectral measure $d\rho_{J}(\lambda)$ by formula \e{eq:UD1}.  It turns out that $d\rho_{J}(\lambda)=d\rho(\lambda)$ if  the  operator $J_{\rm min}$ is essentially self-adjoint (then $J= \clos J_{\rm min}$). However, this is not true in general; an example of such phenomenon is given in Sect.~1.4 -- see the Freud weight \e{eq:Freud} for $\beta<1$. 
     

 This situation  can also  be described   in terms of solutions of the Jacobi  equation  
  \begin{equation} 
  a_{n-1} u_{n-1} (z)+ b_{n} u_{n} (z)+ a_{n}  u_{n+1} (z)=z u_{n} (z),\q n\geq 1. 
  \label{eq:Jy}\end{equation}
  The Weyl theory developed by him for differential equations can be naturally adapted to difference equations (see, e.g., the book \cite{AKH} and reference therein).
  Similarly to differential equations, equation \e{eq:Jy} for $\Im z \neq 0$ always has a non-trivial solution  $f = (f_{n})\in \ell^2 ({\Bbb Z}_{+})$. This solution is either unique (up to a constant factor) or all solutions of equation \e{eq:Jy} belong to $\ell^2 ({\Bbb Z}_{+})$.
The first instance is known as the limit point case and the second one --  as the limit circle case.  It turns out that 
 the operator $J_{\rm min}$ is essentially self-adjoint if and only if the limit point case occurs. In the limit circle case,  the operator $J_{\rm min}$ has deficiency indices $(1,1)$.
 
  Note  also that an operator $J_{\rm min}$ is essentially self-adjoint if and only if the corresponding moment problem is determinate. This is discussed in Appendix~D.
  
  Thus, in the limit point case there is the one-to-one correspondence  between Jacobi coefficients $(a_{n}, b_{n})$ and measures $d\rho(\lambda)$. The reconstruction  of the coefficients $(a_{n}, b_{n})$ of a Jacobi operator $J$ from its spectral measure  $d\rho(\lambda)$  is known as the inverse spectral problem. This procedure (see Proposition~\ref{Fav}) is stated  in terms of the orthonormal polynomials and is quite explicit.  This looks quite different from the case of differential operators where a reconstruction of the coefficient $b(x) $ of a differential operator $H=D^2 + b(x)$, $D=-id/dx$,  from its spectral measure requires a solution of the integral Gelfand-Levitan equation \cite{Ge-Le}. Nevertheless,  finding one-to-one correspondences between classes of the coefficients $(a_{n}, b_{n})$ and the measures $d\rho(\lambda)$ (the characterization problem) is as difficult as for differential operators. 
  
 
Let us, finally,  mention an elementary but useful fact. 

    \begin{proposition}\label{refl}
       If orthonormal polynomials $P_{n}(z)$   are constructed by coefficients $(a_{n}, b_{n})$ and 
$\wt{P}_{n}(z)$  correspond to the coefficients $(a_{n}, - b_{n})$, then 
\[
\wt{P}_{n}(z)= (-1)^{n} P_{n}(-z).
\]
If $\wt{J}$ is the Jacobi operator with matrix elements $a_{n}$, $-b_{n}$, then $\wt{J}=-{\cal U}^* J {\cal U}$ where the unitary operator ${\cal U}$ is defined by $({\cal U} u)_{n}=(-1)^n u_{n}$ for $n\in{\Bbb Z}_{+}$.  In particular, if $b_{n}=0$ for all $n$, then the operators $J$ and $-J$ are unitarily equivalent.  
The corresponding spectral measures are linked by the relation
$d\ti{\rho} (\lambda)= d {\rho} (-\lambda)$.
\end{proposition}

  
 
 The comprehensive presentation of the results  described shortly  here  can be found in  the books \cite{AKH, Chihara,   Schm, Asshe} and the surveys \cite{Lub,Simon, Tot}.



  
    \subsection{Asymptotics of orthogonal polynomials. Stabilizing coefficients}

We are interested in an asymptotic behavior of the polynomials  $P_{n } (z)$ as $n\to\infty$. Asymptotic properties of $P_n (z)$ can be deduced either from the recurrence coefficients $a_{n}$, $b_{n}$ defining $P_{n } (z)$  or from the corresponding orthogonality  measure $d\rho (\lambda)$:
 
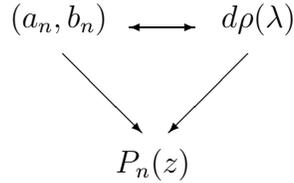
 \begin{picture}(150,80)
 \put(100,55){$(a_{n} , b_{n})$}
\put(180,55){$d\rho(\lambda)$}
\put(170,55){\vector(-1,0){25}}
\put(145,55){\vector(1,0){25}}
\put(120,45){\vector(1,-1){30}}
\put(190,45){\vector(-1,-1){30}}
\put(140, 0){$P_{n} (z)$}
\end{picture}

 \bigskip  \bigskip

Asymptotic formulas are very well known    for the classical, that is, Jacobi, Hermite and Laguerre,  polynomials (see   Appendix~B), but the first general result is probably due to S.~Bernstein 
(see his pioneering paper \cite{Bern} or  Theorem~12.1.4 in the G.~Szeg\H{o} book~\cite{Sz}). These results were stated in terms of the measure $d\rho (\lambda)$.
It was required      that $\supp\rho\subset [-1,1]$, the measure is absolutely continuous,  that is
  \begin{equation}
d\rho (\lambda)= \tau (\lambda) d\lambda ,
\label{eq:Jac}\end{equation}
   and     the weight   $\tau (\lambda)$ satisfies certain regularity conditions, in particular, it  does not tend to zero  too rapidly as $\lambda\to\pm 1$. Bernstein's result can be considered as a far reaching generalization of the asymptotic formulas for Jacobi polynomials determined by  spectral measure   \e{eq:Jac1}.


Alternative line of research where asymptotic formulas for the polynomials $P_n (z)$ are deduced from
assumptions on the coefficients $a_{n}, b_{n}$  was initiated in the book by P.~G.~Nevai \cite{Nev} who considered the case
 \begin{equation}
\lim_{n\to\infty} (a_{n } -1/2) =  \lim_{n\to\infty}b_{n }  = 0 .
\label{eq:comp}\end{equation}
Then
  the ``perturbation" $V =J-J_{0}$  of the ``free'' Jacobi operator $J_{0}$ with the elements $a_{n}= 1/2$, $b_{n}=0$ is compact. Therefore
  the essential spectrum of the operator $J$ coincides with the interval $[-1,1]$, and its discrete spectrum  consists of simple eigenvalues accumulating, possibly, to the points $1$ and $-1$. Under the only assumption 
  \e{eq:comp} the spectral structure of the Jacobi operator $J$ can be quite wild. Therefore Nevai supposed that   \begin{equation}
\sum_{n=0}^\infty \big(| a_{n} -1/2|+|b_{n}|\big)<\infty . 
\label{eq:Tr}\end{equation}
A more general case 
  \begin{equation}
\sum_{n=0}^\infty \big(| a_{n+1} - a_{n} |+|b_{n+1} - b_{n}|\big)<\infty 
\label{eq:LR}\end{equation}
was studied somewhat later in  \cite{ Mate}.
Under assumptions  \e{eq:comp}, \e{eq:LR} the  measure $d\rho(\lambda)$ is absolutely continuous on the interval $(-1,1)$ and the  weight $ \tau (\lambda) $ is a continuous and strictly positive function. The corresponding polynomials  $P_n (z)$  satisfy asymptotic relations generalizing formulas \e{eq:GG}, \e{eq:GGcomp} for the Jacobi polynomials.

The paper \cite{Mate} relies on specific methods of orthogonal polynomials theory.
Its initial point   is the relation 
  \[
 \lim_{n\to\infty} \frac{P_{n-1}(z)}{P_{n}(z)}=
 z -\sqrt{z^{2}-1}, 
 \q \Im z\neq 0, \q \sqrt{z^{2}-1}>0\q\mbox{for}\q z>1,
  \]
    established earlier by P.~Nevai in the book \cite{Nev}, Theorem~4.1.13, and
    improving one of  Poincar\'e's theorems. 
    
    Conditions \e{eq:Tr} and \e{eq:LR} are very precise. Indeed, as shown in \cite{Nab} (see also the preceding paper \cite{Pos}),  there exist coefficients   $b_{n}$  decaying only slightly worse than $n^{-1}$ and oscillating as $n\to\infty$ such that the point spectrum of the corresponding Jacobi operator $J$ with $a_{n}= 1/2$ is dense in $[-1,1]$. In this case the limiting absorption principle for the operator $J$
does not of course hold.

 This discussion is continued in Sect.~8.4.

 
    \subsection{Asymptotics of orthogonal polynomials. Increasing coefficients}

Passing to the case of unbounded Jacobi operators, we first discuss orthogonality measures \e{eq:Jac} with
exponential  weights
   \begin{equation}
   \tau (\lambda)= k_\beta e^{- |\lambda|^{\beta}}, \q \lambda \in {\Bbb R},
   \label{eq:Freud}\end{equation}
       where $\beta>0$   and $k_\beta$ are normalization constants.  Such weights were
         introduced in the paper \cite{Freud} and are known as the Freud weights.   Obviously,  the value $\beta=2$  yields the Hermite polynomials. Let $P_{n} (z)$    be the orthonormal polynomials and let $a_{n}  $, $b_{n}  $ be the recurrence coefficients defined by relations  \e{eq:RP} and \e{eq:RP1}, respectively. It was shown in \cite{LMS,Magnus} that   the off-diagonal   coefficients $a_{n}$ have asymptotics
     \begin{equation}
     a_{n}= \alpha  n^{\ell} (1+ o(1)), \q \ell=1/\beta, \q n\to\infty,
       \label{eq:FreudM}\end{equation}
       with some explicit constants $\alpha=\alpha_{\ell}$. Since  $  \tau (-\lambda)=   \tau (\lambda)$, 
  the diagonal elements $b_{n}=0$. The minimal Jacobi operator $J_{\rm min}$ with these coefficients  is essentially self-adjoint if and only if $\beta\geq 1$.

     An asymptotics of orthonormal polynomials $P_{n}(z)$ defined by measure \e{eq:Jac}, \e{eq:Freud} is given for all $z \in {\Bbb C}$ by  the Plancherel-Rotach formula generalizing the corresponding formula for the Hermite poynomials; see \cite{Nev1}  for $ \beta=4$ and  \cite{Rah, Rakh}  for   all $ \beta\geq 1$. More general results were obtained later  
 with a help of   the Riemann-Hilbert problem  method combined with the steepest descent method  (see the articles \cite{F-I-K, D-Z} and the book   \cite{Deift}).  In particular,  asymptotic  formulas for   $P_{n}(z)$ were extended in \cite{Kriech} to all  $ \beta>0$.  We  emphasize that according to the classical  Nevanlinna's  results 
 \cite{Nevan} in the case $\beta\in (0,1)$,  all self-adjoint extensions $J$ of the minimal   operator $J_{\rm min}$ have purely discrete spectra.  Note that the Plancherel-Rotach formula yields an asymptotics of   $P_{n}(z)$ as $n\to\infty$ for all  $z\in{\Bbb C}$; besides values of $z$  in this formula are not necessarily fixed.
 
 Probably, the first paper where an asymptotics of $P_{n}(z)$ was investigated under conditions  on growing  recurrence coefficients $a_{n}$ (not on the measure $d\rho(\lambda)$) is due to   Janas and Naboko \cite{Jan-Nab}.
 It was assumed in this paper that   condition  \e{eq:FreudM}  holds with $\ell\in (1/2,1)$, $b_{n} =0$ and the spectral parameter $z=\lambda\in {\Bbb R}$.  
  In  \cite{Jan-Nab}, the authors solve equations  \e{eq:RR}  successively starting from $n=0$. This yields a representation for $P_{n}(\lambda)$ in terms of a product of $n$ two-by-two matrices (the transfer matrices) expressed via $a_{n}$ and $\lambda$. Then one has to study an asymptotics of this product as $n\to\infty$  which is a non-trivial problem. The proof of the absolute continuity of the spectrum in this approach requires the   Gilbert-Pearson  subordinacy theory  \cite{GD}     adapted to Jacobi operators in  \cite{KhD}. 
     More general results of this type were obtained in the subsequent  paper \cite{Apt} by Aptekarev and Geronimo where the 
   polynomials $P_{n}(z)$  were  considered for all $z\in{\Bbb C}$. The method of \cite{Apt}    relies on a study of auxiliary Jacobi operators $J^{(N)}$ with the coefficients
$a^{(N)}_{n}=a_{n}$, $b^{(N)}_{n}=b_{n}$  for $n\leq N$ and $a^{(N)}_{n}=a_{N}$, $b^{(N)}_{n}=b_{N}$  for $n\geq N$.  Then one applies to the operators $J^{(N)}$ the results of \cite{Nev, Mate}  and studies the limit $N\to\infty$.
    Asymptotic formulas     found in  papers \cite{Jan-Nab}, \cite{Apt} are consistent with the Plancherel-Rotach formula. 
     Finally, we  note a recent paper  \cite{Sw-Tr} also devoted to Jacobi operators with increasing coefficients. 




 \subsection{Difference versus differential operators}

Our intention is to emphasize and consistently use an analogy 
  between Jacobi operators $J$ defined by formula \e{eq:ZP+1}  in the space  $L^2 ({\Bbb Z}_{+})$
 and differential  operators 
  \begin{equation}
 H=D a (x)D +b(x)
 \label{eq:HHh}\end{equation}
    (with, for example, the boundary condition $u(0)=0$)  in the space  $L^2 ({\Bbb R}_{+})$.  
   For Jacobi operators, the parameter $n\in {\Bbb Z}_{+}$ plays the role of the variable $x\in {\Bbb R}_{+}$ and the coefficients $a_{n}$, $b_{n}$ play the roles of the functions $a(x)$, $b(x)$, respectively.
   In particular,   the    operator $H_{0} = D^2$    
corresponds to the Jacobi operator $J_{0}$ with  the coefficients $a_{n}=1/2$, $b_{n}=0$. It is known as the  ``free"  discrete Schr\"odinger operator.  Of course this analogy
between difference and differential operators is very well known (see, e.g., \cite{Case1}), and we are going to use it in a systematic way.

By construction of spectral theory for the operator $H$, one studies a differential equation 
 \begin{equation}
 - (a(x) u ' (x, z) )'+ b(x) u  (x, z)= z u  (x, z), \q  x>0,
\label{eq:Schr}\end{equation}
and distinguishes its regular $\varphi (x,z)$ and Jost $f (x,z)$ solutions. The regular solution $\varphi (x,z)$ is fixed by initial conditions $\varphi (0,z)=0$,  $\varphi' (0,z)=1$ and the Jost  solution $f (x,z)$ is determined by its asymptotics as $x\to\infty$. In the simplest case when $a(x)=1$ (or  another positive constant)  and $b\in L^1 ({\Bbb R}_{+})$ the asymptotics of the Jost  solution is given by a relation $f (x,z)\sim e^{i\sqrt{z}x}$ (here $\Im \sqrt{z} \geq 0$). For every $x\geq 0$, the functions $f (x,z)$ depend analytically on $z\in {\Bbb C}\setminus [0,\infty)$ and are continuous up to the cut  along $(0,\infty)$.
This implies that the integral kernel of the   resolvent $(H-z)^{-1}$ is a continuous function of $z$ up to   the positive half-axis. This fact is known as the limiting absorption principle. It  follows  that the positive spectrum of the operator $H$ is absolutely continuous.

If $\lambda> 0$, then the solutions $f(x, \lambda+ i 0)$ and $f(x, \lambda - i 0)$ are complex conjugate to each other and are linearly independent. If $z\in {\Bbb C} \setminus [0,\infty)$, then
  the second solution $g (x,z)$ of \e{eq:Schr} can be constructed by an explicit formula stated, for example, in \S 4.1 of the book \cite{YA}. It
exponentially grows as $x\to\infty$.   Thus for all $z\in{\Bbb C}$, one has two linearly independent solutions of equation \e{eq:Schr} with explicit asymptotics as $x\to\infty$. Since the regular solution $\varphi (x,z)$ is a linear combination of $f (x,z)$ and $g (x,z)$,  this gives also an asymptotics  of $\varphi (x,z)$. The scheme described above was realized, for example, in \cite{Y-LR}.

 In more general situations,  asymptotics as $x\to\infty$ of the Jost solutions $f (x,z)$ are given (see, e.g., the book \cite{Olver}) by  the Liouville-Green formula
  \begin{equation}
 f  (x, z)\sim   {\cal G} (x,z)^{-1/2}  \exp  \Big(-\int_{x_{0}} ^x  {\cal G} (y,z) dy\Big)=: {\cal A}  (x, z).
\label{eq:Ans}\end{equation}
  Here $x_{0}$ is some fixed number and
\[
{\cal G} (x,z)= \sqrt{\frac{b(x)-z}{a(x)} } , \q\Re {\cal G} (x,z)  \geq 0.
\]
Otherwise, under quite general assumptions on the coefficients $a(x)$ and $b(x)$  the main steps of the construction of  the spectral theory for the Schr\"odinger operator $H$ remain basically similar to the particular case  $a(x)=1$,  $b\in L^1 ({\Bbb R}_{+})$.

 As far as a relation of discrete and continuous problems is concerned, we note also  the book     \cite{Atk}.

  \subsection{Structure of the paper}
  
  Our objective is to develop 
 the same scheme for Jacobi operators $J$ and the  corresponding difference equations
   \e{eq:Jy}.  We concentrate  on the case of
    recurrence coefficients satisfying the  assumptions
\begin{equation}
 a_{n} \to\infty\q\mbox{and}\q  -\frac{b_{n}}{2\sqrt{a_{n-1}a_{n}} }=:\beta_{n}\to  \beta_{\infty}  \; \mbox{where}\;  |\beta_{\infty}|\neq 1 \;\mbox{as}\;  n\to\infty.
\label{eq:Haupt}\end{equation}

 Sect.~2. is introductory. Here some general facts about Jacobi operators are discussed.

 Sect.~3 and 4 play the central role. In the first of them  we introduce and investigate   a Volterra equation which is an analytical basis of our method. Then we construct the Jost solutions $f_{n} (z)$, $n\geq -1$, and study their basic properties.  Using these results, we find   in   Sect.~5 an asymptotic behavior as $n\to\infty$ of the orthonormal polynomials $P_{n}(z)$ for $z=\lambda\in{\Bbb R}$. At the same time we obtain   spectral results for the operators $J$. In the case $|\beta_{\infty}|<1$, we check that the spectrum of the operator $J$ coincides with the whole real axis and is absolutely continuous. We also find an expression for the spectral measure of the operator $J$ in terms of boundary values $|f_{-1} (\lambda\pm i0)|$.  An asymptotic behavior of the  polynomials $P_{n}(z)$ for $\Im z\neq 0$ is studied in  Sect.~6. 
 Our scheme is essentially the same in the cases $|\beta_{\infty}|<1$ and $|\beta_{\infty}| > 1$, but asymptotic formulas and spectral results obtained in these two cases are drastically different.
 
 In Sect.~3-6, we assume that $a_{n}\to\infty$ but not too rapidly so that
  the condition (introduced in the book \cite{Carleman} and known as the Carleman condition)
\begin{equation}
\sum_{n=0}^\infty a_{n}^{-1}=\infty
\label{eq:Carl}\end{equation}
is satisfied. The  singular case where this condition is violated is studied
   in  Sect.~7. Astonishingly,  in this case  the asymptotic formulas for the Jost solutions and orthonormal polynomials are particularly simple. 
   
    Sect.~8 and 9 are devoted to another limit case of the results of Sect.~3-6 where 
   \begin{equation}
 a_{n} \to a_{\infty} >0 \q \mbox{and}\q b_{n} \to b_{\infty}  \q\mbox{as}\q  n\to\infty. 
 \label{eq:stabi}\end{equation} 
 Here our presentation  is rather sketchy   since   this
  case was  already treated by the methods of the present paper in    \cite{Y/LD, JLR}.  Note that in terms of the analogy with differential operators conditions \e{eq:Tr} and 
  \e{eq:LR} correspond, respectively,   to short-range and long-range perturbations of the operator $D^{2}$.
   
   
 
    We mention that many important topics in the theory of orthogonal polynomials and Jacobi  operators such as 
      the critical case $|\beta_{\infty}|=1$ and various types of oscillating coefficients  
     are completely out of the scope of this text.

 To stress an analogy between differential and difference operators, we often  use the  ``continuous" terminology (Volterra  integral equations, integration by parts, etc.) for sequences labelled by the discrete variable $n$. Below $C$, sometimes with indices,  and $c$ are different positive constants whose precise values are of no importance; $I$ is the identity operator in various spaces.  We usually do not care about precise estimates of various remainders writing $o(1)$.


 
 
   
 \section{Scheme of the approach}

   Our plan here is the following. After introducing necessary notations in Sect.~2.1, we construct  the resolvents of Jacobi operators in Sect.~2.2 and discuss exponentially growing solutions of the Jacobi equations in Sect.~2.3.  In Sect.~2.4 we recall the uniformization $z\mapsto \z$ playing the role of the relation $z\mapsto \z=\sqrt{z}$ for differential operators.  Crucial steps of our approach are described  in Sect.~2.5.  Then we briefly state some particular cases of our main results  in Sect.~2.6.

 \subsection{Standard  relations}

Let us consider  the difference equation \e{eq:Jy}. Note that the values of $u_{N-1}$ and $u_{N }$ for some $N\in{\Bbb Z}_{+}$ determine the whole sequence $u_{n}$ satisfying    \e{eq:Jy}. 
In particular, this is true if $u_{-1} $ and $u_{0}$ are given.  We often start a construction of solutions $u_{n}$ of equation \e{eq:Jy} with large $n$. Then they are   extended to all $n\geq -1$ by formula  \e{eq:Jy}.

 Let $u= ( u_{n} )_{n=-1}^\infty$ and $v= (v_{n} )_{n=-1}^\infty$ be two solutions of the  Jacobi equation \e{eq:Jy}. A direct calculation shows that their Wronskian
  \begin{equation}
\{ u,v \}: = a_{n}  (u_{n}  v_{n+1}-u_{n+1}  v_{n})
\label{eq:Wr}\end{equation}
does not depend on $n=-1,0, 1,\ldots$. In particular, for $n=-1$ and $n=0$, we have
 \[
\{u,v \} = a_{-1} (u_{-1}  v_{0}- u_{0}  v_{-1}) \q {\rm and} \q \{ u,v\} = a_{0}  (u_{0}  v_{ 1}-u_{ 1}  v_{0});
\]
here and below $a_{-1}$ is an arbitrary  fixed positive number.
Clearly, the Wronskian $\{ u,v \}=0$ if and only if the solutions $u$ and $v$ are proportional.

Let  $J_{\rm min}$ be the operator defined in the space $\ell^{2}({\Bbb Z}_{+})$ by formula \e{eq:ZP+1}  on the set ${\cal D}$ of sequences $u=(u_{n})$ such that $u_{n}=0$ for sufficiently large $n$.
In the limit circle case, all solutions of the Jacobi equation \e{eq:Jy}  are in $\ell^{2}({\Bbb Z}_{+})$ so that necessarily
\[
 a_{n}^{-1}= \{ u,v \}^{-1} (u_{n}v_{n+1}- u_{n+1}v_{n} ) \in \ell^{1}({\Bbb Z}_{+}).
\]
It follows that   the Carleman condition
\e{eq:Carl} 
is  sufficient  for   the limit point case  and hence for the essential self-adjointness of the operator $J_{\rm min}$.    




 It is convenient to introduce a notation
 \begin{equation}
u_{n}'= u_{n+1}  - u_{n}
\label{eq:diff}\end{equation}
for the ``derivative" of a sequence $u_{n}$. Note a formula
\begin{equation}
(u_{n} v_{n})'= u_{n}' v_{n+1}  + u_{n} v_{n}'.
\label{eq:diff1}\end{equation}
In particular, this yields the Abel summation formula (``integration-by-parts"):
 \begin{equation}
\sum_{n=N }^ { M} u_{n}  v_{n}' = u_{M}  v_{M+1} - u_{N -1}  v_{N}  -\sum_{n=N } ^{M} u_{n-1}'  v_{n};
\label{eq:Abel}\end{equation}
here $M\geq N\geq 0$ are arbitrary, but we have to set $u_{-1}=0$ so that $u_{-1}'=u_{0}$.



 
  \subsection{Resolvent}
  
 In this subsection we only suppose that the minimal Jacobi operator $J_{\rm min}$ is essentially self-adjoint.  Then $\clos J_{\rm min}=J_{\rm max}=:J$ and  equation \e{eq:Jy}   has (see, e.g., \S 3 of Chapter~1 in the book \cite{AKH}) a unique, up to a constant factor, non-trivial solution $f_{n} (z)$  such that  \begin{equation}
     f_{n}( z )\in \ell^2 ({\Bbb Z}_{+}), \q \Im z\neq 0. 
   \label{eq:asqrf1+}\end{equation}  
     Let us introduce the Wronskian of the solutions $P(z)= ( P_{n} (z) )$ and $f(z)= ( f_{n} (z) )$:
  \begin{equation}
\Omega(z) : = \{ P (z), f (z) \} =  a_{-1} (P_{-1}(z) f_{0} (z)- P_0 (z){ f}_{-1} (z))=- a_{-1} f _{-1} (z) .
\label{eq:J-W}\end{equation} 
Observe that $\Omega(z)\neq 0$ if $\Im z \neq 0$.  Indeed, otherwise $P(z)\in \ell^2 ({\Bbb Z}_{+})$ so that, in view of  equation \e{eq:RR},  $z$ is an eigenvalue of the operator $J$. This is impossible since $J$ is self-adjoint.
 
 Our goal is to construct the resolvent $R(z)= (J-z I)^{-1}$ of the operator $J$ for $\Im z\neq 0$.    The following statement is very close to the corresponding result for differential operators.

 \begin{proposition}\label{res}
 For all $h= (h_{n})\in \ell^{2} ({\Bbb Z}_{+})$, we have
  \begin{equation}
(R (z)h)_{n} = \Omega(z)^{-1} \Big( 
f_{n} (z) \sum_{m=0}^{n} P _{m}(z)h_{m}+ P_{n} (z) \sum_{m=n+1}^{\infty}f_{m}(z) h_{m}\Big),\q \Im z \neq 0. 
\label{eq:RRes}\end{equation}
 \end{proposition} 
 
  \begin{pf} 
  Denote the right-hand side of \e{eq:RRes} by $( {\sf R}(z) h)_{n}$.
    We have to check that 
    \begin{equation}
    ({\cal J}-zI) {\sf R} (z)h=h, \q \Im z \neq 0, 
    \label{eq:RRm}\end{equation}
at least for all sequences $h\in {\cal D}$.  Let us set
   \begin{equation}
A_{n}(z) =\sum_{m=0}^n  P_{m} (z) h_{m},\q    B_{n}(z) =\sum_{m=n+1}^\infty  f_{m} (z) h_{m}.
\label{eq:RR2}\end{equation}
  Then 
    \begin{equation}
  \Omega(z)( {\sf R} (z)h)_{n} =  f_{n}(z) A_{n}(z)+   P_{n}(z) B_{n}(z).
\label{eq:RR11}\end{equation}
  Note that  $A_{n}(z)$ does not depend on $n$ and    $B_{n}(z)=0$  for sufficiently large $n$.
In view of  \e{eq:asqrf1+}, we have
 ${\sf R} (z)h\in \ell^2 ({\Bbb Z}_{+})$.

It follows from definition \e{eq:ZP+1}  of  the Jacobi matrix ${\cal J} $   that
   \begin{multline}
  \Omega ( ({\cal J} -z I) {\sf R} (z) h)_{n} =  a_{n-1}\big(f_{n-1}A_{n-1}+P_{n-1} B_{n-1} \big)
  \\
  + (b_{n}-z) \big(f_{n}A_{n}+ P_{n} B_{n}  \big)+ a_{n}\big(f_{n+1}A_{n+1} +P_{n+1} B_{n+1}    \big)
\label{eq:RR3}
\end{multline}
for $n\geq 1$ and 
   \begin{equation}
  \Omega ( ({\cal J} -z I) {\sf R} (z) h)_0 =    (b_{0}-z) \big(f_{0}A_{0}+ P_{0} B_{0}  \big)+ a_{0}\big(f_{1}A_{1} +P_{1} B_{1}    \big).
\label{eq:RR3+}
\end{equation}
According to \e{eq:RR2} we have
\[
f_{n-1}A_{n-1}+P_{n-1} B_{n-1}=f_{n-1}(A_{n}-P_{n} h_{n})+  P_{n-1}(B_{n}+f_{n} h_{n}), \q n\geq 1,
\]
and
   \begin{equation}
f_{n+1} A_{n+1} +P_{n+1} B_{n+1} =
f_{n+1}A_{n}  +  P_{n+1}B_{n} .
\label{eq:RR3+1}
\end{equation}
Let us substitute these expressions into the right-hand side of \e{eq:RR3} and observe that
the coefficients at $A_{n}$ and $B_{n}$ equal zero by virtue of equation \e{eq:Jy} for $f_{n}$ and $P_{n}$, respectively. It  follows that
   \begin{equation}
 ( ({\cal J}-z I) {\sf R} (z) h)_{n} = \Omega^{-1}  a_{n-1}  (- P_{n}  f_{n-1} +    f_{n} P_{n-1}) h_{n}  =  h_{n} , \q n\geq 1.
\label{eq:RR3+2}
\end{equation}
Next, we consider the right-hand side
of \e{eq:RR3+}. According to  equality \e{eq:RR3+1} for $n=0$ it equals
\[
 (b_{0}-z) \big(f_{0}A_{0}+ P_{0} B_{0}  \big)+ a_{0}\big(f_{1}A_{0} +P_{1} B_{0}    \big)
 =\big( (b_{0}-z) f_{0}+ a_{0} f_{1}\big) A_{0}+ 
 \big( (b_{0}-z)  + a_{0} P_{1}\big) B_{0}.
 \]
 The first term on the right is $-a_{-1} f_{-1} A_{0}= \Omega h_{0}$ by   equation \e{eq:Jy} for $f_{n}$, equality $A_{0}=h_{0}$ and definition  \e{eq:J-W} of the Wronskian $\Omega$. The second term is zero by  equation \e{eq:Jy} for $P_{n}$.  Thus it  follows from \e{eq:RR3+} that 
$ ( ({\cal J}-z) {\sf R} (z) h)_0 =   h_{0}$.  Together with \e{eq:RR3+2}, this implies 
equality  \e{eq:RRm}.  In particular, we see that ${\sf R}  (z) h\in {\cal D} (J_{\rm max}) =  {\cal D} (J)$ so that 
$  (J-zI) {\sf R} (z)h=h$ 
whence  ${\sf R}(z)=(J-zI)^{-1}$.  
    \end{pf}
    
      \subsection{Linearly independent solutions of the Jacobi equation} 
      
      Here we suppose that a solution of the Jacobi equation \e{eq:Jy} is given and describe a general procedure which allows one to construct another solution, linearly independent with the first one. This procedure plays  a crucial role in our study of an asymptotic behavior of the orthonormal polynomials $P_{n}(z)$ as $n\to\infty$  for complex $z$.  Recall that
  the Wronskian of two solutions of equation  \e{eq:Jy} is defined by relation \e{eq:Wr}.

 \begin{theorem}\label{GE}
 Let $f (z)= ( f_{n}(z) )$ be an arbitrary solution of the Jacobi equation \e{eq:Jy} such that $f_{n}(z) \neq 0$ for sufficiently large $n$, say $n\geq  N_{0}$. Put
   \begin{equation} 
 F_{n} (z)=\sum_{m=N_{0} + 1}^n (a_{m-1} f_{m-1}(z) f_{m}(z))^{-1},\q n \geq N_{0} + 1,
\label{eq:GE}\end{equation}
and
  \begin{equation} 
g_{n}(z) =f_{n}(z) F_{n}(z), \q n \geq N_{0} + 1.
\label{eq:GEx}\end{equation}
Then  the sequence $g (z)= ( g_{n}(z) )$ 
 satisfies equation \e{eq:Jy} and
   the Wronskian 
     \[
\{ f(z),g(z)\}=1
\]
    so that
   the solutions $f(z)$ and $g (z)$ are linearly independent.  
    \end{theorem}

\begin{pf}
First, we check equation \e{eq:Jy} for $g_{n} (z)$. According to definition \e{eq:GEx}, we have
 \begin{multline*}
 a_{n-1} g_{n-1}  + (b_{n} -z) g_{n}   +a_{n} g_{n+1} =
a_{n-1} f_{n-1} F_{n-1}+ (b_{n} -z) f_{n} F_{n} +a_{n} f_{n+1}F_{n+1}
\\
=\big(a_{n-1} f_{n-1}  + (b_{n} -z) f_{n}  +a_{n} f_{n+1}\big) F_{n}
+ a_{n-1} f_{n-1} ( F_{n-1}- F_{n}) + a_{n} f_{n+1}(F_{n+1}-F_{n}) .
\end{multline*}
The first term here is zero because    equation \e{eq:Jy} is true for the  sequence $f_{n}$. Since, by  \e{eq:GE},
  \begin{equation} 
F_{n+1}=F_n+ (a_{n} f_{n} f_{n+1})^{-1},
 \label{eq:GExx}\end{equation}
the second and third terms equal $-f_{n}^{-1}$ and $f_{n}^{-1}$, respectively, which proves equation  \e{eq:Jy} for sequence \e{eq:GEx}.

It also follows from definition \e{eq:GEx} and relation  \e{eq:GExx}
that the Wronskian \e{eq:Wr}   equals
  \[
\{ f(z),g(z)\}= a_{n}f_{n}(z)f_{n+1}(z)(F_{n+1}(z)- F_{n}(z))=1,
\]
 whence the solutions $f(z)$ and $g (z)$ are linearly independent.
\end{pf}

As usual,  a sequence $g_{n}(z)$ constructed for large $n$ is extended to all $n\geq -1$ as a solution of equation \e{eq:Jy}. 

\subsection{Uniformization}

We  fix the branch of the analytic function $\sqrt{z^2 -1}$ of 
\[
z\in {\Bbb C}\setminus [-1,1]=:\Pi_{0} 
 \]
  by the condition $\sqrt{z^2 -1}>0$ for $z>1$. Obviously, this function is continuous up to the cut along $[-1,1]$,    it equals $\pm i\sqrt{1-\lambda^2}$ for $z=\lambda\pm i0$ where $\lambda\in (-1,1)$ and $\sqrt{z^2 -1}< 0$ for $z< -1$. Put
\begin{equation}
 \z (z) =z-\sqrt{z^2 -1}= (z + \sqrt{z^2 -1})^{-1} .
 \label{eq:ome}\end{equation}
 The mapping
$\z$ of $ \Pi_{0} $ onto  the unit disc $  {\Bbb D} $ is one-to-one  and holomorphic.
Note also that  $\z(-z)=-\z(z)$ and
$ \z (\bar{z}) = \ov {\z (z)}$.  The function $\z(z)$ maps the half-planes $  {\Bbb C}^{(\pm)}= \{z\in   {\Bbb C} :  \pm \Im z> 0 \}$ onto the half-discs $ {\Bbb D}^{(\mp)}= {\Bbb D}\cap   {\Bbb C}^{(\mp)}$; thus, 
$ \z : {\Bbb C}^{(\pm)}\mapsto {\Bbb D}^{(\mp)}$.    For $\lambda>1$, we have  $\z( \lambda)\in (0,1)$ and $\z( 1)=  1$, $\z( +\infty)=  0$.  Similarly, $\z (\lambda )\in (-1,0)$ for $\lambda< -1$,   $\z( -1)= - 1$ and $\z( -\infty)=  0$.
    For $\lambda\in [-1,1]$, it is common to set $\lambda=\cos \theta$ with $\theta\in [0,\pi ]$. Then 
$\z (\lambda\pm i0)=e^{\mp i \theta}$. Note also that 
 \[
 \z (z) +   \z (z)^{-1}=2 z  
  \]
so that $\z(z)$ is the inverse Zhukovsky function. 
 Function \e{eq:ome} plays the role of the function $\z (z) =\sqrt{z}$ in the theory of the Schr\"odinger operator.
  In all estimates below 
 the values of   $|z|$ are   bounded. Then the values of    function \e{eq:ome} are separated from $0$.

  Let us introduce a notation
   \begin{equation}
\alpha_{n} =\frac{1}{2\sqrt{a_{n-1}a_{n}} } ,  \q \beta_{n} =-\frac{b_{n}}{2\sqrt{a_{n-1}a_{n}} }
\label{eq:aabb}\end{equation}
 and set
\begin{equation}
 z_{n}  = \alpha_{n} z + \beta_{n},\q   \z_{n} = \z (z_{n}).
\label{eq:aabb1}\end{equation}
   
According to assumption \e{eq:Haupt}  we have
\begin{equation}
 \lim_{n\to\infty}z_{n}  =   \beta_{\infty}
\label{eq:aabb2}\end{equation}
whence 
\begin{equation}
 \lim_{n\to\infty}\z (z_{n} )   =   \z (\beta_{\infty}\pm i0) \q \mbox{if}\q \pm \Im z\geq 0.
\label{eq:aab}\end{equation}
Note that
\begin{equation}
 \z (\beta_{\infty} + i0)=  \z (\beta_{\infty} - i0)=   \beta_{\infty}-  \sqrt{\beta_{\infty}^2-1}=: \z_{\infty}\in (-1,1)
\label{eq:aabb3}\end{equation}
if $|\beta_{\infty}| >1$ and
\begin{equation}
 \z (\beta_{\infty}\pm i0)  =   \beta_{\infty}\mp  i \sqrt{1-\beta_{\infty}^2} =: \z_{\infty}^{(\pm)}\in {\Bbb T}\cap{\Bbb C}^{(\pm)}, \q \pm \Im z\geq 0,
\label{eq:aabb4}\end{equation}
if $|\beta_{\infty}| <1$. 
  Relations \e{eq:aabb2} and hence \e{eq:aabb3}, \e{eq:aabb4} are uniform in $z$ from compact subsets of $\Bbb C$.

 We set $\Pi={\Bbb C}\setminus{\Bbb R}$ and denote by $\clos\Pi$ its closure. Thus,  $\clos\Pi$ is the complex plane with the cut along $\Bbb R$.

  \subsection{Main steps} 

  
  Our study of an asymptotic behavior of  the orthonormal polynomials $P_{n} (z)$  as $n\to\infty$ consists of the following steps. We describe them here for Jacobi operators $J$ with coefficients  satisfying    assumptions \e{eq:Haupt}  with $|\beta_{\infty}|<1$ and the Carleman condition \e{eq:Carl}.   For  Schr\"odinger operators with short-range coefficients, the   scheme  used here goes back to  the paper \cite{Jost} by R.~Jost.

A. 
First, we forget about the orthogonal polynomials  
and distinguish solutions (the Jost solutions) $f_{n}  (z)$ of the  difference  equation \e{eq:Jy} 
 by their asymptotics as $n\to\infty$. This requires a construction of an Ansatz ${\cal A}_{n} (z)$   for   the Jost solutions. Actually, ${\cal A}_{n} (z)$   turns out to be the leading term of $f_{n}  (z)$ for $n\to\infty$.

   B. 
    We define ${\cal A}_{n} (z)$ by a formula 
       \begin{equation}
{\cal A}_{n}(z) = a_{n}^{-1/2} \z (z_{N_{0}})\z (z_{N_{0}+1})\ldots \z(z_{n-1}),\q n\geq N_{0}\q z\in \clos\Pi  , 
\label{eq:AnsDD}\end{equation}
with $z_{n}$ defined by \e{eq:aabb1}.
Using \e{eq:aabb2} we choose the number  $N_{0}= N_{0} (z)$    in such a way that $z_{n} $   for $n\geq N_{0} $ are separated from the points $+1$ and $-1$.  Clearly,   $N_{0}$ can be chosen common for all $z$ from a given compact subset of $\clos\Pi$.
It turns out that the relative remainder
\begin{equation}
 r_{n} (z) : =  ( \sqrt{a_{n-1}a_{n}}{\cal A}_{n} (z))^{-1} \Big(a_{n-1}{\cal A}_{n-1} (z) + (b_{n}-z){\cal A}_{n} (z) + a_{n}{\cal A}_{n+1} (z)\Big)
\label{eq:Gry}\end{equation}
belongs to $ \ell^1 $ for $n\to\infty$.  Formula \e{eq:AnsDD} plays the role of the Liouville-Green Ansatz  \e{eq:Ans} for differential operators.
According to  \e{eq:aab}, \e{eq:aabb4}  asymptotics of the products
\begin{equation}
 q_{n} (z)  = \z_{N_{0}} \z_{N_{1}}\cdots  \z_{n-1}  ,\q   \z_{n}= \z (z_{n}),
\label{eq:re1}\end{equation}
contain an oscillating factor $(\z_{\infty}^{(\pm)})^n$. Moreover, it can be deduced from the Carleman condition \e{eq:Carl} that $ |q_{n} (z)|  \to 0$ as     $n\to\infty$  if $\Im z\neq 0$ while $ |q_{n} (z)|  \to c\neq 0$    if $\Im z= 0$.
 


C.
A multiplicative change of variables
 \begin{equation}
   f_{n} (z)= {\cal A}_{n}  (z) u_{n} (z ) 
      \label{eq:Jost}\end{equation} 
    reduces    difference  equation \e{eq:Jy} for  $  f_{n} (z)$  to    a
Volterra ``integral" equation for the sequence $u_{n} (z)$.    This Volterra equation can be standardly solved    by iterations which  allows us to prove  the existence of its solution $ u_{n} (z) $    such that
\begin{equation}
\lim_{n\to\infty} u_{n} (z)=   1.
\label{eq:A12a}\end{equation} 
 Then the Jost solution $ f_{n} (z)$ of equation  \e{eq:Jy}  is  defined by formula
    \e{eq:Jost}. It follows from \e{eq:A12a} that  
     \begin{equation}
   f_{n} (z)= {\cal A}_{n}  (z) (1+ o(1))
      \label{eq:JostAs}\end{equation} 
as $n\to\infty$.
        The functions $ u_{n} (z)$ and therefore $ f_{n} (z)$  turn out to be  analytic in $z$  in the complex plane with the cut along the real axis and are continuous up to the cut.
          
 We emphasize that the sequences  $ {\cal A}_{n} (z)$  and hence $ f_{n} (z)$ are defined up to factors depending, possibly, on $z$. These factors can be chosen at our convenience.

   
   D. 
   If $z=\lambda\in {\Bbb R}$, then equation \e{eq:Jy} has two linearly independent solutions $ f_{n} (\lambda+ i0)$  
   and 
   \[
    f_{n} (\lambda - i0)=\ov{    f_{n} (\lambda + i0)}.
    \]
     Their Wronskian equals
        \[
\{ f (\lambda+i0),  f(\lambda-i0)\} 
= 2 i   \sqrt{1- \beta^2_{\infty}} \neq 0.
 \]
  Therefore the   polynomials $P_{n}  (\lambda)$ are linear combinations of $f(\lambda+i0)$ and $  f (\lambda-i0)$:
  \begin{equation}
P_{n} (\lambda)=\frac{  \Omega (\lambda- i0)  f_{n} (\lambda+i0) - \Omega(\lambda + i0)  f_{n} (\lambda- i0)  }{ 2 i   \sqrt{1- \beta^2_{\infty}}},   \q n=0,1,2, \ldots, \q \lambda\in \Bbb R,
\label{eq:HH4L}\end{equation}
  where    $\Omega (z)$ is  Wronskian \e{eq:J-W}. 
 Note that $\Omega (\lambda\pm i0)\neq 0$ for all $\lambda\in \Bbb R$ so that formulas \e{eq:JostAs} and \e{eq:HH4L} yield an  
   asymptotics of $P_{n} (\lambda)$ as $n\to\infty$.

    If $\Im z\neq 0$,   a solution $ g_{n} (z ) $ of equation \e{eq:Jy} linearly independent with $ f_{n} (z ) $ can be constructed by  explicit formulas
    \e{eq:GE}, \e{eq:GEx}. The sequence $g_{n} (z)$ rapidly  grows     as $n\to\infty$ and a   limit
    \begin{equation} 
  \lim_{n\to\infty } \big( {\cal A}_{n} (z)g_{n} (  z )\big)= \mp \frac{ i}{2\sqrt{1-\beta_{\infty}^2}} , \q \pm \Im z>0,
\label{eq:A2P3}\end{equation}
exists. According to Theorem~\ref{GE}  we have
  \begin{equation}
P_{n} (z)= \omega (z)f_{n} (z) - \Omega(z) g_{n} (z)
\label{eq:Pf}\end{equation}
where $\omega  (z) = \{ P(z),g(z)\}$.
 Therefore the asymptotics of $P_{n}(z)$ is given by the relation
     \begin{equation}
\lim_{n\to\infty}  \big({\cal A}_{n}  (z) P_{n}(z)\big)= \pm i\frac{\{ P(z), f(z)\}}{2 \sqrt{1-\beta_{\infty}^2} } , \q \pm \Im z>0.
\label{eq:HSq5}\end{equation}

  E.
  Our results on the Jost solutions $( f_{n}  (z)) $ directly imply    that  the   spectrum of the Jacobi operator $J$ is absolutely continuous and covers the whole axis   $\Bbb R$. Its spectral measure is given by the formula
     \begin{equation}
  d\rho(\lambda)=  \pi^{-1}\sqrt{1- \beta^2_{\infty}} \, | \Omega (\lambda+i0) |^{-2} d\lambda.
  \label{eq:Gg}\end{equation}
  At  the same time, we obtain the limiting absorption principle for the operator $J$ stating that matrix elements $\la R(z ) u,v \ra$, $\Im z \neq 0$,  of its resolvent  are continuous functions of $z$ up to the cut along the real axis for all $u,v\in{\cal D}$.

Let now assumptions \e{eq:Haupt}  hold with $|\beta_{\infty}|>1$.  Then the scheme described above remains the same; in particular, the Ansatz is again defined by equality  \e{eq:AnsDD}. Nevertheless, asymptotic and spectral results are quite different  from the case  $|\beta_{\infty}|<1$. Thus, the functions $ {\cal A} _{n} (z)$ and hence $f_{n} (z)$ exponentially
tend to zero as $n\to\infty$   for all $z\in{\Bbb C}$.  On the contrary, the solutions $g_{n} (z)$ and the orthonormal polynomials  $P_{n} (z)$
exponentially tend to infinity as $n\to\infty$. The spectra of the corresponding Jacobi operators $J$ are  discrete.

 To a large extent, a construction of an Ansatz for difference equations is similar   to their formal solutions. Such a procedure was suggested by Birkhoff in \cite{Birkh} and substantially developed in \cite{W-L}; see the book  \cite{El}, for a detailed presentation.

We emphasize that the approach of this paper is essentially different from
those of \cite{Jan-Nab} and \cite{Apt}   discussed in  Sect.~1.4.

      
        \subsection{An overview of the main results}
        
        Here we state our results under some simplifying assumptions; see Sect.~4.3 for  more details. We now suppose that the Carleman condition \e{eq:Carl} is satisfied but the coefficients $a_{n}$ do not tend to infinity too slowly.
          More precisely, we require  that
        \begin{equation}
\sum_{n=0}^\infty   a_{n}^{-3}  (1+ |b_{n}|)<\infty.
\label{eq:D}\end{equation} 
Conditions \e{eq:Carl} and \e{eq:D} admit a growth of the off-diagonal coefficients $a_{n}$ as $n^p$ for $p\in (1/2, 1]$ and even for $p\in (1/3, 1]$ if $b_{n}=0$.

Recall that the numbers $\alpha_{n}$, $\beta_{n}$ are defined by equalities   \e{eq:aabb} and $\beta_{\infty}$ is the limit of $\beta_{n}$ as $n\to\infty$. Suppose first that $|\beta_{\infty}| <1$ and set
 \begin{equation}
\phi_{n}= \sum_{m=N_{0}+1}^{n-1} \arccos \beta_{m} \q \mbox{and}\q \psi_{n}= \sum_{m= N_{0}+1}^{n-1} \ \frac{ {\alpha}_m }{\sqrt{|1-\beta_m^{2}|}};
\label{eq:psi}\end{equation}
here $N_{0}$ is so large that $|\beta_{n}|< 1$ for all $n > N_{0}$.
   Under   assumption  \e{eq:D} expression  \e{eq:AnsDD} for the Ansatz ${\cal A}_{n} (z)$ can be simplified which yields asymptotics of the Jost solutions $f_{n} (z)$ of the Jacobi equation \e{eq:Jy} in the following explicit form
  \begin{equation}
f_{n} (z)=a_{n}^{-1/2}e^{\mp i \phi_{n} \pm iz \psi_{n} }(1+ o(1)), \q \pm \Im z \geq 0, \q n\to \infty.
\label{eq:D1}\end{equation} 
 As an example, we note the recurrence coefficients $a_{n}=\nu (n+1)^p$ where $\nu >0$, $p\in (1/3, 1]$ and $b_{n}=0$.  Then  formula  \e{eq:D1} reads as
\[
f_{n}(  z )= \nu^{-1/2}  n ^{-p/2} i^{\mp n}  e^{ \pm i z \gamma n^{1-p}} (1+ o(1)), \q \gamma= ( 2 (1-p)\nu)^{-1}.
\]
In the case $p=1/2$ this is of course consistent with formulas for the Hermite polynomals  when $a_{n }=\sqrt{(n+1)/2}$ and $b_{n}= 0$.

In the general case where $a_{n}$, $b _{n}$ satisfy \e{eq:Haupt}  we have
\[
\phi_{n}= n\arccos\beta_{\infty}+ o(n) \q \mbox{as}\q n\to\infty
\]
and
 $\psi_{n}=O (n^{2/3})$ according to \e{eq:D}, whence the phases $\psi_{n}$ are negligible compared to $\phi_{n}$. However it follows from \e{eq:D1}  that
  \[
|f_{n} (z) | =a_{n}^{-1/2}e^{- |\Im z|  \psi_{n} }(1+ o(1)),  \q n\to \infty,
\]
and hence the behavior of $| f_{n} (z)|$ as $n\to\infty$
  for $\Im z\neq 0$ is determined by   $\psi_{n}$. 
     It is easy to show (see Lemma~\ref{asq2} below) that $f_{n} (z)\in \ell^2 ({\Bbb Z}_{+})$ for $\Im z \neq 0$. For each $n\geq -1$, the functions $f_{n} (z)$ are  analytic in the half-planes
 $\pm \Im z>0$ and are continuous up to the real line.

Consider now the orthonormal    polynomials $P_{n} (z)$. Suppose first that $z =\lambda\in {\Bbb R}$.
   It  follows from formulas \e{eq:HH4L} and \e{eq:D1}   that, as $ n\to\infty$,   the  polynomials $P_{n} (\lambda)$ have asymptotics 
   \begin{equation}
 P_{n} (\lambda)= - a_{n}^{-1/2} \Big(   | \Omega ( \lambda+i0)| (1- \beta^2_{\infty})^{-1/2}     \sin (\phi_{n}  -\lambda\psi_{n} +\arg \Omega(\lambda+i 0) ) + o(1) \Big)  
\label{eq:Sz}\end{equation}
  where   $\Omega (z)$ is  Wronskian \e{eq:J-W}.
    Let us comment on asymptotic coefficients in formula \e{eq:Sz}. For simplicity, we assume that $a_{n}=\nu (n+1) ^p$ where $p\in (1/2,1)$, $\nu >0$, and neglect remainders.  It follows from \e{eq:psi} that
  $  \psi_{n}= \upsilon n^s (1+ o(1)) $ with $ s=1-p$ and $ \upsilon= \big( 2\nu s \sqrt{1 - \beta^2_{\infty} } \big)^{-1}$.
  The amplitude in  \e{eq:Sz} equals $\kappa(\lambda) n^{-r}$ where $\kappa(\lambda)= \nu^{-1/2} | \Omega ( \lambda+i0)| (1- \beta^2_{\infty})^{-1/2} $  and $r=p/2$. Note that $2r+s=1$ which is one of the universal relations observed in \cite{univ}.  Another relation of \cite{univ} links the amplitude factor with the spectral weight $\tau (\lambda)$: 
  \begin{equation}
  \pi \tau (\lambda) \kappa^2(\lambda) = 2s\upsilon.
  \label{eq:Sza}\end{equation}
  Substituting here expressions for $\upsilon$ and $\kappa(\lambda)$, we see that in our case relation \e{eq:Sza} is equivalent to formula  \e{eq:Gg}.

  If $\Im z \neq 0$, then the solution  $g_{n}(z)$ of the Jacobi equation \e{eq:Jy} linear independent with $f_{n}(z)$
    can be constructed by    formulas \e{eq:GE}, \e{eq:GEx}. According to \e{eq:A2P3}
     it has the asymptotics
  \[
   g_{n}(z)=\frac{\mp i}{2 \sqrt{1 - \beta^2_{\infty} } }  \frac{1}{ \sqrt{a_{n}}  }e^{\pm i \phi_{n} \mp  iz \psi_{n}}   (1+o(1)), \q \pm \Im z>0, 
   \q n\to\infty,
\]
so that $ | g_{n}(z) |$ grows   faster than any power of $n$ as $n\to\infty$.
The  asymptotics   of the orthonormal polynomials is given by formula \e{eq:HSq5}:
  \begin{equation}
    P_{n}(z)=\frac{\pm i \Omega(z)}{2 \sqrt{1 - \beta^2_{\infty} } } 
 \frac{1}{ \sqrt{a_{n}}  }  e^{\pm i \phi_{n} \mp iz \psi_{n}}  (1+o(1)), \q \pm \Im z>0,    \q n\to\infty.
\label{eq:Gas}\end{equation}

  Note that formulas \e{eq:Sz} and \e{eq:Gas} are consistent with the classical asymptotic relation \e{eq:H2} for the Hermite polynomials  (see, e.g., Theorems~8.22.6 and 8.22.7 in the G.~Szeg\H{o}'s book \cite{Sz}). 
   Asymptotics \e{eq:Sz} was obtained earlier in \cite{Jan-Nab}, but we are unaware  of  papers where \e{eq:Gas} was deduced from assumptions on the Jacobi coefficients (not from properties of the corresponding spectral measure).
   
   Under the Carleman condition  \e{eq:Carl}  the Jacobi operator $J=\clos J_{\rm min}$ is self-adjoint, and
its resolvent $ (J-zI)^{-1}$   is given by the general formula \e{eq:RRes}. 
 Since the Jost solutions $f_{n} (z)$ depend continuously on  $z$ up to the real axis and 
 $\Omega (\lambda\pm i0)\neq 0$ for   $\lambda\in \Bbb R$, the spectrum of  the Jacobi operator $J$ is absolutely continuous, covers the whole real line,
 and its spectral measure can be constructed by relation  \e{eq:Gg}.

      The main difference between the cases $|\beta_{\infty}|> 1$   and $|\beta_{\infty}|< 1$ is that according to \e{eq:aabb3}, \e{eq:aabb4} $|\z (\beta_{\infty})|< 1$ for  $|\beta_{\infty}|> 1$   while $|\z (\beta_{\infty}\pm i0)|= 1$ for $|\beta_{\infty}|< 1$.
      Technically these cases are rather similar although many estimates are simpler for $|\beta_{\infty}|> 1$. However the asymptotic behavior of orthonormal polynomials and spectral properties of the Jacobi operators are quite different in these cases.  If $|\beta_{\infty}|> 1$, then for all $z\in{\Bbb C}$    
the Jost solution $f_{n} (z)$   of equation \e{eq:Jy}  is  distinguished by the asymptotics 
 \begin{equation}
f_{n} (z)=a_{n}^{-1/2} (\sgn\beta_\infty)^n e^{-  \varphi_{n} - z \psi_{n}}  (1+ o(1)), \q   n\to \infty,
\label{eq:D1+}\end{equation} 
where    
 \begin{equation}
\varphi_{n}= \sum_{m=N_{0}}^{n-1} \arccosh | \beta_{n}|,\q n\geq N_{0}.
\label{eq:psi+}\end{equation}
and  $\psi_{n}$ are, as before,  defined by formula \e{eq:psi}. Thus, $f_{n} (z)\to 0$ exponentially as $n\to\infty$ for all $z\in{\Bbb C}$.
Now the functions   $f_{n} (z)$ are analytic in the whole complex plane.  
As usual, the second solution $g_{n}(z)$ is constructed by  equalities \e{eq:GE}, \e{eq:GEx}.  Its asymptotics as $n\to\infty$ can be deduced from formula \e{eq:D1+}. In view of relation \e{eq:Pf} this yields    an asymptotic formula for the orthonormal polynomials $P_{n} (z)$: 
  \[
  P_{n}(z)=-      \Omega(z)  \frac{ (\sgn \beta_{\infty})^{n +1} }{2\sqrt{\beta^2_{\infty}-1}}  \frac{e^{\varphi_{n}+ z\psi_{n}}  }{  \sqrt{a_{n}}} (1+ o(1)), \q n\to\infty,
\q  \forall z\in{\Bbb C}.
\]
The resolvent  of the Jacobi operator $J$ is again determined by   formula \e{eq:RRes}, but, in contrast to the case $|\beta_{\infty}|<1$,  its singularities are due to zeros of the denominator $\Omega (z)$ only.
Therefore the spectrum of $J$ is discrete.

     \section{Volterra equations}
     
  
  In this   section we begin our study  of   Jacobi operators  with increasing recurrence coefficients  in a systematic way. To be precise, we accept assumptions \e{eq:Haupt} and sometimes distinguish the cases  $|\beta_{\infty}|> 1$ and $|\beta_{\infty}|<1$.  Below   $ z\in {\Bbb C}$ if $|\beta_{\infty}|>1$ and $z\in \clos \Pi$ where $\Pi={\Bbb C} \setminus {\Bbb R}$ if  $|\beta_{\infty}|<1$. The Carleman condition \e{eq:Carl} is not required unless specified explicitly. The function $\zeta (z)$ is defined by formula  \e{eq:ome}.
  
  Here we reduce the Jacobi equation  \e{eq:Jy} for the Jost solutions $f_{n}(z)$ defined by their asymptotics for $n\to\infty$ to a Volterra equation. Then we solve this equation  by iterations. The results of this section give an analytical basis for a construction of the Jost solutions and building of spectral theory of Jacobi operators in subsequent sections.

   \subsection{Ansatz}
   
 

We follow the scheme described in
   Sect.~2.5.   Recall that the numbers $\alpha_{n}$, $\beta_{n}$ and $z_{n}$ are defined by equalities \e{eq:aabb} and \e{eq:aabb1} and $\z_{n}= \z (z_{n})$.
We also set
\begin{equation}
\varkappa_{n}= \sqrt{\frac{a_{n+1} } {a_{n}}} , \q
  k_{n} =\frac{\varkappa_{n-1} } {\varkappa_{n}} =\frac{a_{n}  } {\sqrt{a_{n-1} a_{n+1}} } .
\label{eq:Gr4K}\end{equation}

  Let us choose $\varepsilon_{0}>0$  so small that $\pm 1\not\in (\beta_{\infty}-\varepsilon_{0}, \beta_{\infty} +\varepsilon_{0})$.  We always suppose that the values of the spectral parameter $z$  are bounded, that is, $| z| <\rho_{0}$ for some $\rho_{0}<\infty$ and fix $N_{0}$ in  such a way that $|z_{n}-\beta_{\infty}| <\varepsilon_{0}$ for $n\geq N_{0}$ and $| z| <\rho_{0}$.  Then $z_{n}$ are separated from the singular points $1$ and $-1$:
  \begin{equation}
 |z_{n} \pm 1| \geq \varepsilon >0.
\label{eq:aabb5}\end{equation}


                                                                                                                                                      
Our first goal is to find an Ansatz $ {\cal A}_{n} (z) $ for Jost solutions in the case of recurrence coefficients satisfying  condition \e{eq:Haupt}. Let us seek  it  in the form
\begin{equation}
  {\cal A}_{n} (z)  = p_{n}    q_{n} (z)   
\label{eq:GX}\end{equation}
with $ q_{n} (z) $ defined as  product \e{eq:re1}
 and
  a suitable sequence $p_{n}$. We have to calculate
     relative remainder \e{eq:Gry} in   equation \e{eq:Jy}.
  Since
\[
\frac{ {\cal A}_{n-1}}{ {\cal A}_{n}} =\frac{p_{n-1} q_{n-1}}{p_{n} q_{n}} =\frac{p_{n-1}  }{p_{n}  } \frac{1  }{\zeta _{n-1}  } ,
\]
   expression  \e{eq:Gry} equals
\begin{align}
 r_{n} (z) = & {\bf s}_{n-1}^{-1} \zeta _{n-1}^{-1} -2 z_{n}+
 k_{n} {\bf s}_{n} \zeta _{n}
 \nonumber\\
 = &\big(  {\bf s}_{n-1}^{-1}\zeta_{n-1}^{-1} -\zeta _{n}^{-1} \big)
 +
\big( k_{n} {\bf s}_{n} -1\big)\zeta _{n}
\label{eq:Grr1}\end{align}
where
\[
{\bf s}_{n}=\sqrt{\frac{a_{n+1}}{a_{n}}} \: \frac{p_{n+1}}{p_{n}} .
\]
Since $\zeta_{n-1}^{-1} -\zeta _{n}^{-1}\to 0$ and $k_{n}\to 1$ as $n\to\infty$, in order to estimate expression
\e{eq:Grr1} we have to set ${\bf s}_{n}=1$. Then
 $p_{n}= a_{n}^{-1/2}$,  formula \e{eq:GX} coincides with \e{eq:AnsDD} and
 \begin{equation}
 r_{n} (z)  =   \big(  \z _{n-1}^{-1}-\z _{n}^{-1}\big)+ (k_{n}-1)\z_{n},\q n\geq N_{0}+1.
\label{eq:Gr6}\end{equation}

Let us state this auxiliary result.

 \begin{lemma}\label{Rr}
  Let the Ansatz $ {\cal A}_{n} (z)$ be defined by  formula  \e{eq:AnsDD}. 
    Then  remainder \e{eq:Gry}  admits   representation
\e{eq:Gr6}
where $z_{n}$ and $k_{n}$ are given by   \e{eq:aabb1} and \e{eq:Gr4K}, respectively.
\end{lemma}

 \subsection{Multiplicative substitution}


Let $ {\cal A}_{n} (z)  $ be defined by formula \e{eq:AnsDD}.
We are looking for solutions 
  $f_{n} (z)$  of the   Jacobi equation \e{eq:Jy} satisfying   condition \e{eq:JostAs}.
 We first construct a sequence $f_{n} (z)$ for large $n $. Then  $f_{n} (z)$ is extended to all $n $ as a solution of   difference  equation \e{eq:Jy}. 
 By analogy with the continuous case,   the sequence $f(z)= (f_{n}( z ))_{n=-1}^\infty$   will be called
  the   Jost solution of the Jacobi equation \e{eq:Jy}. 
  
The uniqueness of  solutions with asymptotics  \e{eq:JostAs} is almost obvious.

\begin{lemma}\label{uniq}
Equation \e{eq:Jy} may have only one solution   $f_{n} (z)$ satisfying condition  \e{eq:JostAs}.
 \end{lemma}

\begin{pf}
Let $\ti{f}_{n} (z)$ be another solution of \e{eq:Jy} satisfying   \e{eq:JostAs}. Then the Wronskian \e{eq:Wr} of these solutions
equals
\[
\{f,\ti{f}\}=a_{n} {\cal A}_{n} {\cal A}_{n+1} o(1)=  \sqrt{\frac{a_{n}}{a_{n+1}}} q_{n}q_{n+1} o(1),\q n\to\infty ,
\]
where $q_{n}=q_{n} (z)$ are
  numbers   \e{eq:re1}. Clearly,    $|q_{n}|\leq 1$.  Since $a_{n}\to\infty$, the ratios $a_{n}/a_{n+1}$ are also bounded, at least for some subsequence of $n\to\infty$. Thus, $\{f,\ti{f}\}=0$ whence  $\ti{f} =Cf$. Here $C=1$ by virtue again of condition \e{eq:JostAs}.
\end{pf}

Note  a symmetry relation
    \begin{equation}
f_{n} (\bar{z})=\ov{f_{n} (z)}
\label{eq:AcGG}\end{equation}
which is a consequence of the equality ${\cal A}_{n} (\bar{z})=\ov{{\cal A}_{n} (z)}$ and Lemma~\ref{uniq}.

 
For a construction of $f_{n}  (z)$, we make a multiplicative substitution
 introducing a  sequence
\begin{equation}
 u_{n} (z)=  {\cal A}_{n} (z)^{-1}  f_{n} (z), \q n\in {\Bbb Z}_{+}.
\label{eq:Gs4}\end{equation}
 Then  \e{eq:JostAs} is equivalent to  condition  
\e{eq:A12a}.

  We first derive a difference equation for $ u_{n} (z)$ which will be subsequently reduced to a Volterra ``integral" equation. 


\begin{lemma}\label{Gs}
Let  $ {\cal A}_{n} (z) $ and  $r_{n} (z) $ be given  by formulas \e{eq:AnsDD}  and \e{eq:Gr6}, respectively.
 Then 
 equation  \e{eq:Jy} for a sequence $ f_{n} (z)$ is equivalent to the equation
\begin{equation}
k_{n} \z _{n}  ( u_{n+1} (z)- u_{n} (z)) -     \z _{n-1}^{-1}   ( u_{n} (z)- u_{n-1} (z))=-     r_{n} (z) u_{n} (z) 
\label{eq:Gs5}\end{equation}
for  sequence  \e{eq:Gs4}. 
 \end{lemma}

\begin{pf}
Substituting  expression $f_{n}  = {\cal A} _{n}  u_{n} $ into \e{eq:Jy} and using definition \e{eq:Gry} of $r_{n}$, we see that
\begin{multline}
( \sqrt{a_{n-1}a_{n}}{\cal A}_{n} )^{-1}\Big(  a_{n-1} f_{n-1} + ( b_{n} -z)f_{n} + a_{n} f_{n+1}\Big)
\\= r_{n} u_{n}  +
( \sqrt{a_{n-1}a_{n}}{\cal A}_{n} )^{-1}\Big(  a_{n-1} {\cal A}_{n-1} (u_{n-1}- u_{n})
+ a_{n} {\cal A}_{n+1} (u_{n+1}- u_{n})\Big).
  \label{eq:Gs6} \end{multline}
  In view of the equality 
  \[
\frac{{\cal A}_{n+1}}{{\cal A}_{n}} = \sqrt{\frac{a_{n}}{a_{n+1}}}\z _{n},
\]
the right-hand side here can be written as
\[
 r_{n} u_{n}  + \z_{n-1}^{-1} (u_{n-1} -u_{n} ) + k_{n}  \z _{n}( u_{n+1} -u_{n})    .
\]
Therefore equality \e{eq:Gs6}  implies that equations \e{eq:Jy}  and \e{eq:Gs5} are equivalent.
Finally, we note that representation \e{eq:Gr6} for  the coefficient $r_{n}$   follows from Lemma~\ref{Rr}.
 \end{pf}

 \subsection{Volterra integral equation}

Here  we reduce    difference equation  \e{eq:Gs5}   to a Volterra equation 
     \begin{equation}
 u_{n}(  z )= 1 + \sum_{m=n+1}^\infty G_{n,m} (z)     r_{m} ( z )u_{m} (  z )   ,
\label{eq:A17}\end{equation}
  where   the kernels  $G_{n,m}(z)$ are defined as follows. 
 We set
  \begin{equation}
\sigma_{n}=    \z _{n} \z _{n-1} ,\q S_{n}= \sigma_{1} \sigma_{2}  \cdots \sigma_{n-1}, 
\label{eq:Gp2a}\end{equation}
and
\begin{equation}
G_{n,m} (z)= - \varkappa_{m-1}^{-1}  \z _{m}^{-1} S_{m+1} \sum_{p=n+1}^m \varkappa_{p-1}  S_{p}^{-1}  , \q m>n,\label{eq:Gp6}\end{equation}
where $\z_{n}= \z(z_{n})$ and the numbers $z_{n}$, $\varkappa_n$ are given by \e{eq:aabb1}, \e{eq:Gr4K}, respectively.

 

 In all estimates below,  values of the spectral parameter $z$ are bounded   and $n$
   are sufficiently large. Our
 estimates of   the remainder  $r_{n} (z)$ are practically the same in the case $|\beta_{\infty}|<1$ and $|\beta_{\infty}|>1$.   First we set
  \begin{equation}
\varrho _{n} = 1+  \frac{ z_{n-1}+  z_{n}}{\sqrt{ z_{n-1}^2-1}+\sqrt{ z_{n}^2-1}}
\label{eq:Gs+1}\end{equation}
and note a representation 
 \begin{equation}
\z _{n-1}^{-1} - \z _{n}^{-1}= ( z_{n-1} -  z_{n}) \varrho _{n} 
\label{eq:Gs}\end{equation}
which is a consequence of definitions \e{eq:ome} and   \e{eq:aabb1}.
 It follows from   relation   \e{eq:aabb2} where $\beta_{\infty}^2\neq 1$ that estimate  \e{eq:aabb5} is satisfied.
 Thus   \e{eq:Gs+1},  \e{eq:Gs} yield  an estimate
 \begin{equation}
|\z_{n-1}^{-1} - \z_{n}^{-1} | \leq C | z_{n-1} -  z_{n}| \leq C_{1} (| \alpha_{n-1} -  \alpha_{n}|
+ | \beta_{n-1} -  \beta_{n}|).
\label{eq:Gs2}\end{equation}
In view of representation \e{eq:Gr6},
this leads to the following assertion  where notation \e{eq:diff} is used.

   \begin{lemma}\label{Gr}
   Under assumption 
 \e{eq:Haupt}   remainder \e{eq:Gry} satisfies an estimate
\[
| r_{n} (z) | \leq C  (| \alpha_{n-1} -  \alpha_{n}|+ | \beta_{n-1} -  \beta_{n}| )+  | k_{n} -  1| .
\]
 In particular, if
   \begin{equation}
(\alpha_{n}' )\in \ell^1 ({\Bbb Z}_{+} ), \q  (\beta_{n}' )\in \ell^1 ({\Bbb Z}_{+} )
\label{eq:Gr8}\end{equation}
and
 \begin{equation}
 ( k_{n} -1 )\in \ell^1 ({\Bbb Z}_{+} ),
\label{eq:Gr6b}\end{equation} 
then   
\begin{equation}
 ( r_{n} (z) ) \in \ell^1 ({\Bbb Z}_{+}) .
\label{eq:Gr9}\end{equation}
\end{lemma}


  \begin{remark}\label{Grem1}
  Let the numbers  $\varkappa_n$ be defined by the first equality \e{eq:Gr4K}. 
If inclusion  \e{eq:Gr6b} is true, then there exists a finite limit
\begin{equation}
\lim_{n\to\infty}  \varkappa_{n} =\varkappa_\infty \q \mbox{where}\q \varkappa_\infty\geq 1
\label{eq:Gr6a}\end{equation}
and
\begin{equation}
 (\varkappa_{n}'  ) \in  \ell^1 ({\Bbb Z}_{+}) .
\label{eq:Gs9b}\end{equation}
Moreover, $\varkappa_\infty=1$ if  the Carleman condition \e{eq:Carl}  is satisfied. 

Indeed, by definition \e{eq:Gr4K}, we have
\[
\ln \varkappa_{n} - \ln \varkappa_0=-\sum_{m=1} ^n \ln k_{m}
\]
where the series on the right converges according to \e{eq:Gr6b}.
It follows  that the limit in \e{eq:Gr6a} exists.  
Since
\[
 \varkappa_{n} - \varkappa_{n-1} =\frac{ \varkappa_{n}^2}{ \varkappa_{n} + \varkappa_{n-1} }(1+k_{n})(1-k_{n}),
 \]
 inclusion \e{eq:Gs9b}  is a direct consequence of assumption  \e{eq:Gr6b}.
 If  $\varkappa_\infty<1$, then $a_{n}\to 0$ as $n\to\infty$  which is   excluded by the first assumption  \e{eq:Haupt}.
 If $\varkappa_\infty>1$, then $a_{n}  \geq c l^n$ for some $l>1$ and $c>0$ so that  condition \e{eq:Carl}    is violated.   
\end{remark}
 
\begin{example}\label{GKy} 
Condition    \e{eq:Gr6b}
is of course  satisfied if $a_{n}= \nu (n+1)^p$ for some $\nu >0$ and $p>0$, but it allows  much more rapid growth of $a_{n}$ as $n\to\infty$.  
Consider, for example,   $a_{n}= \nu x^{ n^q}$ where $x>1$.  If $q<1$, then
\[
k_{n}^2= \exp\big( \nu (2n^q-(n+1)^q -(n-1)^q)\big)= \exp\big( O(n^{q-2})\big)=1+   O(n^{q-2})
\]
so that condition    \e{eq:Gr6b}  holds. In this case $\varkappa_\infty=1$. If $q= 1$, then $\varkappa_n=\sqrt{x}$ and  $k_{n}=1$ for all $n$.
On the  contrary,    condition  \e{eq:Gr6b} fails and $\varkappa_{n}\to\infty$ as $n\to\infty$ if $a_{n}= \nu x^{ n^q}$ with  $q>1$.
      \end{example}

  \subsection{Solution by iterations}
   
First,  we estimate   ``kernels" \e{eq:Gp6}. As before,  we use that according to  \e{eq:aabb3} or \e{eq:aabb4} the numbers $\z  (z_{n})$ are separated from the points $1$ and $-1$ for sufficiently large $n$, say $n\geq N_{0}$.
The estimates below are  quite straightforward in the case $|\beta_{\infty}|>1$.

 \begin{lemma}\label{GS2+}
Let  assumptions  \e{eq:Haupt}  where $|\beta_{\infty}|>1$ and  
\begin{equation}
0 < c_{1} \leq \varkappa_{n} \leq c_{2} <\infty
\label{eq:kk}\end{equation}
  be satisfied.   Then
 kernel \e{eq:Gp6} is bounded uniformly in $m>n\geq N_{0}$:
\begin{equation}
|G_{n,m} (z)|\leq    C<\infty.
\label{eq:AbG}\end{equation}
 \end{lemma}

 \begin{pf}
 Let $\sigma_n=\z_n \z_{n-1}$.
  By virtue of \e{eq:aabb3}, we have 
 \[
|\sigma_n|\leq \sigma_{\infty}\q\mbox{where}\q \sigma_{\infty}<1 \q\mbox{for}\q  n\geq N_{0}.
\]
Using also \e{eq:kk} we see that 
\[
|G_{n,m}| \leq C 
  \sum_{p=n+1}^m | \sigma_{p} \cdots \sigma_{m}|\leq C_{1}  
  \sum_{p=n+1}^m \sigma_{\infty}^{m-p+1} =C_{1}  \frac{1- \sigma_{\infty}^{m-n}}  
  {\sigma_{\infty}^{-1}  -1 }   \leq C_{1}  \frac{ \sigma_{\infty} }  
  {1-\sigma_{\infty}  }    
\]
which implies \e{eq:AbG}.  
   \end{pf}

In the case $|\beta_{\infty}|<1$ we have to ``integrate by parts" in \e{eq:Gp6}. This requires some additional assumptions.

 \begin{lemma}\label{GS2}
Let assumptions  \e{eq:Haupt} where $|\beta_{\infty}| < 1$, \e{eq:Gr8}, \e{eq:Gs9b} and \e{eq:kk} 
be satisfied. Then
 kernel \e{eq:Gp6} is bounded uniformly in $m > n\geq N_{0}$, that is, estimate \e{eq:AbG} holds.
 \end{lemma}
 
  \begin{pf} 
 Since $|\z_{p}|\leq 1$, it follows from 
  definition \e{eq:Gp2a} that
  \begin{equation}
|S_{m+1} S_{p}^{-1}| =|\sigma_{p} \cdots \sigma_{m} | \leq  1,\q p \leq m.
\label{eq:ssc}\end{equation}
Since $\sigma_{n}\to (\z^{(\pm)}_{\infty})^2\neq 1$ as $n\to\infty$, we see that
 \begin{equation}
|\sigma_n-1|\geq \varepsilon>0 .
\label{eq:ss}\end{equation}
By definitions \e{eq:diff} and \e{eq:Gp2a} we have
\[
(S_{p } ^{-1})'=S_{p+1} ^{-1}  - S_{p } ^{-1}= ( \sigma_{p}^{-1} -1 )  S_{p}^{-1}, 
\]
and hence integrating by parts (that is, using formula \e{eq:Abel}), we find that
\begin{multline}
\sum_{p=n+1 }^ { m} \varkappa_{p-1} S_{p}^{-1} = \sum_{p=n+1 }^ { m} \varkappa_{p-1}  ( \sigma_{p}^{-1}  -1)^{-1}(S_{p}^{-1})'
 =\varkappa_{m-1} ( \sigma_{m}^{-1}  -1)^{-1}S_{m+1}^{-1}
 \\
 -  \varkappa_{n-1} ( \sigma_{n}^{-1}  -1)^{-1}S_{n+1}^{-1}-\sum_{p=n+1 }^ { m} (\varkappa_{p-2}  (\sigma_{p-2}^{-1}  -1)^{-1})'S_{p}^{-1} .
\label{eq:AbelG}\end{multline}
Note that  
 \begin{equation}
(( \sigma_{p-1}^{-1}  -1)^{-1})'=  ( \sigma_{p-1}  -1)^{-1}( \sigma_{p}  -1)^{-1} \sigma_{p-1}' 
\label{eq:sb}\end{equation}
where $\sigma_{p}'\in  \ell^1 ({\Bbb Z}_{+})$ by  virtue of \e{eq:Gs2}. Using also    \e{eq:ss},  we see that
 sequence \e{eq:sb} belongs to $\ell^1 ({\Bbb Z}_{+})$. Finally, taking \e{eq:Gs9b} into account we find that
\[
(\varkappa_{p-1}  (\sigma_{p-1}^{-1}  -1)^{-1})'\in  \ell^1 ({\Bbb Z}_{+}).
\]
Let us multiply identity \e{eq:AbelG} by $S_{m+1} $.
According to \e{eq:ssc} and \e{eq:ss} all three terms in the right-hand side of the equality obtained are bounded   for $n\geq N_{0}$. 
   \end{pf}

  
Now  we are in a position to   standardly  solve the Volterra equation \e{eq:A17} by iterations. First, we estimate these  iterations. 


  \begin{lemma}\label{GS3p}
  Set $u^{(0)}_n =1$   and 
  \begin{equation}
 u^{(k+1)}_{n}(z)= \sum_{m=n +1}^\infty G_{n,m} (z) r_{m} (z) u^{(k )}_m (z),\q k\geq 0,
\label{eq:W5}\end{equation}
for all $n\in {\Bbb Z}_{+}$. Then the estimates
 \begin{equation}
| u^{(k )}_{n} (z) |\leq \frac{C^k}{k!} \big(\sum_{m=n+1}^\infty |r_{m} (z)|\big)^k,\q \forall k\in{\Bbb Z}_{+}.
\label{eq:W6s}\end{equation}
are true for all sufficiently large $n$ with the same constant $C $ as in  \e{eq:AbG}. 
\end{lemma}

 \begin{pf}
  Suppose that \e{eq:W6s} is satisfied for some $k\in{\Bbb Z}_{+}$. We have to check 
 the same estimate (with $k$ replaced by $k+1$ in the right-hand side)  for $ u^{(k+1)}_{n}$.  
 Set
 \[
 { R}_{m}= \sum_{p=m +1}^\infty  | r_p | .
 \]
  According to definition \e{eq:W5}, it  follows from estimates \e{eq:AbG} and  \e{eq:W6s} that
   \begin{equation}
| u^{(k +1)}_{n} |\leq  C \sum_{m=n +1}^\infty  | r_m | |{ u}_{m}^{(k)} | \leq \frac{C^{k+1}}{k!}  \sum_{m=n +1}^\infty  | r_m | { R}_{m}^k.
\label{eq:V7}\end{equation}
Observe that
 \[
{ R}_{m}^{k+1}+ (k+1)   | r_{m}|  { R}_{m}^k  
\leq
{ R}_{m-1}^{k+1},
\]
and hence, for all $N\in{\Bbb Z}_{+}$,
   \[
 (k+1)  \sum_{m=n +1}^N  | r_m | {  R}_{m}^k 
 \leq 
 \sum_{m=n +1}^N  ( {  R}_{m-1}^{k+1}-{  R}_{m}^{k+1})
= { R}_{n}^{k+1}-{  R}_{N}^{k+1}\leq  {  R}_{n}^{k+1}.
 \]
Substituting this bound into   \e{eq:V7}, we obtain estimate \e{eq:W6s} for $u^{(k +1)}_{n}$.
    \end{pf}
    
    This lemma enables us to construct a solution $u _{n} (z)$
of  equation \e{eq:A17} as a convergent series.

  \begin{theorem}\label{GS3}
  Let the assumptions  \e{eq:Haupt}    as well as   \e{eq:Gr8},  \e{eq:Gr6b}  be satisfied.
  Suppose that $|z|<\rho_{0}$ for some $\rho_{0}<\infty$ and $n\geq N_{0}(\rho_{0})$ for sufficiently large 
  $N_{0}(\rho_{0})$.
  Then equation \e{eq:A17} has a unique bounded solution $u_{n} (z)$, and condition \e{eq:A12a} is satisfied.
 If $|\beta_{\infty}| >1$, the functions $u_{n} (z)$  are analytic in $z$ for $|z|<\rho_{0}$ and $n\geq N_{0}(\rho_{0})$. If  $|\beta_{\infty}| <1$, the same is true in the complex plane  cut along $\Bbb R$; in this case the functions $u_{n} (z)$ are continuous up to the cut.
\end{theorem}

 \begin{pf}  Set
     \begin{equation}
   u_{n} =\sum_{k=0}^\infty u^{(k)}_{n} 
\label{eq:W8}\end{equation}
where $u^{(k)}_0=1$ and $u^{(k)}_{n}$, $k\geq 1$,  are defined by recurrence relations \e{eq:W5}.
Estimate \e{eq:W6s} shows that this series is absolutely convergent. Using the Fubini theorem to interchange the order of summations in $m$ and $k$, we see that
\[
   \sum_{m=n+1}^\infty G_{n,m}     r_{m}  u_{m}    =   \sum_{k=0}^\infty\sum_{m=n+1}^\infty G_{n,m}   r_{m}  u_{m}^{(k)} = \sum_{k=0}^\infty  u_{n}^{(k+1)}=- 1+ \sum_{k=0}^\infty  u_{n}^{(k)}= -1+ u_{n}.
\]
This is equation  \e{eq:A17} for sequence \e{eq:W8}. Condition \e{eq:A12a}  for this sequence   also follows from \e{eq:W6s} because  $r_{n}  (z)$ satisfies \e{eq:Gr9}. 
\end{pf} 


      \section{Jost solutions}
      
      We here use the results of the preceeding section to construct the Jost solutions for the Jacobi equation \e{eq:Jy}. 
      
            \subsection{An auxiliary  difference equation}
            
            The first step is to transform `` integral"  equation   \e{eq:A17}
      to  difference equation  \e{eq:Gs5}.
 
  \begin{lemma}\label{GS4}
  Let $r_{n}(z)$ and $G_{n,m} (z)$ be given by formulas  \e{eq:Gr6} and  \e{eq:Gp6}, respectively. Then the solution $u_{n}  (z)$ of     integral  equation   \e{eq:A17} satisfies an identity
   \begin{equation}
 \varkappa_{n}^{-1} (u_{n+1} - u_{n} )=   \sum_{m=n+1}^\infty      S_{n+1}^{-1} S_{m+1} \varkappa_{m-1}^{-1}  \z _{m} ^{-1}  r_{m} u_{m} ,
  \label{eq:A17Mb}\end{equation}
 and  difference equation  \e{eq:Gs5}.
 \end{lemma}
 
  \begin{pf}
  It follows from \e{eq:A17}    that
  \begin{equation}
 u_{n+1} - u_{n}=   \sum_{m=n+2}^\infty (G_{n+1, m}-G_{n, m})r_m u_{m}  -
 G_{n, n+1} r_{n+1} u_{n+1}.
  \label{eq:A17Ma}\end{equation}
  Since according to  \e{eq:Gp6}  
\[
G_{n+1, m}-G_{n, m}= \varkappa_{n} \varkappa_{m-1}^{-1} \z _{m} ^{-1} S_{n+1}^{-1} S_{m+1} \q \mbox{and}\q
G_{n, n+1}= - \z _{n+1}^{-1}  S_{n+2} S_{n+1}^{-1} ,
\]
equality \e{eq:A17Ma} can be rewritten as
 \e{eq:A17Mb}.

Putting together   equality \e{eq:A17Ma} with the same equality for $n$ replaced by $n-1$, we see that
 \begin{align*}
  \varkappa_{n}^{-1} ( u_{n+1} - u_{n}) - \sigma_{n}^{-1}  \varkappa_{n-1}^{-1} ( u_{n} - u_{n-1})&
  \\
  =  \sum_{m=n+1}^\infty    S_{n+1}^{-1} S_{m+1}     \varkappa_{m-1}^{-1}  \z _{m} ^{-1} r_{m} u_{m}-
\sigma_{n}^{-1} \sum_{m=n}^\infty  S_{n}^{-1} S_{m+1}  &   \varkappa_{m-1}^{-1} \z _{m} ^{-1} r_{m} u_{m}.
  \end{align*}
Since  $S_{n+1}=\sigma_{n} S_{n}$, the right-hand side here equals $ -  \varkappa_{n-1}^{-1} \z _n ^{-1} r_{n}  u_{n}$, and hence the equation obtained coincides with \e{eq:Gs5}.  
   \end{pf}
   
     \subsection{Main result}
  
  Now we are in a position to construct solutions of  the Jacobi equation \e{eq:Jy} with asymptotics \e{eq:JostAs} as $n\to\infty$.
  We call them the Jost solutions. 
  Recall that, by Lemma~\ref{uniq},  equation \e{eq:Jy} may have only one solution with such asymptotics.  We fix some $\rho_{0}> 0$,  suppose that $|z|< \rho_{0}$ and choose $N_{0} = N_{0}(\rho_{0} )$ such that condition \e{eq:aabb5} is satisfied
for $|z|\leq \rho_{0}$ and $n\geq N_{0}$. 
    Then using Theorem~\ref{GS3} we construct $u_{n}  (z)$ for $n\geq N_{0}$ as solutions of  the   integral  equation   \e{eq:A17}. According to Lemma~\ref{GS4} this sequence  satisfies  also the difference equation  \e{eq:Gs5}. Now we define $f_{n} (z)$ for $n\geq N_{0}$ by   formula  \e{eq:Jost}.
      In view of Lemma~\ref{Gs} this sequence satisfies the Jacobi equation \e{eq:Jy}. Then we extend $  f_{n}( z )$ as a solution of  equation \e{eq:Jy} to all $n\in {\Bbb R}_{+}$.  Thus we have the following result.

    \begin{theorem}\label{GSS}
   Suppose that  assumptions  \e{eq:Haupt}    as well as    \e{eq:Gr8} and \e{eq:Gr6b}  are satisfied.
Let $u_{n}( z )$ be  the sequence constructed in Theorem~\ref{GS3}, and let $q_{n}( z )$ be product \e{eq:re1}. 
Suppose thar   $ z\in {\Bbb C}$ if $|\beta_{\infty}|>1$ and $z\in \clos \Pi$ where $\Pi={\Bbb C} \setminus {\Bbb R}$ if  $|\beta_{\infty}|<1$.  
   Then  the sequence $f_{n}( z )$ defined by the equality
   \[
    f_{n}( z )= a_{n} ^{-1/2}q_{n}( z ) u_{n} (z )
    \]
     satisfies equation \e{eq:Jy},  and it has  asymptotics
  \begin{equation}
f_{n}(  z )= a_{n} ^{-1/2}q_{n}( z )\big(1 +  o(1) \big), \q n\to\infty.
\label{eq:A22G}\end{equation} 
If $|\beta_{\infty}| >1$, then all functions $f_{n} (z)$  are analytic in $z\in{\Bbb C}$. If  $|\beta_{\infty}| <1$, the same is true for $z\in \Pi$; in this case the functions $f_{n} (z)$ are continuous up to the cut along $\Bbb R$.   
  \end{theorem}
 
 
 We  emphasize that the  Carleman condition \e{eq:Carl} is not required in Theorem~\ref{GSS}.

Recall that the sequence $f_{n }(z)$ obeys condition \e{eq:AcGG} and, in particular,  
  \begin{equation}
 f_{n}( \lambda- i0 )=\ov{f_{n}(  \lambda+ i0 )}, \q \lambda\in {\Bbb R} .
 \label{eq:AcG}\end{equation}
 
  We now define the Jost function $\Omega (z)$  by formula \e{eq:J-W}  where $f(z)$ is the Jost solution.
     The following result is a direct consequence of Theorem~\ref{GSS}.

\begin{corollary}\label{JOSTG}
If $|\beta_{\infty}| >1$, then   $\Omega (z)$  is an analytic function of $z\in{\Bbb C}$. If  $|\beta_{\infty}| <1$, then  $\Omega(z)$  is 
analytic in $z\in \Pi$ and is continuous up to the cut along 
 $\Bbb R$. Its  limit values on the cut satisfy the identity
   \[
 \Omega( \lambda- i0 )=\ov{\Omega (  \lambda+ i0 )}, \q \lambda\in {\Bbb R} .
 \]
\end{corollary}

 

   \subsection{Explicit formulas} 
   
   We here transform asymptotics  \e{eq:A22G} of the Jost solutions $f_{n}(z)$ to a simpler form.  
   This will require an assumption  that $a_{n}\to\infty$ not too slowly, i.e.
    \begin{equation}
\sum_{n=0}^\infty \alpha_{n}^{K+1} <\infty
 \label{eq:ss2}\end{equation}
 for some $K\in{\Bbb Z}_{+}$.

    We first discuss a behavior of 
    the numbers 
    \begin{equation}
     \z_{n}= \z (z_{n})= z_{n}-\sqrt{z_{n}^2-1}  
     \label{eq:zz}\end{equation}
     as $n\to\infty$.  Recall that  $z_{n}$    are defined by equalities \e{eq:aabb}.
           Differentiating relation  \e{eq:ome}, we see that
      \begin{equation}
 -\big( \ln\z(z)\big)'= (z^2-1)^{-1/2}=: \vartheta (z).
 \label{eq:Te}\end{equation}
 Here the function $ \vartheta (z) $ is defined on ${\Bbb C}\setminus [-1,1]$ and $\vartheta (z) >0$   for $z>1$.  Then $\vartheta (z) <0$   for $z<-1$ and
 $\vartheta (\lambda \pm i 0)  =\mp i (1-\lambda^{2})^{-1/2}$ if $\lambda\in (-1,1)$.
 We  fix the branch of the function $ \ln\z(z)$ by the condition $\ln \z (z) < 0$ for $z>1$. The functions $ \vartheta (z)$ and $\ln \z (z)$ are holomorphic
 on the set ${\Bbb C}\setminus [-1,1]$.   
  
 We start with an elementary assertion. 
 
     \begin{lemma}\label{E1}
     Let condition  \e{eq:Haupt} be satisfied, and let $n$ be sufficiently large.   
  If $|\beta_{\infty}|< 1$ and $\pm\Im z \geq 0$, then
       \begin{equation}
   \z (\alpha_{n}z +\beta_{n})= \z(\beta_{n}\pm i0) e^{\psi (\alpha_{n}z, \beta_{n}\pm i0)}
 \label{eq:Te1J}\end{equation}
 where, for an  arbitrary $K\geq 1$,   
   \begin{equation}
 \psi ( \alpha_{n}z , \beta_{n}\pm i0)= -\sum_{k=1}^K (k!)^{-1}  \vartheta ^{(k-1)} ( \beta_{n}\pm i0) (\alpha_{n}z)^k + O (\alpha_{n}^{K+1})
    \label{eq:Te2J}\end{equation}
     as $n\to\infty$.
   If $|\beta_{\infty}|> 1$, then relations \e{eq:Te1J} and \e{eq:Te2J} are true for all $z\in{\Bbb C}$ with $\beta_{n}\pm i0$ replaced by $\beta_{n}$.
     \end{lemma}

    \begin{pf}
    Let us set
         \begin{equation}
         \psi (t, \beta)=
  \ln \z (t+\beta)-   \ln \z (\beta )
 \label{eq:Te3}\end{equation}
 where  either $\beta \in (-1, 1)$ and $\pm \Im t\geq 0$ or $\beta\in (-\infty, -1)\cup (1,\infty)$ and $t\in{\Bbb C}$.
 Then $\psi (0,\beta)=0$  and  
     \begin{equation}
   \z (t+\beta)= \z(\beta) e^{\psi (t, \beta)} .
 \label{eq:Te1}\end{equation}
  Differentiating \e{eq:Te3} in $t$ and using definition \e{eq:Te}, we see that   
  \begin{equation}
  \psi_{t} (t, \beta)= - \vartheta  (t+\beta).
 \label{eq:Te4}\end{equation}
  Let  us take $K-1$ terms of the Taylor expansion of the function $\vartheta (\beta+t)$ at the point $t=0$  and then integrate \e{eq:Te4} over $t$ using that $  \psi (0, \beta)=0$. This yields the relation 
    \begin{equation}
 \psi (t,\beta)= -\sum_{k=1}^K (k!)^{-1}  \vartheta ^{(k-1)} (\beta\pm i0) t^k + O (|t|^{K+1})\q \mbox{as}
\q t\to 0.
   \label{eq:Te2}\end{equation}
      Moreover, in the case $|\beta|>1$, the numbers $\beta\pm i0$   can be here replaced by $\beta$.
       To obtain relations   \e{eq:Te1J} and \e{eq:Te2J}, we only have to apply \e{eq:Te1} and \e{eq:Te2}
   to $t=\alpha_{n}  z$  and $\beta=\beta_{n}$.
      \end{pf}
      
       \begin{corollary}\label{Te1c}
       If $|\beta_{n}| <1$, then
       \[
|   \z (\alpha_{n}z +\beta_{n})| = e^ {-|\Im z| \alpha_{n} (1-\beta_{n}^2)^{-1/2}+ O(\alpha_{n}^{2})}.
 \]
  \end{corollary}

   Now we are in a position to  find  an asymptotics as $n\to\infty$ of the product $q_{n} (z)$ defined by formula  \e{eq:re1}.   Let us set
    \begin{equation}
 \Theta_{n}^{(\pm)}    = \prod_{m=N_{0}}^{n-1} \z(\beta_{m}\pm i0),\q n\geq N_{0}+1,
  \label{eq:LM}\end{equation}
  and, for $K\geq 1$,   
   \begin{equation}
 L_{n}^{(\pm)}  (z; K) =- \sum_{k=1}^K (k!)^{-1} z^k \sum_{m=N_{0}}^{n-1}  \vartheta ^{(k-1)} (\beta_{m}\pm i0)\alpha_{m}^k 
 \q \mbox{for} \q \pm \Im z \geq 0 .
  \label{eq:L2}\end{equation}
    If $|\beta_{\infty}| > 1$, then all objects labelled by $``+ "$ and $``- "$  coincide with each other so that $\beta_{n}\pm i0$ can be replaced   by $\beta_{n} $. In this case \e{eq:L2} is true for all $z\in{\Bbb C}$.
    We often write   $L_{n}^{(\pm)}  (z)$ instead of $    L_{n}^{(\pm)}  (z;K) $ omitting the dependence  on  $K$.
  
   The following statement is a direct consequence of Lemma~\ref{E1}.
  
        \begin{lemma}\label{Te1}
        Let assumption \e{eq:ss2} be satisfied for some $K\in{\Bbb Z}_{+}$.
       Define the sequences $\Theta_{n}^{(\pm)}$ and $L_{n}^{(\pm)}  (z; K)$   by formulas \e{eq:LM} and \e{eq:L2}.
      Then         the limits
       \begin{equation}
\lim_{n\to\infty} q_{n} (z) / \Theta_{n}^{(\pm) } , \q \pm \Im z \geq 0,\q K=0,
  \label{eq:LM2-}\end{equation}
  and
 \begin{equation}
\lim_{n\to\infty} \frac{q_{n} (z)} { \Theta_{n}^{(\pm)}  \exp \big(  L_{n}^{(\pm)}  (z; K)\big)} , \q \pm \Im z \geq 0, \q K\geq 1,
  \label{eq:LM2}\end{equation}
  exist and are not zero. Moreover, if $|\beta_{\infty}| > 1$, then \e{eq:LM2-} an \e{eq:LM2}  hold for all $z\in{\Bbb C}$ and the indices $``\pm" $ may be omitted.
  \end{lemma}
  

Using this result  we can
 renormalize definition \e{eq:Jost}  of the Jost solutions $f_{n} (z)$ multiplying them by a constant non-zero factor
 that may depend on $z$ but not on $n$.   This allows us to state Theorem~\ref{GSS} in a more explicit form.
  
     \begin{theorem}\label{GSS+}
     In addition to the assumptions of Theorem~\ref{GSS} suppose  that condition  \e{eq:ss2} is satisfied  for some $K \in{\Bbb Z}_{+}$. Let the functions  $\Theta_{n}^{(\pm)}  $  be defined by formula \e{eq:LM}.  The functions $L_{n}^{(\pm)} (z)$ are defined by formula
       \e{eq:L2} for $K\geq 1$ and $L_{n}^{(\pm)} (z)=0$ for $K=0$. Then for all $z\in \clos \Pi$, equation \e{eq:Jy} has a unique solution $f_{n} (  z )$ with asymptotics
  \begin{equation}
f_{n}(  z )= a_{n} ^{-1/2} \Theta_{n}^{(\pm)}  \exp \big(  L_{n}^{(\pm)}  (z)\big) \big(1 + o( 1)\big),\q \pm \Im z\geq 0 ,  \q n\to\infty.
\label{eq:A22G+}\end{equation} 
Moreover, if $|\beta_{\infty}| > 1$, then \e{eq:A22G+}  is true  for all $z\in{\Bbb C}$ and the indices $``\pm" $ may be omitted.
 \end{theorem}
 
 \begin{remark}\label{GSrem}
     Under the assumptions of Theorem~\ref{GSS+}  we can simplify definition \e{eq:AnsDD} of the Ansatz setting
      \[
{ \cal A}_{n}^{(\pm)} (z)= a_{n} ^{-1/2} \Theta_{n}^{(\pm)}  \exp \big(  L_{n}^{(\pm)}  (z)\big).
\]
Then expression \e{eq:Gr6} for remainder \e{eq:Gry} is changed, but otherwise the proof  of Theorem~\ref{GSS+}  works with this Ansatz  without significant modifications. 
 \end{remark}
 
  \begin{remark}\label{carl-b}
In view of Remark~\ref{Grem1}, under assumption  \e{eq:Gr6b} condition  \e{eq:ss2} for $K=0$ is satisfied if and only if the Carleman condition  \e{eq:Carl} is violated.
 \end{remark}

According to Theorem~\ref{GSS+} the asymptotic formulas are simplest in the most singular case when  the Carleman condition \e{eq:Carl} is violated. Indeed, if  
assumption \e{eq:ss2} is satisfied for   $K =0$, then 
$   L_{n}^{(\pm)}  (z)=0$   and the asymptotics of $f_{n}(  z )$ is determined by  $ \Theta_{n}^{(\pm)}  $ only. 
Let us  give explicit formulas for   this term.  
  In view of definitions \e{eq:ome} and \e{eq:LM}   we have
   \begin{align}
\Theta_{n} = &(\sgn \beta_{\infty})^{n}\prod_{m= N_{0}}^{n-1}\big(|\beta_{m} | +|\sqrt{ \beta_{m}^2-1}\, | \big)^{-1}
    \nonumber\\
     =& \big(\sgn \beta_{\infty})^{n}\exp(-\sum_{m= N_{0} }^{n-1}\arccosh |\beta_{m}| \big),\q | \beta_{\infty}|>1.
     \label{eq:ss5o}\end{align}
     and
     \begin{equation}
   \Theta_{n}^{(\pm)} =  \prod_{m= N_{0}}^{n-1}\big( \beta_{m}\mp i |\sqrt{1- \beta_{m}^2}\, | \big) 
   =  
     \exp(\mp i\sum_{m= N_{0}}^{n-1}\arccos\beta_{m} \big),\q | \beta_{\infty}|<1.
     \label{eq:ss6o}\end{equation}
     In particular, for $b_{n}=0$, \e{eq:ss6o} reduces to $ \Theta_{n}^{(\pm)} = i^{\pm N_{0}} i^{\mp n}$. 
         As usual, the constant factor $i^{\pm N_{0}}$ can be neglected here.

In the Carleman case formula \e{eq:A22G+}
           contains the correction  $  \exp \big(  L_{n}^{(\pm)}  (z)\big)$ where $L_{n}^{(\pm)}  (z)$ is defined by  \e{eq:L2}.  This correction gets more complicated as $K$ increases. Let us first consider the leading term of $  L_{n}^{(\pm)}  (z)$. This yields  a simplified formulation  of Theorem~\ref{GSS+} where an estimate of the remainder is not very precise.
 
      \begin{theorem}\label{GHS+}
     In addition to the assumptions of Theorem~\ref{GSS} suppose  that the Carleman condition \e{eq:Carl}  is satisfied.
      Then  formula \e{eq:A22G+} holds true with $   \Theta_{n} $, $   \Theta_{n}^{(\pm)}$ given by \e{eq:ss5o}, \e{eq:ss6o} and
      \begin{equation}
{\sf L}_{n} (z) =   - z \sgn \beta_{\infty}\sum_{m= N_{0}}^{n-1} (\beta_{m}^2 -1)^{-1/2}\alpha_{m}    + o \big( \sum_{m=0}^{n-1}\alpha_{m}\big),\q  |\beta_{\infty}| > 1,
\label{eq:A3H+}\end{equation}
\begin{equation}
{\sf L}_{n}^{(\pm)} (z) =  \pm i z \sum_{m= N_{0}}^{n-1} (1-\beta_{m}^2)^{-1/2}\alpha_{m}  + o \big( \sum_{m= 0}^{n-1}\alpha_{m}\big),\q  |\beta_{\infty}|< 1,\q \pm \Im z\geq  0 .
\label{eq:A2H+}\end{equation}
 \end{theorem}

  Note that  if $K=1$, then the error terms in \e{eq:A3H+} and \e{eq:A2H+}  can be omitted. 
           In this case Theorem~\ref{GHS+} yields formulas of Sect.~2.6.

            Let us  write down explicitly two terms of expression  \e{eq:L2} for $K=2$.
             If $|\beta_{\infty}| > 1$, then  
  \begin{equation}
L_{n}   (z) = - z \sgn \beta_{\infty}\sum_{m= N_{0}}^{n-1} (\beta_{m}^2 -1)^{-1/2}\alpha_{m} 
 +    2^{-1}z^2 \sgn \beta_{\infty} \sum_{m= N_{0} }^{n-1}   (\beta_{m}^2 -1)^{-3/2}\beta_{m}\alpha_{m}^2.
\label{eq:L21+}\end{equation} 
If $|\beta_{\infty}| <1$ and $ \pm \Im z \geq 0$, we have
      \begin{equation}
 L_{n}^{(\pm)}  (z) = \pm i z \sum_{m= N_{0}}^{n-1} (1-\beta_{m}^2)^{-1/2}\alpha_{m} 
 \pm i    2^{-1}z^2 \sum_{m= N_{0} }^{n-1}   (1-\beta_{m}^2)^{-3/2}\beta_{m}\alpha_{m}^2 .
\label{eq:L21}\end{equation}

Moreover, if 
\[
\sum_{m=0 }^\infty   |\beta_{m}|\alpha_{m}^2 <\infty
\]
(in particular, if \e{eq:ss2}  is true with $K= 1$), then the terms with $z^2$ in the right-hand sides of  \e{eq:L21+} and \e{eq:L21} can be also omitted.   Thus we again recover formulas of Sect.~2.6.

In Sect.~5.2,  we need the following
  
      \begin{corollary}\label{GSSn}
      If $|\beta_{\infty}|<1$, then
   \begin{equation}
    |f_{n} (z)|^2\leq \alpha_{n}\exp\Big(- c \sum_{m= 0}^{n-1} \alpha_m\Big)
 \label{eq:L22}\end{equation}
 for all $z\in{\Bbb C}$ with $\Im z\neq 0$ and some $c=c(z)>0$.
        \end{corollary}
        
          \begin{pf}
         Observe that $|\Theta_{n}^{(\pm)}|=1$ and  according to Corollary~\ref{Te1c}
          \[
   \Re L_{n}^{(\pm)}  (z) =- |\Im z  |  \sum_{m= N_{0}}^{n-1} \frac{\alpha_{m}} {\sqrt{1-\beta_{m}^2}}+ o( \sum_{m= 0}^{n-1}  \alpha_{m}).
    \]
    Therefore \e{eq:L22} is a direct consequence of \e{eq:A22G+}.
        \end{pf}

  \section{Small diagonal elements: orthogonal polynomials  on the real axis and the spectral measure}
  
  In this section we accept    assumption  \e{eq:Haupt}  where $|\beta_{\infty}|< 1$ and the Carleman condition \e{eq:Carl}.  Our goal   is to     find an asymptotic behavior of the orthonormal polynomials $P_{n} (z)$ as $n\to\infty$  for $z=\lambda\in{\Bbb R}$. 
At  the same time we will show that the spectrum of the Jacobi operator $J$ is absolutely continuous, covers the whole real line and obtain an expression for the spectral measure of  $J$. 
 

           \subsection{Asymptotics on the continuous spectrum} 
    
   Let us proceed  from Theorem~\ref{GSS} where $z=\lambda\pm i0$ for $\lambda \in{\Bbb R}$. Recall that $q_{n}(\lambda\pm i0)$ is product \e{eq:re1} and
     \begin{equation}
\z_{n} (\lambda \pm i0)= \lambda_{n}\mp i\sqrt{ 1-\lambda_{n}^2}= e^{\mp i \arccos\lambda_{n}}
       \label{eq:wwm}\end{equation} 
where
       \[
 \lambda_n  = \alpha_{n} \lambda+\beta_{n}
\]
with the numbers   $\alpha_{n} $,  $\beta_{n} $  defined by formulas \e{eq:aabb}. 
     According to \e{eq:aabb2} we have
 \begin{equation}
\lim_{n\to\infty} \lambda_n=\beta_{\infty}
\label{eq:aabx}\end{equation}
so that $\lambda_{n}\in (-1,1)$ if $n$ is sufficiently large.  It follows from   \e{eq:wwm} that
      \begin{equation}
q_{n}(\lambda\pm i0)= e^{\mp i \varphi_{n } (\lambda)} \q \mbox{where}
\q \varphi_{n } (\lambda)= \sum_{m=N_{0}}^{n-1}\arccos\lambda_{m} .
      \label{eq:ww3}\end{equation}

 

Let us now state a particular case  of  Theorem~\ref{GSS}   for  the   case
 $z=\lambda\pm i0$ where $\lambda\in \Bbb R$.

 \begin{theorem}\label{As+}
 Let the assumptions of Theorem~\ref{GSS}   be satisfied and $|\beta_{\infty}| <1$.   Then  for all $ \lambda\in \Bbb R$, the equation
  \[
  a_{n-1} f_{n-1} + b_{n}  f_{n}  + a_{n}   f_{n+1}  =\lambda  f_{n} 
  \]
  has solutions $f  (\lambda+ i0)=(f_{n} (\lambda+ i0))$  and $f (\lambda- i0)= (f_{n} (\lambda- i0))
   $  with asymptotics 
   \begin{equation}
 f_{n}(\lambda\pm i0)= a_{n}^{-1/2}  q_{n}(\lambda\pm i0) (1+o(1)), \q n\to\infty.
   \label{eq:Jost1+}\end{equation}
   These sequences are complex conjugate to each other, that is, relation \e{eq:AcG} holds.
 \end{theorem}
 
 
   The following result shows that these solutions  are linearly independent.  Recall that the Wronskian of  $f (\lambda+ i0)$  and $f  (\lambda- i0)$ is defined by the relation
       \begin{equation}
\{ f (\lambda+i0), f(\lambda-i0)\}= a_{n} \big(f_{n}(\lambda+i0)  f_{n+1}(\lambda-i0)-f_{n}(\lambda-i0)  f_{n+1}(\lambda+i0)\big). 
       \label{eq:wwn}\end{equation} 

\begin{lemma}\label{Asc}  
The Wronskian  of $f (\lambda+ i0)$  and $f  (\lambda- i0)$  
equals
   \begin{equation}
\{f (\lambda+i0),  f (\lambda-i0)\}=  2 i   \sqrt{1- \beta^2_{\infty}}\neq 0,
      \label{eq:ww}\end{equation} 
and hence these solutions  are linearly independent.
\end{lemma}

  \begin{pf}
The right-hand side of \e{eq:wwn} does not depend
on $n$, and so we can calculate this expression for $n\to\infty$.    Therefore using \e{eq:Jost1+}, we find  that
    \[
\{ f (\lambda+i0), f(\lambda-i0)\}= \sqrt{\frac{a_{n}}{a_{n+1}}} \big(q_{n}(\lambda+i0)  q_{n+1}(\lambda-i0)-q_{n}(\lambda-i0)  q_{n+1}(\lambda+i0)+ o(1)\big) 
\]
whence
     \[
\{ f (\lambda+i0), f(\lambda-i0)\}= \sqrt{\frac{a_{n}}{a_{n+1}}} \big(\z (\lambda_{n}-i0)- \z (\lambda_{n}+i0)+ o(1)\big) 
   \]
      as $ n\to \infty $.
Observe that $a_{n} a_{n+1}^{-1}\to 1$ under the assumptions of Theorem~\ref{GSS}.  Thus, relation \e{eq:ww} is a direct consequence of \e{eq:wwm}  and \e{eq:aabx}.
 \end{pf}

Now we are in a position to find an asymptotic behavior of the polynomials $P_{n} (\lambda)$ for $\lambda\in \Bbb R$,  that is, on the continuous spectrum of the Jacobi  operator $J$.
Since the  Jost solutions $f(\lambda\pm i0)=( f_{n} (\lambda\pm i0))$  are linearly independent  and $\ov{P_{n} (\lambda)} =P_{n} (\lambda)$, 
  we see that
 \[
P_{n} (\lambda)= \ov{c(\lambda)} f_{n} (\lambda+i0)+ c(\lambda)  f_{n} (\lambda-i0)
\]
for some complex constant $c (\lambda)$. Taking the Wronskian of this equation with $ f  (\lambda+i0)$, we can  express $c(\lambda)$ via Wronskian     \e{eq:J-W}: 
  \[
- c (\lambda)\{ f (\lambda+i0),  f (\lambda-i0)\}=\{P (\lambda),  f (\lambda+ i0)\}  =  \Omega (\lambda+i0) .
\]
In view of Lemma~\ref{Asc} this yields the following result.

\begin{lemma}\label{HH+}
 For all $ \lambda\in \Bbb R$, we have the representation 
  \begin{equation}
P_{n} (\lambda)=\frac{  \Omega(\lambda- i0) f _{n} (\lambda+i0) - \Omega (\lambda + i0) f _{n} (\lambda- i0)  }{ 2 i   \sqrt{1- \beta^2_{\infty}}},   \q n \in {\Bbb Z}_{+}. 
\label{eq:HH4L+}\end{equation}
 \end{lemma}
 
 Properties of  the Wronskians  $\Omega(\lambda \pm i0)$ are summarized in the following statement.
 
 \begin{theorem}\label{HX+}
 Let the assumptions of Theorem~\ref{GSS} be satisfied and $|\beta_{\infty}|<1$.
 Then the Wronskians $\Omega(\lambda + i0)$ and $\Omega (\lambda - i0)=
 \ov{ \Omega (\lambda + i0) }$  are continuous functions of $\lambda\in \Bbb R$ and
 \begin{equation}
\Omega (\lambda\pm i0) \neq 0 ,\q \lambda\in \Bbb R  .
\label{eq:HH5+}\end{equation}
 \end{theorem}
 
  \begin{pf}
  The functions 
 $ \Omega (\lambda\pm i0)$ are    continuous    by Corollary~\ref{JOSTG}.
  If $ \Omega (\lambda\pm i0)=0$, then according to \e{eq:HH4L+}
    $P_{n} (\lambda)=0$ for   all $n\in{\Bbb Z}_{+}$.  However, 
 $P_0 (\lambda)=1$ for all $\lambda$.
    \end{pf}

Let us set
 \begin{equation}
\kappa (\lambda) = | \Omega(\lambda+i0)|, \q
- \Omega (\lambda\pm i0) =  \kappa (\lambda) e^{\pm i \eta (\lambda) } .
\label{eq:AP+}\end{equation}
In the theory of short-range perturbations of the Schr\"odinger operator, the functions $\kappa (\lambda) $ and $\eta (\lambda)$ are known as the limit amplitude and the limit phase, respectively; the function   $\eta (\lambda)$ is also called the scattering  phase or the   phase shift.    Definition \e{eq:AP+}     fixes the phase $\eta (\lambda)$  only up to a term $2\pi m$ where $m\in{\Bbb Z}$.  We emphasize that the amplitude $\kappa (\lambda) $ and the phase $\eta (\lambda)$
 depend on the   values of the coefficients $a_{n}$ and $b_{n}$ for all $n$,  and hence they are not determined by an asymptotic behavior of $a_{n}$, $b_{n}$ as $n\to\infty$. 

Combined together,   relations \e{eq:Jost1+} and \e{eq:HH4L+} yield    asymptotics of   the  orthonormal  polynomials  $P_{n} (\lambda)$.  
  
 \begin{theorem}\label{Sz+}
 Let the assumptions of Theorem~\ref{GSS}   be satisfied  and $|\beta_{\infty}|< 1$.    Then, as $n \to\infty$,  the  polynomials $P_{n} (\lambda)$ have asymptotics 
   \begin{equation}
 P_{n} (\lambda)=   \kappa ( \lambda) (1- \beta^2_{\infty})^{-1/2}  a_{n}^{-1/2} \sin (\varphi_{n}  (\lambda) + \eta(\lambda) ) +  o(1) , \q  \lambda\in \Bbb R,
\label{eq:Sz+}\end{equation}
where the phases $\varphi_{n}  (\lambda) $, $ \eta(\lambda)$ and the amplitude $  \kappa ( \lambda)$ are given by formulas  \e{eq:ww3} and   \e{eq:AP+}.
 Relation \e{eq:Sz+} is uniform in $\lambda$ on compact subintervals of $\Bbb R$.
 \end{theorem}
 
 \subsection{Resolvent and spectral measure}


Under the assumption that the operator $\clos J_{\rm min} =: J$ is self-adjoint,  its resolvent $R(z)= (J-zI)^{-1}$  was constructed in Proposition~\ref{res}. To use formula \e{eq:RRes}, we only need to identify the solution of equation  \e{eq:Jy} satisfying condition \e{eq:asqrf1+}
with  the Jost solution.

  \begin{lemma}\label{asq2}
    In addition to the assumptions of Theorem~\ref{GSS} suppose  that the Carleman condition \e{eq:Carl}  is satisfied  and $|\beta_{\infty}|< 1$.  Then,
 for the Jost solution $f_{n} (z)$,  inclusion \e{eq:asqrf1+} holds true.   
 \end{lemma}
 
   \begin{pf} 
   In view of inequality \e{eq:L22}, we only have to show that
        \begin{equation}
   \sum_{n=1}^\infty \alpha_{n}\exp\Big(-  \sum_{m=0}^{n-1} \alpha_m\Big)<\infty.
   \label{eq:asqrf2}\end{equation} 
   For a proof of the convergence of this series, we can replace here $\alpha_{n}$ by $1-e^{-\alpha_{n}}$ and then observe that
      \begin{align*}
      \sum_{n=1}^N (1-e^{-\alpha_{n}})\exp\Big(-  \sum_{m=0}^{n-1} \alpha_m\Big)&
      \\
      =  \sum_{n=1}^N  \exp\Big(-  \sum_{m=0}^{n-1} \alpha_m\Big)- \sum_{n=2}^{N+1}  \exp\Big(-  \sum_{m=0}^{n-1} \alpha_m\Big)
     & =e^{-\alpha_{0}}- \exp\Big(-  \sum_{m=0}^N \alpha_m\Big)< e^{-\alpha_{0}}.
    \end{align*}
    This proves   \e{eq:asqrf2} and hence in view of  estimate \e{eq:L22}, inclusion \e{eq:asqrf1+}. 
          \end{pf}

According to Theorem~\ref{GSS},   $f_{n} (z)$, $n=-1, 0,1,\ldots$,  and, in particular,  $\Omega (z)$ are analytic     functions of $z\in{\Bbb C}\setminus {\Bbb R}$ continuous  up to the cut along ${\Bbb R}$. Taking also \e{eq:HH5+} into account, we obtain the following result.  Recall that the set ${\cal D}\subset \ell^2({\Bbb Z}_{+})$ consists of vectors with only a finite number of non-zero elements.

 \begin{theorem}\label{AC}
  In addition to the assumptions of Theorem~\ref{GSS} suppose  that the Carleman condition \e{eq:Carl}  is satisfied  and $|\beta_{\infty}|< 1$.
  Then
  \begin{enumerate}[{\rm(i)}]
 \item
For $\Im z\neq 0$, the  resolvent $R(z)=(J-zI)^{-1}$ of the Jacobi operator $J$ acts by formula \e{eq:RRes}.

 \item
 For all $u,v\in {\cal D}$, the functions $\la R(z) u,v \ra$ are  continuous in $z$ up to the cut along ${\Bbb R}$.
  \end{enumerate}
 \end{theorem}
 
 
 Note that the functions $f_{n} (z)$ depend on the choice of the parameter $N_{0}$ in definition \e{eq:re1}  of $q_{n} (z)$  but due to the factor $\Omega (z)^{-1} $ expression \e{eq:RRes} for the resolvent does not depend on it.
 
  Statement (ii)   is known as the limiting absorption principle. It implies
 
 \begin{corollary}\label{Rc}
The spectrum of the operator $J$   is absolutely continuous.
 \end{corollary}

Let us now consider the spectral projector $E(\lambda)$ of the operator $J$. By the Cauchy-Stieltjes-Privalov formula for $u,v\in {\cal D}$, its matrix elements satisfy the identity
 \begin{equation}
2\pi i \frac{d\la E (\lambda)u,v\ra} {d\lambda}=  \la R (\lambda+ i0)u,v\ra - \la R (\lambda- i0) u,v\ra.
\label{eq:Priv}\end{equation}
The following assertion is a direct consequence of  Theorem~\ref{AC}, part (ii).

\begin{corollary}\label{RE}
 For all $u,v\in {\cal D}$, the functions $\la E(\lambda) u,v \ra $ are  continuously differentiable in $\lambda\in {\Bbb R}$.
   \end{corollary}


Now we are in a  position to calculate  the spectral family  $d E(\lambda)$ in  terms of the Jost function.
Let us proceed from the identity \e{eq:Priv}.  Let $e_{0}, e_{1},\ldots, e_{n},\ldots$ be the canonical basis in the space $\ell^2 ({\Bbb Z}_{+})$. It follows from representation  \e{eq:RRes} that
\begin{equation}
\la R (z)e_{n}, e_{m}\ra =   \Omega (z)^{-1} P_{n} (z) f_{m}(z),\q n\leq m,
\label{eq:Rrpm}\end{equation}
whence
\[
\la R (\lambda\pm i0)e_{n}, e_{m}\ra =  \Omega (\lambda\pm i 0)^{-1} P_{n} (\lambda) f_{m}(\lambda\pm i 0),\q n\leq m.
\]
Substituting  this expression  into \e{eq:Priv} and using relation  \e{eq:AcG},
 we find that
  \begin{equation}
2\pi i \frac{d \la E (\lambda)e_n, e_m \ra } {d\lambda}= P_{n} (\lambda) \frac{ \Omega (\lambda-i 0)  f_{m}(\lambda+i 0)    -  \Omega (\lambda+i 0) f_{m}(\lambda-i 0)    }{|  \Omega (\lambda\pm i 0) |^2}
\label{eq:RE1}\end{equation}
 (note that $|\Omega (\lambda+i 0) |=|\Omega (\lambda-i 0) |$). Combining this representation with  formula \e{eq:HH4L} for $P_m (\lambda) $, we obtain    the following result.

 \begin{theorem}\label{SF}
  Let the assumptions of Theorem~\ref{AC}  be satisfied. Then the spectrum of the operator covers the whole line and, for all $n,m\in{\Bbb Z}_{+}$ and all $\lambda\in {\Bbb R}$, we have the representation
  \[
 \frac{d \la E (\lambda)e_n, e_m\ra } {d\lambda}= (2\pi)^{-1}\sqrt{1- \beta^2_{\infty}}  |\Omega (\lambda\pm i0) |^{-2} P_{n} (\lambda) P_m (\lambda)   ,\q \lambda\in \Bbb R.
\]
 In particular,
the spectral measure of the operator $J$ equals
  \begin{equation}
d\rho(\lambda) : =d \la E (\lambda)e_0, e_0 \ra= \tau(\lambda) d\lambda 
\label{eq:SFx}\end{equation}
where the weight  $\tau(\lambda)$ is given by the formula  
 \begin{equation}
\tau(\lambda)=  (2\pi)^{-1}\sqrt{1- \beta^2_{\infty}}  \, |  \Omega ( \lambda\pm i0) |^{-2}   .  
\label{eq:SF1}\end{equation}
 \end{theorem} 
 
 Putting together Theorem~\ref{HX+} and formula \e{eq:SF1}, we obtain
 
  \begin{theorem}\label{SFr}
  Under the assumptions of Theorem~\ref{AC}
 the weight $\tau (\lambda)$ is a continuous strictly positive function of $\lambda\in \Bbb R$.
 \end{theorem}

 Note that this result was earlier obtained in \cite{Apt} and \cite{Jan-Nab}  by specific methods of the orthogonal polynomials theory.


 In view of   \e{eq:SF1} the amplitude $\kappa(\lambda)$ defined by \e{eq:AP+} can be  expressed via the weight $\tau(\lambda)$, 
     \[
 \kappa (\lambda) =  (2\pi)^{-1/2}  (1- \beta^2_{\infty})^{1/4} \tau(\lambda)^{-1/2}  ,
\]
and hence asymptotic formula   \e{eq:Sz+} can be rewritten as
    \begin{equation}
 P_{n} (\lambda)=  (2\pi)^{-1/2}  (1- \beta^2_{\infty})^{-1/4} \tau(\lambda)^{-1/2}   a_{n}^{-1/2} \sin (\varphi_{n}  (\lambda) + \eta(\lambda) ) +  o(1)   
\label{eq:SW}\end{equation}
as $n\to\infty$.  This form seems to be more common for the orthogonal polynomials literature (cf. Theorem~3 in \cite{Apt}).


   \section{Asymptotics in the complex plane} 
   
   In this section, we  find (growing) asymptotics of the orthonormal  polynomials $P_{n} (z)$ as $n\to\infty$ for $\Im z\neq 0$.  To that end, we first solve the same problem for the solutions $g_{n}  (z)$ of the Jacobi equation \e{eq:Jy}  defined by formulas 
 \e{eq:GE}, \e{eq:GEx}.  Since $P_{n} (z)$ are linear combinations of the solutions $f_{n}  (z)$  and $g_{n}  (z)$, this yields an asymptotics of $P_{n} (z)$. We start in Sect.~6.1 with some general arguments which apply to all Jacobi coefficients satisfying condition \e{eq:Haupt}. Then we consider the cases of  large (Sect.~6.2) and small (Sect.~6.3)     diagonal elements $b_{n}$ separately. 
 
        \subsection{A representation for growing solutions} 
       
        Recall that the Jost solution $f_{n}  (z)$ of   equation \e{eq:Jy}  was defined in Theorem~\ref{GSS} by the formula
        \begin{equation}
  f_{n} (z) = a_{n}^{-1/2} q_{n} (z) u_{n} (z)  
\label{eq:JostG}\end{equation}
 where  $q_{n} $ is   product \e{eq:re1} and
  $u_{n}\to 1$ as $n\to\infty$ according to Theorem~\ref{GS3}. 

 We define the second solution $g_{n}  (z)$ of \e{eq:Jy} by equalities \e{eq:GE},  \e{eq:GEx}   which in view of    \e{eq:JostG} yields a representation
  \begin{equation}
  g_{n} (z) = a_{n}^{-1/2} q_{n} (z) u_{n} (z)   F_{n} (z)
\label{eq:TS}\end{equation}
 where
  \begin{equation}
  F_{n} (z) =   \sum_{m=N_{0}+1}^n \sqrt{\frac{a_{m}}{a_{m-1}}}  (q_{m-1} (z)q_{m} (z))^{-1}   (u_{m-1} (z)u_{m} (z))^{-1}  .
\label{eq:TS1}\end{equation}

Our goal is to find   an  asymptotics of this sum. 
  First, integrating  by parts we transform representation  \e{eq:TS1}   to a more convenient form.
    Recall that the numbers $\alpha_{n}$, $\beta_{n} $ were defined by equalities \e{eq:aabb} and the numbers $z_{n}$, $\z_{n} $ were defined by equalities \e{eq:zz}.

  \begin{lemma}\label{GM1}
    Let us set
\begin{equation}
 v_{n}= \varkappa_{n-1} \z_{n-1}\z_{n} (1- \z_{n-1}\z_{n} )^{-1}(    u_{n-1} u_{n})^{-1}  
\label{eq:hv}\end{equation}
where the numbers $\varkappa_n$ are given by formula \e{eq:Gr4K} and
\begin{equation}
  t_{n}=(q_{n-1} q_{n})^{-1} =\z_{0}^{-2} \z_{1}^{-2} \cdots \z_{n-2}^{-2}\z_{n-1}^{-1}.
\label{eq:hv2}\end{equation}
Then  sum \e{eq:TS1} can be written as 
 \begin{equation}
F_{n} =   v_{n}  t_{n+1} - v_{N_{0}}  t_{N_{0}+1} +\wt{F_{n}}  \label{eq:GE1}\end{equation}
where
\begin{equation}
 \wt{F_{n}}
=-   \sum_{m= N_{0} + 1}^n v_{m-1}'   t_{m}.
\label{eq:GE1a}\end{equation}
  \end{lemma}

   \begin{pf}
   Calculating the derivative of  product \e{eq:hv2}, we see that
   \[
      t_{n}' = ( \z_{n-1}\z_{n})^{-1} (1- \z_{n-1}\z_{n} )  t_{n}, 
\]
   whence using notation \e{eq:hv}, we obtain
   \[
    \sqrt{\frac{a_{n}}{a_{n-1}}}  (q_{n-1} q_{n}  )^{-1}   (u_{n-1} u_{n})^{-1}= ( \z_{n-1}\z_{n})^{-1} (1- \z_{n-1}\z_{n} ) v_{n} t_{n}=  v_{n}  t_{n}' .
    \]
    We can now rewrite representation  \e{eq:TS1}  as
    \[
F_{n} =  \sum_{m= N_{0} + 1}^n v_{m}   t_{m}'.
\]
Applying here integration-by-parts  formula \e{eq:Abel}, we arrive at relations \e{eq:GE1}, \e{eq:GE1a}.
    \end{pf}

       We will see that an asymptotics of the sequence $F_{n}$ is determined by  the first term in the right-hand side of \e{eq:GE1}. Let us calculate it.  Recall that  $\z_{\infty}$ is given by \e{eq:aabb3}  for  all $z\in{\Bbb C}$ with $\Im z \neq 0$  if $|\beta_{\infty}| > 1$.  If  $|\beta_{\infty}| < 1$, then $\z_{\infty}=\z_{\infty}^{(\pm)}$ is given by \e{eq:aabb4}  for  $ \pm \Im z > 0$.  The number $\varkappa_\infty$ is defined as limit \e{eq:Gr6a}. 
       
        \begin{lemma}\label{GM2+}
     The asymptotic relation
  \begin{equation}
\lim_{n\to\infty}    q_{n}^2 (z)    v_{n}(z)  t_{n+1} (z) = \varkappa_\infty \z_{\infty}  (1-\z_{\infty}^2)^{-1}  
\label{eq:GM1+}\end{equation}
holds. 
 \end{lemma}
 
   \begin{pf}
Since $q_{n+1}=\z_{n} q_{n}$, it follows from equalities \e{eq:hv} and  \e{eq:hv2} that
\[
  q_{n}^2  v_{n}  t_{n+1} = q_{n}    q_{n+1}^{-1}  v_{n} = \varkappa_{n-1} \z_{n-1} (1-\z_{n-1}\z_{n})^{-1}(   u_{n-1} u_{n})^{-1}.
\]
This yields   \e{eq:GM1+} because  $\z_{n}\to \z_{\infty} $, $\varkappa_n\to \varkappa_\infty$  and  $u_n\to 1$   as $n\to\infty$.  
      \end{pf}
      
      It remains to show that the second and third terms in the right-hand side of  \e{eq:GE1}  give no contribution to the asymptotics of $F_{n}$. This is obvious for
        $v_{N_{0}}  t_{N_{0}+1}$    because $q_{n} (z)\to 0$ as $n\to\infty$. 
        
        A relatively difficult part of the proof is
    to show that the same is true for the sum $\wt{F_{n}}$. This requires an inclusion
      \begin{equation}
( u_{n}' ) \in \ell^1 ( {\Bbb Z}_{+})
 \label{eq:TT1}\end{equation}
      for the  derivatives $u_{n}'$ of $u_{n}$.  
By its proof, we have to distinguish the cases $|\beta_{\infty}| >1$ and $|\beta_{\infty}| <1$.  An important  difference between   them  is that  $|\z_{\infty} |<1$ for  $|\beta_{\infty}| >1$ while $|\z_{\infty}^{(\pm)}|=1$ 
 for  $|\beta_{\infty}| <1$.  This essentially simplifies estimates in the first case.

      \subsection{Large diagonal elements} 
      
            Here we suppose that assumption \e{eq:Haupt} is satisfied with $|\beta_{\infty}| >1$ and consider arbitrary $z\in{\Bbb C}$. 
   

 \begin{lemma}\label{TT+}
Inclusion
     \e{eq:TT1} is true if the assumptions of Theorem~\ref{GSS} are satisfied.
\end{lemma}

 \begin{pf}
 It follows from identity \e{eq:A17Mb} that
    \begin{equation}
|  u_{n}'|\leq C \sum_{m=n+1}^\infty   | S_{n+1}^{-1} S_{m +1}|  |r_{m}|.
 \label{eq:TT1z}\end{equation}
 Since
       \begin{equation}
|\z_{n}|\leq \varepsilon< 1
 \label{eq:TT1y}\end{equation}
  for sufficiently large $n$, by definition \e{eq:Gp2a}, we have
  \[
  | S_{n+1}^{-1} S_{m +1}| =| \sigma_{n+1} \cdots \sigma_{m}| \leq C \varepsilon^{2(m-n)}, \q m> n,
  \]
  so that according to \e{eq:TT1z}
   \[
\sum_{n=0}^\infty |  u_{n}'|\leq C \sum_{n=0}^\infty   \sum_{m=n+1}^\infty \varepsilon^{2(m-n)} |r_{m}|
=C\sum_{m=1}^\infty \big(  \sum_{n=0}^{m-1} \varepsilon^{2(m-n)} \big) |r_{m}|
\leq C_{1}\sum_{m=1}^\infty  |r_{m}|.
\] 
Since $( r_{n})\in \ell^1 ( {\Bbb Z}_{+})$, this implies inclusion \e{eq:TT1}.   
   \end{pf}
   
   Now we are in a position to estimate  $\wt{F_{n}}(z)$.   
   
 \begin{lemma}\label{GM3+}
 Let  $\wt{F_{n}}(z)$ be defined by formula \e{eq:GE1a}. Then
     \begin{equation}
\lim_{n\to\infty}   q_{n}^2(z)  \wt{F_{n}}(z)=  0.
\label{eq:gasW}\end{equation}
\end{lemma}

  \begin{pf} 
  By definition \e{eq:re1} of $q_{n}$, we have
    \begin{equation}
 q_{n}^2  \wt{F_{n}}
=-   \sum_{m= N_{0} + 1}^n (\z_{m}^{2}\cdots \z_{n-1}^{2}) (\z_{m-1}v_{m-1}' )
\label{eq:TR}\end{equation}
where  the sequence $v_{n}$ is defined by relation  \e{eq:hv}. Note  that $\varkappa_{n}' \in \ell^1$ according to Remark~\ref{Grem1}  and $\z_{n}' \in \ell^1$ according to  conditions   \e{eq:Gr8}. Therefore it follows from  Lemma~\ref{TT+} that
   \begin{equation}
| v_{n}'| \leq C  | u_{n}'|\q \mbox{whence} \q v_{n}'\in \ell^1 ({\Bbb Z}_{+}).
\label{eq:TT2}\end{equation}
 Using now \e{eq:TT1y}, we find  that
 \begin{align*}
| q_{n}^2  \wt{F_{n}}|
&\leq  \sum_{m= N_{0} + 1}^{n/2} \varepsilon^{2 (n-m)}|v_{m-1}' | + \sum_ {m\geq n/2}^n \varepsilon^{2 (n-m)} |v_{m-1}' |\nonumber\\
&\leq \varepsilon^{n}  \sum_{m\geq N_{0} + 1} |v_{m-1}' | + \sum_ {m\geq n/2}   |v_{m-1}' |.
\end{align*}
Inclusion  \e{eq:TT2}  implies that
both  terms on the right tend to zero as $n\to\infty$.
 \end{pf} 
 
 Putting now together equality \e{eq:GE1}  with Lemmas~\ref{GM2+} and \ref{GM3+} and using that
 \[
 \z_{\infty}  (1-\z_{\infty}^2)^{-1}=\frac{ 1}{2\sqrt{\beta^2-1}}, 
 \]
 we find an asymptotics of the sequence $g_{n} (z)$.
 
  \begin{theorem}\label{GEL}
  Let $|\beta_{\infty}| >1$, and let $\varkappa_\infty$ be defined by \e{eq:Gr6a}.  
  Under the assumptions of Theorem~\ref{GSS}  the relation 
   \[
\lim_{n\to\infty} \sqrt{a_{n}}  q_{n}(z)  g_{n}(z)=\frac{ \varkappa_\infty}{2  \sqrt{\beta_{\infty}^2-1}}
\]
is true for all  $z\in{\Bbb C} $ with convergence uniform on compact subsets of ${\Bbb C} $.
 \end{theorem}
 
 Let us now use equality \e{eq:Pf}. In view of Theorem~\ref{GSS} the term $\omega(z) f_{n} (z)$ is negligible unless $\Omega (z)=0$.  Therefore Theorem~\ref{GEL}  yields an asymptotics of the orthonormal polynomials.

 \begin{theorem}\label{GEL1}
   Let $|\beta_{\infty}| >1$. 
  Under the assumptions of Theorem~\ref{GSS}  the relation 
   \[
\lim_{n\to\infty} \sqrt{a_{n}}  q_{n}(z)  P_{n}(z)=-\frac{\varkappa_\infty \Omega(z)}{2 \sqrt{\beta_{\infty}^2-1}}
\]
is true for all  $z\in{\Bbb C} $ with convergence uniform on compact subsets of ${\Bbb C}\setminus {\Bbb R}$.
 Moreover, if  $ \Omega(z)=0$, then 
   \[
\lim_{n\to\infty} \sqrt{a_{n}}  q_{n}(z)^{-1}  P_{n}(z)=  \{P(z), g(z)\}
\]
$($note that   $ \{P(z), g(z)\}\neq 0$ if  $ \Omega(z)=0)$.
  \end{theorem}

      \subsection{Small diagonal elements} 
      
      Here we suppose that assumption \e{eq:Haupt} is satisfied with $|\beta_{\infty}| <1$  and $z\in{\Bbb C}\setminus\Pi$.   In contrast to the previous subsection, the Carleman condition   \e{eq:Carl}  is also assumed. We again use identity \e{eq:A17Mb} but integrate by parts in its right-hand side which requires additional assumptions. 
      First, we improve estimate  \e{eq:Gr9} of the remainder $ r_{n} (z)$. Recall notations \e{eq:aabb} and \e{eq:Gr4K}.
      
      \begin{lemma}\label{TY}
   Suppose that, for some $\d >0$, 
    \begin{equation}
  \alpha'_{n} =O (n^{-1-\d}),\;  \alpha''_{n} =O (n^{-2-\d}) \q
\beta'_{n} =O (n^{-1-\d}),\; \beta''_{n} =O (n^{-2-\d})
 \label{eq:asb}\end{equation}
 and
   \begin{equation}
 k_{n}-1=O (n^{-1-\d}),\q    k_{n}'    =O (n^{-2-\d}).
 \label{eq:ask}\end{equation}
Then  remainder \e{eq:Gr6}   satisfies  estimates
    \begin{equation}
  r_{n} (z) =O (n^{-1-\d})\q \mbox{and}\q    r_{n}' (z)= O (n^{-2-\d}).
 \label{eq:as}\end{equation}
 \end{lemma}
 
  \begin{pf}
  Using \e{eq:Gs}, we rewrite
 relation \e{eq:Gr6}  as
   \begin{equation}
r _{n} = - z_{n-1}'  \varrho _{n}  + ( k_{n}-1)\z_{n}
\label{eq:Gs+}\end{equation}
where $z_{n}$,  $\z_{n}$ and $\varrho _{n} $ are defined by equalities \e{eq:aabb1} and \e{eq:Gs+1}.
  It follows from assumptions \e{eq:asb}   that
 \[
   z'_{n} =O (n^{-1-\d}),\q    z''_{n} =O (n^{-2-\d}).
 \]
 In view of \e{eq:aabb5} this also implies that
 \[
 \varrho_{n}= O (1) \q \mbox{and}\q
 \varrho_{n}' = O (n^{-1-\d}). 
 \]
 Thus, the first estimate \e{eq:as} is a direct consequence of relation \e{eq:Gs+}.  
 
 Using \e{eq:diff1}, we differentiate \e{eq:Gs+}. Since $   |\z_{n}'| \leq C | z_{n}'| $ according to  \e{eq:Gs},  the second bound \e{eq:as} is also a consequence of the estimates obtained.
      \end{pf}
 
 An asymptotics of $g_{n}(z)$ is stated in the following assertion.

 \begin{theorem}\label{GEe}
 Let $|\beta_{\infty}| <1$. 
 Let assumptions \e{eq:Haupt}, \e{eq:Carl} as well as  \e{eq:asb} and  \e{eq:ask} be satisfied.
  Then
     \begin{equation}
\lim_{n\to\infty}  \sqrt{a_{n}}   q_{n} (z)  g_{n}(z)  =\pm \frac{1}{2i \sqrt{1 - \beta^2_{\infty} } }, \q \pm \Im z>0. 
\label{eq:gas+}\end{equation}
 \end{theorem}

According to representation  \e{eq:TS} where $u_{n}(z)\to 1$  as $n\to\infty$
  relation \e{eq:gas+} can be rewritten as
   \begin{equation}
\lim_{n\to\infty}     q_{n} (z)^2  F_{n}(z)  =\pm \frac{1}{2i \sqrt{1 - \beta^2_{\infty} } }, \q \pm \Im z>0. 
\label{eq:gas++}\end{equation}

For a proof of \e{eq:gas++}, we again proceed from equality   \e{eq:GE1}
  and consider the terms on the right separately. According to Lemma~\ref{GM2+} the asymptotics of $q_{n}^2 v_{n }t_{n+1}$ is   given by relation
\e{eq:GM1+}  where according to  \e{eq:aabb4} we now have
   \begin{equation}
 \z_{\infty}^{(\pm)}  \big(1-(\z_{\infty}^{(\pm)} )^2\big)^{-1}=\pm \frac{ 1}{2i\sqrt{1-\beta^2_{\infty}}}, 
 \label{eq:ga}\end{equation}
 The term $v_{N_{0}}  t_{N_{0}+1}$ in \e{eq:GE1}  gives no contribution to the asymptotics of $F_{n}$. 
 
It remains to show that the same is true for the sum $\wt{F_{n}}$.  
We start with estimates of derivatives $u_{n}'$. This requires bootstrap arguments.

 \begin{lemma}\label{TT}
 Let the  assumptions of Theorem~\ref{GEe}  be satisfied.
 Then inclusion
     \e{eq:TT1} holds.
\end{lemma}

 \begin{pf}
Let us proceed from identity \e{eq:A17Mb} and integrate by parts in its right-hand side.  It follows from definition  \e{eq:Gp2a} that  
\[
S_{n}'= S_{n}(\sigma_{n}-1).
\]
Setting  $\ti{r}_{m}= (\sigma_{m+1}-1)^{-1} \z _{m} ^{-1}  r_{m}$, we see that
     \begin{align}
   \varkappa_{n}^{-1}   u'_{n} = &-  S_{n+1}^{-1} \sum_{m=n+1}^\infty        S'_{m+1} \ti{r}_{m} u_{m} 
 \nonumber  \\ 
  = &\ti{r}_{n-1} u_{n-1}-    \sum_{m=n+1}^\infty      S_{n+1}^{-1} S_{m+1} \ti{r}_{m-1}' u_{m}-  \sum_{m=n+1}^\infty      S_{n+1}^{-1} S_{m+1} \ti{r}_{m-1} u_{m-1}'. 
 \label{eq:TT1x}\end{align}
 
 Under our assumptions we have  estimates  \e{eq:as}.  Since $u_{n}\in \ell^\infty ( {\Bbb Z}_{+})$, the first term in the right-hand side of \e{eq:TT1x} is   bounded by $n^{-1-\d}$.
 By definition  \e{eq:Gp2a}, we have $| S_{n+1}^{-1} S_{m+1}|\leq 1$ if $m\geq n+1$. Therefore in view of   estimate \e{eq:as} on $\ti{r}_{n}' $,
 the second term in the right-hand side of \e{eq:TT1x} is also bounded by $n^{-1-\d}$. Thus idenity \e{eq:TT1x} implies that
  \begin{equation}
 |u_{n}'| \leq C n^{-1-\d}+ C \sum_{m=n}^\infty   m^{-1-\d}    | u_{m}'|.
 \label{eq:as1}\end{equation}
 
 Suppose now that
$  u_{n}'= O( n^{-\sigma})  $
  for some $\sigma \in [0,1]$. Then it follows from \e{eq:as1} that actually $  u_{n}'= O( n^{-\sigma-\d})  $. Since 
  $  u_{n}'= O( 1)  $, repeating this argument a sufficient number of times,
  we find that $  u_{n}'= O( n^{-1-\d})  $. This implies \e{eq:TT1}.
         \end{pf}


Now we are in a position to estimate sum \e{eq:GE1a}. As in the case $ |\beta_{\infty}| >1$, we essentially use inclusion \e{eq:TT1}.

 \begin{lemma}\label{GM3}
 The function  $\wt{F_{n}}(z)$ defined by formula  satisfies relation
     \e{eq:gasW}.
\end{lemma}

  \begin{pf} 
  By definition \e{eq:re1}, we  have relation 
    \e{eq:TR} where $v_{n}$ are given by \e{eq:hv}.   Lemma~\ref{TT} again  implies inclusion
   \e{eq:TT2}. Compared to the case $|\beta_{\infty}| >1$, the  problem is that now inequality \e{eq:TT1y} is no longer true.  Instead, we have a weaker inequality
    \begin{equation}
|\z_{n}|\leq e^{-c \alpha_{n}}, \q c=c(z)>0,
\label{eq:ca}\end{equation}
following from Corollary~\ref{Te1c}.
   
   Let us now write \e{eq:TR} as
 \begin{equation}
 q_{n}^2  \wt{F_{n}}
=-   \sum_{m= N_{0} + 1}^\infty X_{m,n}  \z_{m-1}v_{m-1}' )
\label{eq:TR1}\end{equation}
where 
 \[
  X_{m,n}  = \z_{m}^{2}\cdots \z_{n-1}^{2}\q \mbox{for}\q   m<n  \q \mbox{and}\q X_{m,n}   = 0 \q \mbox{for}\q m\geq n .
\]
According to \e{eq:ca} we have an estimate
\[ 
| X_{m,n}  |\leq \exp \big(-c \sum_{k=m}^{n-1}\alpha_{k}\big).
\]
 Therefore $  X_{m,n} \to 0$ as $n\to\infty$ for every fixed $m$ by virtue of the  Carleman condition \e{eq:Carl}. Now inclusion  \e{eq:TT2} allows us to use the dominated convergence theorem to show that the limit of sum \e{eq:TR1} as $n\to\infty$ is zero.
     \end{pf}
     
     Let us come back to the proof of Theorem~\ref{GEe}. As in Sect.~6.2, we proceed from relation \e{eq:GE1}.  Combining Lemma~\ref{GM2+} and
       equality  \e{eq:ga}, we find that
       \[
\lim_{n\to\infty}    q_{n}^2 (z)    v_{n}(z)  t_{n+1} (z) = \pm \frac{1}{2 i \sqrt{1-\beta^2_{\infty}}}.
\]
       Therefore using Lemma~\ref{GM3}, we arrive at relation  \e{eq:gas++} 
     or, equivalently, \e{eq:gas+}. $\q \Box$

In view of equality \e{eq:Pf} the following result is a direct corollary
of Theorems~\ref{GSS}   and \ref{GEe}. We recall that the Wronskian $\Omega(z)$ is defined by relation   \e{eq:J-W} and $\Omega(z)\neq 0$ for $\Im z\neq 0$.

\begin{theorem}\label{GE1+}
  Under the assumptions of Theorem~\ref{GEe}    the relation 
   \[
\lim_{n\to\infty} \sqrt{a_{n}}  q_{n}(z)  P_{n}(z)=\frac{i \Omega(z)}{2\sqrt{1-\beta^2_{\infty}}}
\]
is true for all  $z\in{\Bbb C}\setminus \Bbb R$  
with convergence uniform on compact subsets of ${\Bbb C}\setminus {\Bbb R}$.
 \end{theorem}
 
               \section{Non-Carleman case} 
               
               In this section, we specially consider the case when off-diagonal entries $a_{n} $ grow so rapidly that
                \begin{equation}
\sum_{n=0}^\infty a_{n}^{-1}<\infty;
   \label{eq:nonCarl}\end{equation}
   thus the Carleman condition  \e{eq:Carl}  fails. Asymptotic properties of orthonormal polynomials are discussed in Sect.~7.1 and spectral properties of the corresponding Jacobi operators - in the following subsections. Some proofs will be omitted since they were already published in the papers \cite{nCarl}  and \cite{Jacobi-LC}.
 

     \subsection{Jost solutions and orthogonal polynomials}
     
     In the case $|\beta_{\infty} |>1$, there is almost no difference between the Carleman and non-Carleman cases.
     Thus, Theorems~\ref{GSS}, \ref{GEL}  and \ref{GEL1}   remain  true, and their formulations may even be  simplified. Indeed, the sequence $q_{n}  (z)$ in all asymptotic formulas may be replaced (cf. Sect.~4.3) by a simpler sequence $\Theta_{n}$  defined by formula     \e{eq:ss5o}. We emphasize that $\Theta_{n}$ do not depend on $z\in{\Bbb C}$. Moreover,  in the  non-Carleman case, the first condition \e{eq:Gr8} can be omitted. Actually,  in view of Remark~\ref{Grem1} under assumption  \e{eq:Gr6b}  conditions \e{eq:nonCarl} and   
      \[
\sum_{n=0}^\infty \alpha_{n} <\infty
   \]
   are equivalent. The last condition implies that the series of $|\alpha_{n}' |$ converges.
   
   Let us summarize the results of Theorems~\ref{GSS}, \ref{GEL}  and \ref{GEL1} 
     for the non-Carleman case.  Recall that the sequence $g_{n} (z)$ is defined by formulas \e{eq:GE},  \e{eq:GEx} and $\varkappa_\infty$ -- by formulas \e{eq:Gr4K},  \e{eq:Gr6a}.
   
   \begin{theorem}\label{GSnC1}
Let condition \e{eq:nonCarl} hold true and $|\beta_{\infty} |>1$.
  Suppose that  assumptions  \e{eq:Haupt},  \e{eq:Gr6b}  and
  \begin{equation}
   (\beta_{n}' )\in \ell^1 ({\Bbb Z}_{+} )
\label{eq:Gr8B}\end{equation}
are satisfied.
        Then for   all $z\in{\Bbb C}$, equation \e{eq:Jy} has a unique solution $f_{n}    (  z ) $ with asymptotics
  \begin{equation}
f_{n}    (  z )= a_{n} ^{-1/2} \Theta_{n}   \big(1 + o( 1)\big), \q n\to\infty. 
\label{eq:A22g-}\end{equation} 
Besides,   the relations 
   \begin{equation}
\lim_{n\to\infty} \sqrt{a_{n}}  \Theta_{n} g_{n}(z)=\frac{ \varkappa_\infty}{2  \sqrt{\beta_{\infty}^2-1}}
\label{eq:GEGnC}\end{equation}
and
 \[
\lim_{n\to\infty} \sqrt{a_{n}}  \Theta_{n}    P_{n}(z)=-\frac{\varkappa_\infty \Omega(z)}{2 \sqrt{\beta_{\infty}^2-1}}
\]
are true for all  $z\in{\Bbb C} $ with convergence uniform on compact subsets of ${\Bbb C} $.
Moreover, if  $ \Omega(z)=0$, then 
   \[
\lim_{n\to\infty} \sqrt{a_{n}}  \Theta_{n}^{-1}  P_{n}(z)=  \{P(z), g(z)\}
\]
$($note that   $ \{P(z), g(z)\}\neq 0$ if  $ \Omega(z)=0)$.
All functions $f_{n}   ( z )$    are analytic in $z\in  \Bbb C$.
 \end{theorem}
 
 On the contrary, for $|\beta_{\infty} | < 1$, the Carleman and non-Carleman cases are significantly different. In particular,  for all $z\in{\Bbb C}$, we now have two linearly independent solutions $f_{n}^{(+)}   (  z )$  and $f_{n}^{(-)}   (  z )$. The leading terms of their asymptotics do not depend on $z$  and are complex conjugate to each other.


   \begin{theorem}\label{GSnC}
Let condition \e{eq:nonCarl} hold true and $|\beta_{\infty} |<1$.
  Suppose that  assumptions  \e{eq:Haupt},  \e{eq:Gr6b}  and
  \e{eq:Gr8B} are satisfied. Define the sequences  $\Theta_{n}^{(\pm)}  $ by formula
       \e{eq:ss6o}.
   Then for each of the signs $``\pm"$  and all $z\in{\Bbb C}$, equation \e{eq:Jy} has a unique solution $f_{n}^{(\pm)}   (  z )$ with asymptotics
  \begin{equation}
f_{n}^{(\pm)}   (  z )= a_{n} ^{-1/2} \Theta_{n}^{(\pm)}    \big(1 + o( 1)\big), \q n\to\infty. 
\label{eq:A22g+}\end{equation} 
These solutions are linked by the relation  
\begin{equation}
f^{(-)} (\bar{z}) =\ov{ f^{(+)} (z)}
\label{eq:WR1}\end{equation}  
and their Wronskian equals
\begin{equation}
\{f^{(+)} (z), f^{(-)} (z)\}= - 2i \varkappa_\infty^{-1} \sqrt{1-\beta_{\infty}^2}\neq 0
\label{eq:WR}\end{equation} 
so that the    solutions $f_{n}^{(+)}   (  z )$ and $f_{n}^{(-)}   (  z )$ are   linearly independent.   
 For all $n\in {\Bbb Z}_{+}$, the functions   $f_{n}^{(\pm)}   ( z )$   are analytic in $z\in  \Bbb C$.
  \end{theorem}


Let us make some comments on the {\it  proof}.  Under assumption  \e{eq:nonCarl} we can   define the Ansatz   by one of the formulas $ {\cal A}^{(\pm)} = a_{n} ^{-1/2} \Theta_{n}^{(\pm)}   $. Note that it does not depend on $z$.  An easy calculation shows that, for this Ansatz,  remainder  \e{eq:Gry}  equals
\[
 r_{n}^{(\pm)}    (z)  =   \big(  \z(\beta_{n-1}  \pm i 0) ^{-1}-\z(\beta_{n}   \pm i 0)^{-1}\big)+ (k_{n}-1)\z(\beta_{n}  \pm i 0) -2z\alpha_{n}.
\]
Then $( r_{n}^{(\pm)}    (z))\in \ell^1$, and the whole scheme of the proof of Theorem~\ref{GSS} works without any changes. Moreover, all iterations \e{eq:W5} are entire functions of $z$ because $ r_n^{(\pm)} (z)$ is a linear function of $z$ and the kernels  $G_{n,m}$ do not depend on $z$. This ensures that the solutions $u_{n}^{(\pm)}( z )$ of Volterra equation   \e{eq:A17}  and hence  functions \e{eq:Jost} are analytic in $z\in  \Bbb C$.

  It follows from formulas   \e{eq:ss6o} and \e{eq:A22g+} that the Wronskian of the solutions $f^{(+)} (z)$ and $ f^{(-)} (z)$  equals
\begin{align*}
\{f^{(+)} (z), & f^{(-)} (z)\}=\lim_{n\to\infty}\varkappa_{n}^{-1} (\Theta_{n}^{(+)}\Theta_{n+1}^{(-)} -\Theta_{n+1}^{(+)}
\Theta_{n}^{(-)}) 
\nonumber\\
= &\lim_{n\to\infty}\varkappa_{n}^{-1} ( e^{- i\arccos \beta_{n}}-e^{i\arccos \beta_{n}})
 =  - 2i \varkappa_\infty^{-1} \sin \arccos \beta_\infty 
\end{align*} 
which yields  relation  \e{eq:WR}.  $\q\q\q \Box$

All solutions of equations  \e{eq:Jy}   are linear combinations of the Jost solutions $f^{(+)}   (  z )$ and $f ^{(-)}   (  z )$, and hence their asymptotics are determined by Theorem~\ref{GSnC}.  In particular, this is true for the polynomials $P_{n} (z)$ defined by   conditions \e{eq:RR1} and the polynomials of the second kind $Q_{n} (z)$ defined by the boundary conditions
 \begin{equation}
 Q_{0} (z)=0, \q Q_{1} (z)=a_{0}^{-1}.
\label{eq:RQ}\end{equation} 
Set $P(z)= (P_{n}(z))$,  $Q(z)= (Q_{n}(z))$.
Note that
 \begin{equation}
  {\cal J} P(z) =z P(z) \q \mbox{and} \q   {\cal J} Q(z) =z Q(z)  +e_{0}.
  \label{eq:PQzz}\end{equation}

    Thus, we have 
  \begin{equation}
P_{n}(z)= \sigma_{+} (z)f_{n}^{(+)}(  z ) + \sigma_{-} (z) f_{n}^{(-)}(  z )
\label{eq:LC}\end{equation}
and
  \begin{equation}
Q_{n}(z)= \tau_{+} (z)f_{n}^{(+)}(  z ) + \tau_{-} (z) f_{n}^{(-)}(  z ),
\label{eq:LCq}\end{equation}
where the coefficients $ \sigma_{\pm} (z)$ and $ \tau_{\pm} (z)$
 can be expressed via the Wronskians:
  \begin{equation}
  \sigma_{+}(z)=  \frac{ i \varkappa_{\infty} } {2 \sqrt {1-\beta^2_{\infty}}}  \{P(z), f^{(-)}(  z )\}, \q
  \sigma_{-}(z)=  - \frac{ i \varkappa_{\infty} } {2 \sqrt {1-\beta^2_{\infty}}}  \{P(z), f^{(+)}(  z )\}
\label{eq:st}\end{equation}
and 
  \begin{equation}
 \tau_{+}(z)=  \frac{ i \varkappa_{\infty} } {2 \sqrt {1-\beta^2_{\infty}}}  \{Q(z), f^{(-)}(  z )\}, \q
  \tau_{-}(z)=  - \frac{ i \varkappa_{\infty} } {2 \sqrt {1-\beta^2_{\infty}}}  \{Q(z), f^{(+)}(  z )\} .
\label{eq:ts}\end{equation}

Observe that
\[
\sigma_{-} (\bar{z})= \ov{\sigma_{+}(z)}\q \mbox{and}\q \tau_{-} (\bar{z})= \ov{\tau_{+}(z)}
\]
 because $P_{n}(\bar{z})= \ov{P_{n}(z)}$, $Q_{n}(\bar{z})= \ov{Q_{n}(z)}$  and $ f_{n}^{(\pm)}(  z )$ satisfy   \e{eq:WR1}. 
Of course, all coefficients $ \sigma_{\pm}(z)$ and $ \tau_{\pm}(z)$  are  entire functions of $z$.

According to \e{eq:LC}  and \e{eq:LCq}   the following result is a direct consequence of Theorem~\ref{GSnC}. 

 \begin{theorem}\label{LC}
    Under the assumptions of Theorem~\ref{GSnC} the orthogonal polynomials $P_{n}(z)$ and $Q_{n}(z)$
 have  the  asymptotics, as $n\to\infty$,
  \begin{equation}
P_{n} (z)= a_{n} ^{-1/2} \big(\sigma_{+} (z) \Theta_{n}^{(+)} + \sigma_{-}(z) \Theta_{n}^{(-)} + o(1)\big)      
\label{eq:LC1P}\end{equation} 
and
 \begin{equation}
Q_{n} (z)= a_{n} ^{-1/2} \big(\tau_{+} (z)\Theta_{n}^{(+)} + \tau_{-} (z) \Theta_{n}^{(-)}+ o(1)\big)   .
\label{eq:LC2q}\end{equation} 
   \end{theorem}

    In view of conditions \e{eq:RR1}  and \e{eq:RQ} the Wronskian $\{ P(  z ), Q(z)\} =1$. On the other hand, we can calculate this Wronskian using relations   \e{eq:WR} and \e{eq:LC}, \e{eq:LCq}. This yields an identity
   \begin{equation}
  2i \varkappa_{\infty}^{-1} \sqrt {1-\beta_{\infty}^2} \big(\sigma_{-} (z)\tau_{+} (z)-\sigma_{+} (z) \tau_{-} (z)\big)=1,\q \forall z\in{\Bbb C}.
\label{eq:Wro}\end{equation}
We also note the identity  
  \begin{equation}
|\sigma_{+} (z)|^{2} - |\sigma_{-} (z)|^{2}  =\Im z \: \varkappa_{\infty}  (1-\beta_{\infty}^2)^{-1/2}  \sum_{n=0}^{\infty}  |P_{n} (z)|^{2}
 \label{eq:LC2p}\end{equation} 
  established in Theorem~4.4 of \cite{nCarl}.

        \subsection{Essential self-adjointness}
   
   Here we consider self-adjoint extensions of the minimal Jacobi operator $J_{\rm min}$ defined in the space $\ell^2({\Bbb Z}_{+})$ by formula  \e{eq:ZP+1}  on the set ${\cal D}$ of elements with a finite number of non-zero components.
   
   Recall first of all that  the Carleman condition \e{eq:Carl} is sufficient but not necessary for the essential self-adjointness of the   operator $J_{\rm min}$. Nevertheless it is close to necessary for small diagonal elements $b_{n}$. Indeed, according to  the Berezanskii theorem (see, e.g., page 26
 in the book \cite{Ber}) if $b_{n}=0$ (or, more generally, $(b_{n})\in \ell^{\infty} ({\Bbb Z}_{+})$), then  the Carleman condition is equivalent to the essential self-adjointness of the   operator $J_{\rm min}$   provided     $a_{n-1}a_{n+1}\leq a_{n}^2$. 
     The following example shows that the last condition is very essential.

 \begin{example}\label{LCa}
Suppose that $b_{n}=0$ and that $a_{n}=n^p (1+c_{1}n^{-1})$ if $n$ is odd and $a_{n}=n^p (1+c_2n^{-1})$ if $n$ is even, at least
for  sufficiently large $n$.  
As shown in \cite{Kost, J-M}, the corresponding Jacobi  operator  $J_{\rm min}$ is essentially self-adjoint if $p>1$ and $|c_{2}-c_{1}| \geq p-1$.
Since
\begin{equation}
\frac{a_{n}}{\sqrt{a_{n-1} a_{n+1}}}=1+ (-1)^n  \,\frac{c_{2}-c_{1}}{n} +O (\frac{1}{n^2}), \q n\to\infty,
\label{eq:Ber}\end{equation} 
the condition $a_{n-1}a_{n+1}\leq a_{n}^2$ fails in this example.
 \end{example}
 
 In addition  to  the Berezanskii theorem, we have the following result where the case $b_{n}\to\infty$ is not excluded.
 
    \begin{proposition}\label{S-A+}
 Under the assumptions of Theorem~\ref{GSnC}   
  the operator $J_{\rm min}$   has deficiency indices $(1,1)$.
  \end{proposition}

  \begin{pf}
  Let us consider the solutions $f_{n}^{(\pm)} (z)$  of the Jacobi equation with asymptotics \e{eq:A22g+}. Since $a_{n}$ satisfy condition \e{eq:nonCarl}  and $|\Theta_{n}|=1$, it follows from  \e{eq:A22g+} that $f_{n}^{(\pm)} (z)\in \ell^2 ({\Bbb Z}_{+})$. The solutions   are linearly independent so that we are in the limit circle case. 
   \end{pf}
   
    Note that  this result does not apply to Example~\ref{LCa} because  relation \e{eq:Ber} excludes condition \e{eq:Gr6b}.
  
All self-adjoint extensions of the operator $J_{\rm min}$ will be described in Sect.~7.4.

Next, we consider the case  $|\beta_{\infty}| >1$.

    \begin{proposition}\label{S-A}
 Under the assumptions of Theorem~\ref{GSnC1}   the operator $J_{\rm min}$ is essentially self-adjoint if and only if
 \begin{equation}
\sum_{n=0}^\infty a_{n}^{-1}  \Theta_{n}^2=\infty.
\label{eq:S-A1}\end{equation}  
Otherwise, $J_{\rm min}$ has deficiency indices $(1,1)$.
  \end{proposition}
  
     \begin{pf}
       According to \e{eq:A22g-} under assumption \e{eq:nonCarl} the Jost solution $f_{n}(z)\in \ell^2({\Bbb Z}_{+})$ for all $z\in{\Bbb C}$.  Therefore the operator $J_{\rm min}$ is essentially self-adjoint if and only if the solution $g_{n}(z)$ linearly independent  with $f_{n}(z)$ is not in $\ell^2({\Bbb Z}_{+})$.  By virtue of  \e{eq:GEGnC} this is equivalent to condition \e{eq:S-A1}.
         \end{pf}

  Relation \e{eq:S-A1} generalizes the Carleman condition  \e{eq:Carl}. Proposition~\ref{S-A}  shows that  in the case  $|\beta_{\infty}| >1$  the operator $J_{\rm min}$ is essentially self-adjoint unless the sequence $a_{n}$ grows very rapidly. Indeed, it follows from definition \e{eq:psi+}  that
  \[
  \varphi_{n}= n \arccosh | \beta_{\infty}| + o(n) 
  \]
  as $n\to\infty$ whence
  \[
 \exp (n \arccosh | \beta_{\infty}| -\varepsilon n) \leq \Theta_{n}\leq \exp (n \arccosh | \beta_{\infty}| +\varepsilon n)
 \]
 for an arbitrary $\varepsilon>$ and sufficiently large $n$.  Therefore  the operator $J_{\rm min}$ is essentially self-adjoint if $a_{n} \leq c n^p$ for some $c> 0$, $p< \infty$ and all $n\geq 1$.  On the contrary,   $J_{\rm min}$ has deficiency indices $(1,1)$ if $a_{n}\geq c x^{n^p}$ for some $x>1$ and $p>1$.


      \subsection{Quasiresolvent}

In this subsection we do not make any  specific assumptions about the coefficients $a_{n}$ and $b_{n}$ supposing only  that the minimal Jacobi operator $J_{\rm min}$ is not essentially self-adjoint.    Thus, we are in the limit circle case so that 
\begin{equation}
\clos J_{\rm min}\neq J_{\rm max}=J_{\rm min}^*.
  \label{eq:cl}\end{equation}
  and the   inclusions 
  \begin{equation}
  P(z)\in  \ell^2 ({\Bbb Z}_{+}) \q \mbox{and} \q    Q (z)\in  \ell^2 ({\Bbb Z}_{+})  
  \label{eq:PQz}\end{equation}
  are satisfied
 for all $z\in {\Bbb C}$    (see, for example, Theorem~6.16 in \cite{Schm}).  
 
 In the limit circle case, the operator $J_{\rm max}$  is   of course not symmetric.  
     For all $u, v \in {\cal D} (J_{\rm max})$, we have the identity  (the Green formula)
       \begin{equation}
\la {\cal J} u,v \ra - \la u,  {\cal J} v \ra = \lim_{n\to\infty}  a_{n } (u_{n +1}\bar{v}_{n }- u_{n }\bar{v}_{n +1})
\label{eq:qres4}\end{equation}
where the limit in the right-hand side exists.
 Indeed, a direct calculation shows that
      \[
\sum_{m=0}^n ( {\cal J} u)_{m}\bar{v}_{m} -\sum_{m=0}^n u_{m} ( {\cal J} \bar{v})_{m} =     a_{n } (u_{n +1}\bar{v}_{n }- u_{n }\bar{v}_{n +1}).
\]
Passing here to the limit $n\to\infty$ and using that ${\cal J} u  \in \ell^2 ({\Bbb Z}_{+})$ for $u \in {\cal D} (J_{\rm max})$, we obtain \e{eq:qres4}.

    Let us define an operator ${\cal R}  (z)$ playing the role of the resolvent of the operator $J_{\rm max}$ by   an equality 
        \begin{equation}
 ( {\cal R} (z)h)_{n} =  Q_{n}(z) \sum_{m=0}^n  P_{m} (z) h_{m}+   P_{n}(z)   \sum_{m=n+1}^\infty  Q_{m} (z) h_{m}.\label{eq:RR11C}\end{equation} 
 In view of inclusions \e{eq:PQz} the operators ${\cal R}  (z)$ considered in the space $\ell^2 ({\Bbb Z}_{+})$ belong to the Hilbert-Schmidt class for all $z\in {\Bbb C}$.  Note also that ${\cal R}  (z)$ is a holomorphic operator-valued function of $z\in {\Bbb C}$ and ${\cal R} (z)^*={\cal R} (\bar{z})$.
   It follows from \e{eq:RR1}, \e{eq:RQ} that  
  \begin{equation}
  ({\cal R} (z)h)_{0}=\la h, Q(\bar{z})\ra 
  \label{eq:r01}\end{equation} 
and
 \begin{equation}
   ({\cal R} (z)h)_{1}= h_{0}  a_{0}^{-1}+  (z-b_{0})a_{0}^{-1} \la h, Q (\bar{z})\ra 
\label{eq:r02}\end{equation} 
for all  $h\in {\ell}^2  ({\Bbb Z}_{+})$.

Our proof of the following statement is  close to the construction  of the resolvents for    self-adjoint Jacobi operators (see Proposition~\ref{res}).    
   
    \begin{theorem}\label{resC}
    Let   \e{eq:cl} be satisfied.
 For all $z \in {\Bbb C}$, we have
  \begin{equation}
{\cal R}  (z): \ell^2 ({\Bbb Z}_{+})\to  {\cal D} (J_{\rm max} )
\label{eq:qres}\end{equation}
and
  \begin{equation}
(J_{\rm max} -zI) {\cal R}  (z)=I.
\label{eq:qres1}\end{equation}
   \end{theorem} 
 
  \begin{pf} 
  Recall that the operator ${\cal J}$  was defined by equalities \e{eq:ZP+1}. We will  check that 
    \begin{equation}
  (  ({\cal J}-z I) {\cal R} (z)h)_{n}=h_{n} 
    \label{eq:RRmc}\end{equation}
for all $n \in {\Bbb Z}_{+}$ and $h= ( h_{n})\in \ell^2 ({\Bbb Z}_{+})$.  For $n=0$, we have
\[
(  ({\cal J}-z I) {\cal R} (z)h)_0= (b_{0}-z)({\cal R} (z)h)_0  + a_{0}({\cal R} (z)h)_1=h_{0}
\]
according to formulas \e{eq:r01} and \e{eq:r02}. For $n\geq 1$, we rewrite definition  \e{eq:RR11}   as
 \[
 ( {\cal R} (z)h)_{n} =  Q_{n} (z) x_{n}  (z) +   P_{n} (z)   y_{n} (z)
 \]
where
   \begin{equation}
x_{n}(z) =\sum_{m=0}^n  P_{m} (z) h_{m},\q    y_{n}(z) =\sum_{m=n+1}^\infty  Q_{m} (z) h_{m}.
\label{eq:RR2c}\end{equation} 
  It now follows from definition \e{eq:ZP+1}     that
   \begin{multline}
 ( ({\cal J} -z  I) {\cal R} (z) h)_{n} =  a_{n-1}\big(Q_{n-1}x_{n-1}+P_{n-1} y_{n-1} \big)
  \\
  + (b_{n}-z) \big(Q_{n}x_{n}+ P_{n} y_{n}  \big)+ a_{n}\big(Q_{n+1}x_{n+1} +P_{n+1} y_{n+1}    \big), \q n\geq 1.
\label{eq:RR3c}
\end{multline}
 According to \e{eq:RR2c}  we have
\[
Q_{n-1}x_{n-1}+ P_{n-1} y_{n-1}= Q_{n-1}(x_{n}-P_{n} h_{n})+  P_{n-1}(y_{n}+Q_{n} h_{n}) 
\]
and
 \begin{align*}
Q_{n+1} x_{n+1} +P_{n+1} y_{n+1} & =
Q_{n+1}(x_{n}  +P_{n+1}h_{n+1}) 
\\
& +  P_{n+1} (y_{n} - Q_{n+1}h_{n+1})=
Q_{n+1}x_{n}  +  P_{n+1}y_{n} .
\end{align*}
 Substituting these expressions into the right-hand side of \e{eq:RR3c},  we see that
 \begin{align*}
  \big( ({\cal J} -z I) {\cal R} (z) h\big)_{n}& = a_{n-1}\Big(Q_{n-1}(x_{n}- P_{n} h_{n})+  P_{n-1}(y_{n}+Q_{n} h_{n}) \Big)
  \\
 & + (b_{n}-z) \big(Q_{n}x_{n}+ P_{n} y_{n}  \big)+ a_{n}\big(Q_{n+1}x_{n}  +  P_{n+1}y_{n}    \big).
\end{align*}
Let  us now collect together all terms containing $x_{n}$, $y_{n}$ and $h_{n}$.  Then
 \begin{multline}
  \big( ({\cal J} -z I) {\cal R} (z) h\big)_{n} = \Big(a_{n-1}Q_{n-1}+  (b_{n}-z) Q_{n} + a_{n} Q_{n+1}\Big) x_{n}
  \\
    + \Big(a_{n-1}P_{n-1}+  (b_{n}-z)P_{n} + a_{n}P_{n+1}\Big) y_{n}+a_{n-1}\big( -P_{n} Q_{n-1}+  P_{n-1} Q_{n}\big) h_{n}  .
\label{eq:RR3A}
\end{multline}
The coefficients at $x_{n}$ and $y_{n}$ equal zero by virtue of  the Jacobi equation \e{eq:Jy}  for $( Q_{n} )$ and $(P_{n})$, respectively. Since $\{P,Q \} =1$, the right-hand side of 
\e{eq:RR3A} equals $ h_{n}$.   
This proves \e{eq:RRmc}  whence $   ({\cal J}-z I) {\cal R} (z)h=h $. In particular, we see that ${\cal R}  (z) h\in {\cal D} (J_{\rm max})  $, and hence ${\cal J}$ can be replaced here by $J_{\rm max}$.
This yields both \e{eq:qres} and \e{eq:qres1}. 
    \end{pf}
    
    \begin{remark}\label{res3}
    In definition  \e{eq:RR11C},
    one can replace   $Q_{n}(z)$ by the polynomials $\wt{Q}_{n}(z)= Q_{n}(z)+ c P_{n}(z)$ for an arbitrary  $c\in {\Bbb C}$. Then $\wt{\cal R} (z)= {\cal R} (z) +c \la \cdot, P(\bar{z})\ra P(z)$ and formulas \e{eq:qres}, \e{eq:qres1}  remain true.
      \end{remark}

    Since $P(z)$ is a unique (up to a constant factor)  solution of the homogeneous equation $(J_{\rm max} -zI) u =0$, we can also state 
    
     \begin{corollary}\label{res1}
     Let $z\in {\Bbb C} $ and $ h\in \ell^2 ({\Bbb Z}_{+})$.
     Then all solutions of the equation
       \[
(J_{\rm max} -zI) u =h \q \mbox{where}  \q  
\]
  for $u\in {\cal D} (J_{\rm max} )$ are  given by the formula
 \begin{equation}
  u = \Gamma P(z)+ {\cal R}  (z) h \q \mbox{for some} \q \Gamma  \in{\Bbb C}.
\label{eq:qres3}\end{equation}
   \end{corollary}

The following  asymptotic relation for $ ( {\cal R} (z)h)_{n}$ is a direct consequence of definition  \e{eq:RR11C}
   and condition  \e{eq:PQz}:
      \begin{equation}
 ( {\cal R} (z)h)_{n} =  Q_{n}(z) \la    h, P (\bar{z}) \ra+   o(|P_{n} (z)| +|Q_{n}(z)| )\q {\rm as}\q n\to\infty
 \label{eq:Ras}\end{equation}
 for all $h\in \ell^2 ({\Bbb Z}_{+})$.
 This relation   can be supplemented by the following result.

     \begin{proposition}\label{AS}
     Set
   \begin{equation}
     {\sf D}=\clos \cal D (J_{\rm min}).
 \label{eq:MIN}\end{equation}
    For all $z\in {\Bbb C}$ and  all $h\in \ell^2 ({\Bbb Z}_{+})$, we have
    \[
 u:= {\cal R} (z)h - Q (z) \la    h, P (\bar{z}) \ra\in {\sf D}.
 \]
        \end{proposition} 
        
              \begin{pf}
              Let first  $h\in {\cal D}$. Then $ ( {\cal R} (z)h)_{n} =  Q_{n}(z) \la    h, P (\bar{z}) \ra$ for sufficiently large $n$, whence $u\in {\cal D}\subset {\sf D}$.

Let now $h$ be an arbitrary vector in  $  \ell^2 ({\Bbb Z}_{+})$.     Observe that $ u \in {\sf D}$ if and only if there exists a sequence $u^{(k)}\in {\cal D}$ such that 
                \begin{equation}
u^{(k)}\to u \q \mbox{and} \q {\cal J}u^{(k)}\to {\cal J}u \q \mbox{as} \q k\to\infty.
 \label{eq:Ras2}\end{equation}
            Let us take any sequence $h^{(k)}\in {\cal D}$ such that   $h^{(k)}\to h$ and set
                 \[
 u^{(k)}= {\cal R} (z)h^{(k)} - Q (z) \la    h^{(k)}, P (\bar{z}) \ra .
 \]
 Then $u^{(k)} \in {\cal D}$ and  $u^{(k)}\to u$ as $k\to\infty$ because the operator $ {\cal R} (z)$ is bounded.  It follows from equalities  \e{eq:PQzz}  and \e{eq:qres1}  that
 \[
({\cal J}  -z) u^{(k)}=  h^{(k)} - e_{0} \la    h^{(k)}, P (\bar{z}) \ra \to h - e_{0} \la    h , P (\bar{z} ) \ra =({\cal J}  -z) u 
\]
as $k\to\infty$. This proves relations \e{eq:Ras2} whence $u\in {\sf D}$. 
                              \end{pf}
                              

     \subsection{Self-adjoint extensions}

First, we extend  the asymptotic formulas of Theorem~\ref{LC} to all vectors  $ u\in{\cal D}(J_{\rm max})$.
Using Corollary~\ref{res1}, we define the number $\Gamma(z;  u)$ by  relation \e{eq:qres3}, that is, 
  \[
  \Gamma(z;  u) P(z)= u- {\cal R} (z) (J_{\rm max}-zI) u, \q z\in {\Bbb C}.
 \]



  \begin{theorem}\label{LCR}
    Let the assumptions of Theorem~\ref{GSnC}  be satisfied. Then all
  sequences  $ (u_{n})\in{\cal D}(J_{\rm max})$  have asymptotics
  \begin{equation}
u_{n} = a_{n} ^{-1/2} \big(s_{+}   \Theta_{n}^{(+)} +s_{-}  \Theta_{n}^{(-)} + o(1)\big)      , \q n\to\infty,
\label{eq:LC2}\end{equation} 
with some coefficients $s_{\pm} = s_{\pm} (u) $. They   can be constructed by the relations
 \begin{equation}
  \begin{split}
s_{+}(u)=  \Gamma(z;  u) \sigma_{+} (z)+ \la ( J_{\rm max}  -zI) u, P(\bar{z})\ra \tau_{+} (z),
  \\
 s_{-}(u)= \Gamma(z;  u) \sigma_{-} (z)+ \la (J_{\rm max} - zI) u, P(\bar{z})\ra \tau_{-} (z), 
  \end{split}
  \label{eq:LC31}\end{equation}  
  where the number $z\in {\Bbb C}$ is arbitrary.

  Conversely, for arbitrary $ s_{+},s_{-}  \in {\Bbb C} $, there exists a vector  
$u \in{\cal D}(J_{\rm max})$ such that    asymptotics \e{eq:LC2}  holds.
 \end{theorem}
 
 
   \begin{pf}
   According to Corollary~\ref{res1}  a vector   $u \in{\cal D}(J_{\rm max})$  admits representation  \e{eq:qres3} where the operator ${\cal R} (z)$  is defined by equality \e{eq:RR11}.
      In view of   relation \e{eq:Ras}  and asymptotics   \e{eq:LC2q} we have
    \begin{equation}
 ( {\cal R} (z)h)_{n} =    a_{n} ^{-1/2}  \big( \tau_{+} (z) \Theta_{n}^{(+)} + \tau_{-} (z)  \Theta_{n}^{(-)} \big) \la h, P(\bar{z})\ra     + o(a_{n} ^{-1/2}), \q n\to\infty,
\label{eq:RR1r}\end{equation}
  for all vectors $h \in \ell^2 ({\Bbb Z}_{+})$.  Therefore it follows from \e{eq:qres3}  and \e{eq:LC1P} that
    \begin{align*}
u_{n} = &a_{n} ^{-1/2}  \Gamma(z; (J_{\rm max} -zI)u)  \big(\sigma_{+} (z)  \Theta_{n}^{(+)}  + \sigma_{-} (z)
 \Theta_{n}^{(-)}\big) 
\\
&+ a_{n} ^{-1/2} \big(\tau_{+} (z)  \Theta_{n}^{(+)}  + \tau_{-} (z)  \Theta_{n}^{(-)}  \big) \la  ( J_{\rm max}-zI) u, P(\bar{z})\ra    + o(a_{n} ^{-1/2})   
 \end{align*} 
as $ n\to\infty$. This yields  relation \e{eq:LC2} with the coefficients $s_{\pm}$ defined by \e{eq:LC31}.

Conversely, given $ s_{+} $ and $ s_{-} $ and fixing some $z\in{\Bbb C}$, we consider the system of equations 
 \begin{equation}
  \begin{split}
s_{+}=  \Gamma\sigma_{+} (z)+ \la h , P(\bar{z})\ra \tau_{+} (z),
  \\
 s_{-} = \Gamma \sigma_{-} (z)+ \la h , P(\bar{z})\ra \tau_{-} (z).
  \end{split}
  \label{eq:LC3x}\end{equation}  
   for  $\Gamma$ and $ \la h, P(\bar{z})\ra $. According to \e{eq:Wro} the determinant  of this system is not zero,  so that  $\Gamma$ and $ \la h, P(\bar{z})\ra $ are uniquely determined by  $s_{+} $ and $ s_{-} $.
 Then we take any $h$ such that its scalar product with $ P(\bar{z})$ equals the found value of  $ \la h, P(\bar{z})\ra$.
  Finally, we define $u$ by formula  \e{eq:qres3}. The asymptotics as $n\to\infty$ of $P_{n}(z)$   and $( {\cal R} (z)h)_{n} $ are given by formulas \e{eq:LC1P} and \e{eq:RR1r}, respectively. In view of equations \e{eq:LC3x} this leads to  asymptotics \e{eq:LC2}.
       \end{pf}
       
       
      Theorem~\ref{LCR}  yields a mapping $ {\cal D}(J_{\rm max})\to{\Bbb C}^2$ defined by the formula
       \begin{equation}
u\mapsto (s_{+}  (u) , s_{-}(u)).
\label{eq:mapping}\end{equation}
      The construction of      Theorem~\ref{LCR} depends on the choice of $z\in{\Bbb C}$, but this mapping is defined intrinsically.  In particular, we can set $z=0$ in all formulas of  Theorem~\ref{LCR}.  Note that mapping \e{eq:mapping} is surjective.

Evidently,   \e{eq:mapping} plays the role of the mapping  $f\mapsto (f(0) , f' (0)) $  for the differential operator $-d^2/ dx^2$ in the space $L^2 ({\Bbb R}_{+})$ and formula \e{eq:qres4} plays the role of the integration-by-parts formula
\[
\int_{0}^\infty f'' (x) \ov{g} (x) dx-\int_{0}^\infty f (x) \ov{g''} (x) dx= f(0) \bar{g}' (0)- f'(0) \bar{g} (0).
\]
 
     Under the assumptions of Theorem~\ref{GSS} the right-hand side of \e{eq:qres4} can be expressed  in terms of the coefficients $s_{+}$ and $s_{-}$.  Recall that the numbers $\alpha_{n}$, $\beta_{n}$ are defined by equalities \e{eq:aabb} and  $\alpha_{\infty}$, $\beta_{\infty}$ are their limits as $n\to\infty$.

 \begin{proposition}\label{resAS}
     For all $u, v \in {\cal D} (J_{\rm max})$, we have the identity
    \begin{equation}
\la J_{\rm max}u,v \ra - \la u, J_{\rm max}v \ra = 2 i \alpha_{\infty}^{-1} \sqrt{ 1-\beta_{\infty}^2}\big( s_{+}(u) \ov{s_{+}(v) }- s_{-}(u) \ov{s_{-}(v) }\big)  .
\label{eq:ABS}\end{equation}
 \end{proposition} 
   
    \begin{pf}
    It follows from formula \e{eq:LC2}  that
       \begin{align*}
\sqrt{a_{n+1}  a_{n } } (u_{n +1}\bar{v}_{n } &- u_{n }\bar{v}_{n +1})
\\
= &\big(s_{+}(u)   \Theta_{n+1}^{(+)} + s_{-} (u)  \Theta_{n+1}^{(-)} \big)\big(\ov{s_{+}(v) }  \Theta_{n}^{(-)} + \ov{ s_{-} (v)}  \Theta_{n}^{(+)}\big)
\\
- & \big(s_{+}(u)   \Theta_{n}^{(+)} + s_{-} (u)  \Theta_{n}^{(-)}\big)\big(\ov{s_{+}(v) }   \Theta_{n+1}^{(-)}  + \ov{ s_{-} (v)}   \Theta_{n+1}^{(+)} \big) + o(1)
\\
=&\big( s_{+}(u) \ov{s_{+}(v) }-s_{-}(u) \ov{s_{-}(v) }\big) (e^{i\arccos\beta_{n}} -e^{-i \arccos\beta_{n}}) + o(1).
 \end{align*}
Passing here to the limit $n\to\infty$ and using equality \e{eq:qres4}, we obtain identity \e{eq:ABS}.
       \end{pf}
       
    We can now characterize set \e{eq:MIN}.
       
       
       \begin{proposition}\label{ABS}
       A vector
    $v  \in {\cal D} (J_{\rm max})$ belongs to     $ {\sf D}  $ if and only if
    $v_{n}= o(a_{n}^{-1/2})$, that is, 
     \begin{equation}
    s_{+} (v) = s_{-} (v) =0.
      \label{eq:ABS4}\end{equation}
    \end{proposition} 
     
     \begin{pf}    
     A vector $v$ belongs to     $ {\cal D} ( J_{\max}^{*})$ if and only if
       \begin{equation}
    \la J_{\max}u,v \ra= \la u, J_{\max}v \ra
  \label{eq:ABS2}\end{equation}
    for all   $u  \in {\cal D} (J_{\rm max})$.  According to Proposition~\ref{resAS},   equality \e{eq:ABS2} is equivalent to
     \begin{equation}
 s_{+}(u) \ov{s_{+}(v) }- s_{-}(u) \ov{s_{-}(v) }=0.
  \label{eq:ABS3}\end{equation}
  This  is of course true if \e{eq:ABS4} is satisfied. Conversely, if  \e{eq:ABS3} is satisfied for all  $u\in {\cal D} ( J_{\max})$, we use that according to Theorem~\ref{LCR}  the numbers $ s_{+}(u)$ and $s_{-}(u)$  are arbitrary. This implies \e{eq:ABS4}.
         \end{pf}       
    
This result shows that  \e{eq:mapping} considered as a mapping of the factor space  $  {\cal D} (J_{\rm max})/ {\sf D}$ onto ${\Bbb C}^2$ is injective.

    All self-adjoint extensions ${ J}_{\omega} $ of the operator $J_{\min}$ are now parametrized by complex numbers $\omega\in{\Bbb T}\subset {\Bbb C}$.  Let the set $ {\cal D} (  { J}_\omega)\subset {\cal D} (    J_{\max})$ of vectors $u$  be distinguished by  the condition    
    \begin{equation}
   s_{+} (u) =\omega s_{-} (u) , \q |\omega|=1 .
      \label{eq:gen2}\end{equation}

      \begin{theorem}\label{ABX}
           Let the assumptions of Theorem~\ref{GSS}  be satisfied. Then
     for all $\omega\in{\Bbb T}$, the operators ${ J}_{\omega}$ are  self-adjoint. Conversely, every   operator $J$ such that 
     \begin{equation}
  J_{\min} \subset  J =J^*\subset J_{\max} 
      \label{eq:ABX2}\end{equation}
      equals ${ J}_{\omega}$ for some $\omega\in{\Bbb T}$.
    \end{theorem} 
    
    \begin{pf}
    We proceed from Proposition~\ref{resAS}.   If $u,v \in  {\cal D} (  { J}_\omega)$,  it follows from condition \e{eq:gen2} that   $s_{+}(u) \ov{s_{+}(v) }= s_{-}(u) \ov{s_{-}(v) }$. Therefore   according to equality \e{eq:ABS} $\la { J}_{\omega}u,v\ra = \la u, { J}_{\omega}v \ra$ whence ${ J}_{\omega}\subset { J}_{\omega}^*$. If $ v\in {\cal D} (  { J}_\omega^*)$, then 
   $\la { J}_{\omega}u,v \ra= \la u, { J}_{\omega}^* v \ra $
    for all $ u\in {\cal D} (  { J}_\omega)$ so that in view of \e{eq:ABS} equality \e{eq:ABS3} is satisfied.
    Therefore $s_{-}(u) (\omega\ov{s_{+}(v) }- \ov{s_{-}(v) })=0$. Since $s_{-}(u)$ is arbitrary, we see  that  $ \omega\ov{s_{+}(v) }- \ov{s_{-}(v) }=0$, and hence $ v\in {\cal D} (  { J}_\omega)$.
    
    Suppose that an operator $J$  satisfies conditions  \e{eq:ABX2}.  Since  $J$ is symmetric, it follows from  Proposition~\ref{resAS} that  equality   \e{eq:ABS3}   is true  for all $ u, v\in {\cal D} (  J )$.  Setting here $u=v$, we see that $|  s_{+}(v)|=|   s_{-}(v)|$.   There exists a vector $v_{0} \in {\cal D} (  J )$ such that
    $s_{-}(v_{0})\neq 0$ because  $J\neq \clos J_{\rm min}$. Let us set $\omega= s_{+}(v_{0})/ s_{-}(v_{0})$. Then $|\omega|=1$  and relation  \e{eq:gen2}  is a direct consequence of \e{eq:ABS3}.
        \end{pf}

       \subsection{Resolvent}

     Now it easy to construct the   resolvent  of the operator ${ J}_\omega$  defined in the previous subsection.
     We previously note that, by definition \e{eq:LC2}, 
     \[
     s_{\pm} (P(z))= \sigma_{\pm} (z) \q \mbox{and} \q      s_{\pm} (Q(z))= \tau_{\pm} (z).
     \]
     
     \begin{theorem}\label{ABX1}
            Let the assumptions of Theorem~\ref{GSS}  be satisfied. Then
     for all $z\in {\Bbb C}$ with $\Im z\neq 0$ and all $h\in\ell^2 ({\Bbb Z}_{+})$, the resolvent ${ R}_{\omega} (z)= ({ J}_{\omega}-zI)^{-1}$ of the operator ${ J}_\omega$ is given by   the equality 
        \begin{equation}
{ R}_{\omega} (z) h = {\gamma}_{\omega} (z) \la h,  P(\bar{z})\ra P(z)+ {\cal R}(z)h 
      \label{eq:ABY1}\end{equation}
      where   
           \begin{equation}
 {\gamma}_{\omega}(z)= -\frac{\tau_{+}(z)-\omega \tau_{-}(z)} {\sigma_{+}(z)-\omega \sigma_{-}(z)}
    \label{eq:gam}\end{equation}
    and $ \sigma_{\pm}(z)$,  $ \tau_{\pm}(z)$  are  entire functions of $z$ defined by relations \e{eq:st}, \e{eq:ts}.
    \end{theorem} 
    
         \begin{pf}
           According to      Corollary~\ref{res1} the vector  $u={ R}_{\omega}(z) h$ is given by formula \e{eq:qres3} where the coefficient $\Gamma$ is determined by condition 
            \e{eq:gen2}. It follows from Theorem~\ref{LCR}  that the components $u_{n}$ of $u$ have  asymptotics \e{eq:LC2}  with the coefficients $s_{\pm}$ defined by relations 
            \e{eq:LC3x}.  Thus, 
          $u\in {\cal D}({ J}_{\omega})$ if and only if  
    \[
    \Gamma \sigma_{+} (z) + \tau_{+} (z) \la h,  P(\bar{z})\ra=\omega \big( \Gamma \sigma_{-} (z) + \tau_{-} (z)\la h,  P(\bar{z})\ra\big)
    \]
    whence
    \[
\Gamma= -\frac{\tau_{+}(z)-\omega \tau_{-}(z)} {\sigma_{+}(z)-\omega \sigma_{-}(z)}\: \la h,  P(\bar{z})\ra.
      \]
Substituting this expression into \e{eq:qres3}, we arrive at  formulas
    \e{eq:ABY1},   \e{eq:gam}.
          \end{pf}
          
          We emphasize that, for various $\omega \in{\Bbb T}$, the resolvents  $R_\omega(z)$ of the operators $J_\omega$ differ from each other only by the coefficient $\gamma_\omega(z)$ at the rank-one operator $\la\cdot, P(\bar{z})\ra P(z)$.
   Observe also that
    $  \ov{\gamma_\omega(z)}= \gamma_\omega (\bar{z})$.

     Since  $\la {\cal R}(z) e_{0}, e_{0}\ra=0$, we see that $\la R_\omega(z) e_{0}, e_{0}\ra=\gamma_\omega (z)$. Thus, Theorem~\ref{ABX1}  implies the classical Nevanlinna representation obtained in \cite{Nevan}  for the Cauchy-Stieltjes transform of the spectral measures $d\rho_\omega (\lambda)= d (E_\omega (\lambda) e_{0}, e_{0})$ of the operators $J_\omega$.
     
     \begin{corollary}\label{RESc}
     For all $z\in {\Bbb C}$ with $\Im z\neq 0$,  we have
     \[
     \int_{-\infty}^\infty
 (\lambda-z)^{-1} d\rho_\omega(\lambda)= \gamma_\omega (z).
 \]
    \end{corollary} 
    
      \begin{corollary}\label{RESb}
     If $z\in {\Bbb C}$  is  not an eigenvalue of the operator $J_\omega$, then its resolvent $R_\omega  (z)$ is in the Hilbert-Schmidt class.  Thus, the spectra of the operators $J_\omega$   are discrete and consist of the points  $z$ where
          \begin{equation}
          \sigma_{+}(z)-\omega \sigma_{-}(z)=0.
         \label{eq:ABY3}\end{equation}
 \end{corollary}

          Note also that according to \e{eq:LC2p} $\sigma_{+} (z)\neq \sigma_{-} (z)$ if $\Im z\neq 0$, and therefore the zeros $z$ of equation \e{eq:ABY3} lie on the real axis.  This result has  of course to be expected since   $z$ are eigenvalues of the self-adjoint operator ${ J}_{\omega}$.  We finally note that the discreteness of the spectrum of the operators ${ J}_{\omega}$ is quite natural because their domains ${\cal D}({ J}_{\omega})$ are distinguished by boundary conditions at the point $n=0$ 
 and for $n\to\infty$. Therefore  ${ J}_{\omega}$ acquire some features of regular operators.

 \section{Jacobi operators with stabilizing coefficients}  
 
 We here study the case of stabilizing recurrence coefficients  satisfying condition \e{eq:stabi}. Then the Jacobi operator $J$ 
 differs by a compact term from the operator  with the coefficients $a_{n}=a_{\infty} > 0$, $b_{n}=b_{\infty}$ for all $n\in{\Bbb Z}_{+}$. Under additional assumptions  \e{eq:Tr} or \e{eq:LR}  the corresponding orthogonal polynomials are rather close to Jacobi polynomials (see Appendix~B).  We follow the scheme  already presented for the case of increasing coefficients satisfying condition  \e{eq:Haupt}  where $|\beta_{\infty}  | < 1$ and  \e{eq:Carl}.  Since the case of stabilizing recurrence coefficients was investigated  in \cite{JLR}, some technical details will be omitted.
 
 
    \subsection{Compact perturbations}
   
   As a preliminary remark, we note that, without a loss of generality, we may suppose that $a_{\infty}=1/2$ and 
  $b_{\infty}=0$ in \e{eq:stabi}, that is, condition \e{eq:comp} is satisfied. 
  Indeed, let $\wt{J}$ be the Jacobi operator with off-diagonal elements $\ti{a}_{n}=(2a_{\infty})^{-1} a_{n}$ and diagonal elements $\ti{b}_{n}=(2a_{\infty})^{-1} (b_{n}-b_{\infty})$. Then
    \begin{equation}
 J= 2a_{\infty} \wt{J} + b_{\infty }I
   \label{eq:Stabi1}\end{equation} 
  and $\ti{a}_{n}\to 1/2$, $\ti{b}_{n}\to 0$  as $n\to\infty$ according to  \e{eq:stabi}. In terms of the orthonormal polynomials, relation \e{eq:Stabi1} means that
  \[
    \wt{P}_{n}(z)= P_{n} (2a_{\infty} z + b_{\infty }).
    \]
    
    Let $J_{0}$ be the ``free'' Jacobi operator. It has the coefficients $a_{n}=1/2$, $b_{n}=0$ for all $n\in{\Bbb Z}_{+}$. Under assumption \e{eq:comp} the operator
    $J-J_{0}$ is compact so that the essential spectrum of $J$ coincides with the interval $[-1,1]$.
         Note that   \e{eq:comp} corresponds to the assumptions $a(x)\to 1/2$, $b(x)\to 0$ as $x\to\infty$  for the differential operator \e{eq:HHh}. We also require conditions \e{eq:Tr} or \e{eq:LR} corresponding, respectively,  to  ``short-range''   and ``long-range'' 
perturbations  of the  ``free"  differential   operator $D^2$.

      Our notation is very close to that used in  the case of increasing coefficients. 
     Recall that $\Pi={\Bbb C}\setminus{\Bbb R}$.
       As before, we  suppose that  the spectral parameter $z$ belongs to a bounded subset of $\clos \Pi$, but  we now additionally assume that $z\neq \pm 1$.              For $z\in\clos\Pi$, we set
 \begin{equation}
 z_{n}  = \frac{z-b_{n}} {2a_{n}}
\label{eq:aa}\end{equation}
which is now slightly more convenient than \e{eq:aabb1}. According to \e{eq:comp} we have
\begin{equation}
\lim_{n\to\infty} z_{n}  = z.
\label{eq:Stabi2}\end{equation}
It follows that, for sufficiently large $N_{0}$, the set of points $z_{n}$ where $n\geq N_{0}$ is separated from the points $1$ and $-1$. 
Relations  \e{eq:Stabi2} and \e{eq:aabb2}   are quite different but play  similar roles in our presentation.

 We start with the  general case when the   long-range condition \e{eq:LR} is satisfied.  Short-range perturbations will be specially discussed in Sect.~9.

 \subsection{Jost solutions and orthogonal polynomials}

  Here we suppose that conditions \e{eq:comp} and  \e{eq:LR}   are satisfied.  The first step is to distinguish the Jost solutions $f_{n}(z)$ of the Jacobi equation \e{eq:Jy}    by their asymptotic behavior  for $n\to\infty$. The corresponding Ansatz
     ${\cal A}_{n} (z)$ can be again  defined  by formula \e{eq:AnsDD}, but we will omit the factor $ a_{n}^{-1/2}$  since it tends to a   constant. 
   We always suppose that the values of the spectral parameter $z$  are separated from the points $\pm 1$ and fix $N_{0}$ in  such a way that estimate \e{eq:aabb5} holds true  for all $n\geq N_{0}$.
Now we set
       \[
{\cal A}_{n}(z) =   \z_{N_{0}}\z_{N_{0}+1}
  \cdots \z_{n-1} \q\mbox{where}  \q \z_{m}=\z (z_m),
\]
so that ${\cal A}_{n}(z)$ coincides with the sequence $q_{n}(z)$ defined by equality \e{eq:re1}.
Then the remainder  equals
\begin{align}
 r_{n} (z) : &=   {\cal A}_{n} (z)^{-1} \Big(a_{n-1}{\cal A}_{n-1} (z) + (b_{n}-z){\cal A}_{n} (z) + a_{n}{\cal A}_{n+1} (z)\Big)
 \nonumber\\
 &=a_{n-1} \z_{n-1}^{-1} +a_{n} \z_{n} + b_{n}-z
  \nonumber\\
   &=  (a_{n-1} -a_{n})\z_{n-1}^{-1} +a_{n} (\z_{n-1}^{-1} - \z_{n}^{-1}).
\label{eq:Grst}\end{align}

The following statement plays the role of Lemma~\ref{Gr}.

\begin{lemma}\label{re}
For $z\in\clos{\Pi}\setminus\{-1,1\}$, an estimate
\[
 | r_{n} (z) | \leq C  \frac{ |a_{n}-a_{n-1}| + |b_{n}-b_{n-1}|}
 {\sqrt{z_{n-1}^2-1}+\sqrt{z_{n}^2-1}} 
\]
holds.
In particular, $( r_{n} (z) )  \in \ell^1 ({\Bbb Z}_{+})$ if assumptions  \e{eq:comp} and \e{eq:LR} are satisfied.
\end{lemma}

\begin{remark}\label{Grem}
The crucial difference between the cases of stabilizing and increasing recurrence coefficients is that
$a_{n}a_{n-1}^{-1}-1 \in \ell^1 ({\Bbb Z}_{+})$ in the first case while $a_{n}a_{n-1}^{-1}=1+ p n^{-1}+ O(n^{-2})$
if, for example,  $a_{n}=(n+1)^p$ for some $p>0$. In particular, this is the reason why we had to introduce the factor $a_{n}^{-1/2}$ in Ansatz \e{eq:AnsDD}. Indeed, without this factor in the case of  increasing coefficients $a_{n}$,  remainder \e{eq:Gry}
    differs from \e{eq:Gr6} by a term of order $n^{-1}$ and hence does not belong to $\ell^1 ({\Bbb Z}_{+})$.
  \end{remark}

Making again  the multiplicative change  of variables  
         \begin{equation}
f_{n}(  z )= {\cal A}_{n}(z)u_{n}( z ) \q\mbox{where}\q {\cal A}_{n}(z)=q_{n} (z),
\label{eq:JJs}\end{equation}
we see that
 equation  \e{eq:Jy} for   $ f_{n} (z)$ is equivalent to the equation
\begin{equation}
 a_{n} \z_{n} ( u_{n+1} (z)- u_{n} (z)) - a_{n-1} \z_{n-1}^{-1} ( u_{n} (z)- u_{n-1} (z))=- r_{n} (z) u_{n} (z), \q n\in {\Bbb Z}_{+},
\label{eq:A12}\end{equation}
for  the sequence  $ u_{n} (z)$. This equation is quite similar to  \e{eq:Gs5} and plays the same role.  The condition 
 \begin{equation}
f_{n}(  z )=q_{n}( z )\big(1 + o( 1)\big), \q n\to\infty,
\label{eq:A22}\end{equation}
means of course that $u_{n} (z) \to 1$ as $n\to\infty$.

Next, we reduce the  difference equation \e{eq:A12} for $u_{n} (z)$
  with  this condition 
         to a  ``Volterra integral" equation \e{eq:A17}
 with the kernel 
\begin{equation}
G_{n,m} (z) =  q_{m}(z)^2 \sum_{p=n}^{m-1} ( a_{p} \z_{p })^{-1}  q_{p }(z)^{-2},\q n,m\in{\Bbb Z}_{+} , \q m\geq n+1.
\label{eq:Gg1}\end{equation}
Note that $G_{n,m} (\bar{z})=\ov{G_{n,m} (z)}$.
The functions $G_{n,m} (z) $ are analytic in $z\in\Pi$ and are continuous up to the real axis. 

 The following assertion plays the crucial role in our analysis of equation \e{eq:A17}, in particular, for $z$ lying on the cut along $ (-1,1)$. It shows that   kernels \e{eq:Gg1}  are bounded   uniformly in $n$ and $m$  provided  some neighborhoods of the points $ \pm 1$ are excluded.

 \begin{lemma}\label{eik1}
There exist constants $C(z)$ and $N(z)$ such that  an estimate 
\[
| G_{n,m} (z) |\leq  C (z)<\infty, \q m-1\geq n \geq N(z),
\]
 is true for all $z\in \clos\Pi\setminus \{-1,1\}$. The constants $C(z)$ and $N(z)$ are common for $z$ in  compact subsets of $\clos\Pi\setminus \{-1,1\}$, that is, for all
$ z\in \clos\Pi$ such that $|z^2-1|\geq \epsilon$ and $|z|\leq \rho$ where   $\epsilon>0$ and $\rho <\infty$ are some fixed numbers. 
 \end{lemma}

  Lemmas~\ref{re} and \ref{eik1} allow us to solve the Volterra equation \e{eq:A17} with the coefficients  \e{eq:Grst},  \e{eq:Gg1} by iterations. Similarly to Sect.~3.4, we set $u^{(0)}_n (z)=1$ and define  $u^{(k)}_n (z)$
  recursively by relation \e{eq:W5}.  Then estimates \e{eq:W6s} hold true and a solution $u_{n}$ of \e{eq:A17}  is built as series \e{eq:W8}. The following assertion plays the role of Theorem~\ref{GS3}.


  \begin{theorem}\label{eik2}
 Let assumptions \e{eq:comp} and \e{eq:LR} be satisfied. 
  For  $z\in \clos \Pi\setminus\{-1,1\}$,
equation  \e{eq:A17} has a $($unique$)$ bounded solution $\{u_{n}( z )\}_{n=0}^\infty$. This sequence obeys an estimate
 \begin{equation}
| u_{n}( z )-1|\leq C \sum_{m=n}^\infty (|a_{m}'| + |b_{m}'| )
\label{eq:A20}\end{equation} 
where the constant $C$ is common for $z$ in  compact subsets of $\clos\Pi\setminus\{-1,1\}$. 
For all $n\in {\Bbb Z}_{+}$, the functions $u_{n}( z )$ are analytic in $z\in  \Pi$  and are continuous up to the cut along $\Bbb R$
with possible exception of the points $z=-1$ and $z=1$.
  \end{theorem}

 
 Similarly to Lemma~\ref{GS4}, we check that $u_{n} (z)$ satisfies also 
  the  difference equation \e{eq:A12}.     
 Then we define $f_{n}(  z )$ by formula \e{eq:JJs}.
Since equations  \e{eq:A12} for $u_{n}(z)$ and  \e{eq:Jy}  for $f_{n}(z)$ are equivalent, $f_{n}(  z )$ satisfies the Jacobi equation  \e{eq:Jy}.
 Obviously,      estimate  \e{eq:A20} on  $u_{n}( z )$ implies   asymptotics \e{eq:A22} of $f_{n}(z)$.
  Thus we arrive at the following result.
    
    \begin{theorem}\label{EIK}
   Let assumptions \e{eq:comp} and \e{eq:LR} be satisfied,   let  $z\in \clos\Pi  \setminus \{-1,1\}$, and let       $q_{n}(z)$ be defined by equality \e{eq:re1}.
Denote by $u_{n}( z )$   the sequence constructed in Theorem~\ref{eik2}. Then  the sequence $f_{n}( z )$ defined by equality \e{eq:JJs} satisfies equation \e{eq:Jy},  and it has  asymptotics
\e{eq:A22}.  For all $n\in {\Bbb Z}_{+}$, the functions $f_{n}( z )$ are analytic in $z\in  \Pi$  and are continuous up to the cut along $\Bbb R$
with possible exception of the points $z=-1$ and $z=1$. Relation  \e{eq:AcGG} is satisfied. Asymptotics \e{eq:A22} is  uniform in $z$ from compact subsets    of the set $  \clos\Pi\setminus\{-1,1\} $.
 \end{theorem}

   Recall that the polynomials   $P_{n}(z)$  are solutions of equation  \e{eq:RR} satisfying  conditions  \e{eq:RR1}.
   Put $ P(z)= ( P_{n}(z) )_{n=-1}^\infty$, $ f(z)= ( f_{n}(z) )_{n=-1}^\infty$. As before,   the sequence $ (f_{n}( z ))_{n=-1}^\infty$   will be called
  the   Jost solution of equation \e{eq:Jy}, and  the   Jost function $ \Omega (z)$ is defined as Wronskian  \e{eq:J-W}. For the operator $J_{0}$,     the Jost solution is $ (\z(z)^n )_{n=-1}^\infty$ and  the corresponding Wronskian
          \begin{equation}
     \Omega_{0} (z) = -(2\z (z) )^{-1}.
     \label{eq:WRfr}\end{equation}

     The following result is a direct consequence of Theorem~\ref{EIK}.

\begin{corollary}\label{JOST}
The Wronskian $\Omega(z)$   depends analytically on $z\in\Pi$, and it is a continuous function of  $z$ up to the cut along $\Bbb R$ except, possibly, the points $\pm 1$. 
\end{corollary}

  \begin{remark}\label{AR}
  Suppose that $z \in (-\infty,- 1)\cup (1,\infty)$. Then according to  \e{eq:Stabi2} we also have $z_{n} \in (-\infty,- 1)\cup (1,\infty)$ if $n$ is sufficiently large. Since the function $\z(z)$ is real for $z \in (-\infty,- 1)\cup (1,\infty)$, it follows that the numbers $\z(z_{n})$, $q(z_{n})$ and hence $f_{n}(z)$ are real.  Using now  relation  \e{eq:AcGG}, we see that values of $f_{n}(z)$ on the upper and lower edges of the cut along $  (-\infty,- 1)\cup (1,\infty)$ coincide. Therefore the functions $f_{n} (z)$ are actually analytic on the whole set 
  \[
  \Pi_{0}={\Bbb C}\setminus [-1,1].
  \]
   \end{remark}


 An asymptotics as $n\to\infty$ of  product \e{eq:re1}  can be found rather explicitly.  
     
      \begin{lemma}\label{asq}
   Let assumption \e{eq:comp}  be satisfied,  let $z\in\clos\Pi\setminus\{-1,1\}$,  and let $\z=\z (z)$ be given by equality
   \e{eq:ome}. Then 
     \begin{equation}
q_{n} (z)= e^{n(\ln\z + o (1))}, \q n\to\infty.
 \label{eq:asq}\end{equation} 
 \end{lemma}
  
   \begin{pf}
   Let the sequence $z_{n}$   be defined by formula  \e{eq:aa} and   $\z_{n} =\z (z_{n})$.
   It follows from \e{eq:comp}  that $z_{n}=z+ o(1)$ and hence
   $
\z_{n}  = \z (1 + \epsilon_{n})  
$
 where $\epsilon_{n}\to 0$ as $ n\to\infty$. For product \e{eq:re1}, this yields
    \[
\ln q_{n}  =n \ln \z+ \sum_{m=N_{0}}^{n-1} \ln (1 + \epsilon_m )  = n \ln \z+ o(n)
 \]
 which is equivalent to \e{eq:asq}.
        \end{pf}
        
 

{\bf Asymptotics in the complex plane.} 
 Supposing that $z\not\in{\Bbb C}\setminus [-1, 1]$,   we
 introduce a solution $g_{n} (z)$ of equation  \e{eq:Jy} exponentially growing as $n\to\infty$.  
 A  proof of the following  result  is similar to that of Theorem~\ref{GEe}. Details
  can be found in \cite{JLR}.

\begin{theorem}\label{GEx}
Let   $z\in{\Bbb C}\setminus [-1,1]$, and let  assumptions \e{eq:comp} and \e{eq:LR}  be satisfied. If
$f_{n} (z)$ is the Jost solution of equation \e{eq:Jy}, then the solution $g_{n} (z)$ of the same equation defined by 
      \e{eq:GE},      \e{eq:GEx}
 satisfies a relation
     \[
\lim_{n\to\infty} q_{n}(z)  g_{n}(z)=\frac{1}{\sqrt{z^2-1} }.
\]
 \end{theorem}

Now it is easy to find asymptotics of the orthogonal polynomials $P_{n} (z)$  for $z\in{\Bbb C}\setminus [-1,1]$.
By Theorem~\ref{GE}, the Wronskian $ \{ f(z),g(z)\}$  of $f (z)= ( f_{n}(z) )$ and $g (z)= ( g_{n}(z)   )$ equals $1$,
whence these  solutions of equation  \e{eq:Jy} are linearly independent.
This yields relation \e{eq:Pf}  where $\Omega(z)$ is given by \e{eq:J-W} and $\omega(z)= \{ P(z),g(z)\}$.    Obviously,  
$\omega(z)\neq 0$ if $\Omega(z)= 0$, that is, $z$ is an eigenvalue of the operator $J$. Therefore Theorems~\ref{EIK}   and \ref{GEx} imply the following result (cf. Theorems~\ref{GEL1}. and \ref{GE1+}).

\begin{theorem}\label{GE1}
  Under assumptions \e{eq:comp} and \e{eq:LR}  the relation  
   \begin{equation}
\lim_{n\to\infty} q_{n}(z)  P_{n}(z)=-\frac{\{ P(z), f(z)\} }{\sqrt{z^2-1}}
\label{eq:GEGEy}\end{equation}
is true for all  $z\in{\Bbb C}\setminus [-1,1]$  
with convergence uniform on compact subsets of $z\in{\Bbb C}\setminus [-1,1]$.
Moreover, if $\Omega(z)=0$, then
 \begin{equation}
\lim_{n\to\infty}  q_{n}(z)^{-1} P_{n}(z)=\{ P(z),g(z)\} \neq 0.
\label{eq:GEGEx}\end{equation}
 \end{theorem}

{\bf Asymptotics on the continuous spectrum.} 
Suppose  that $z=\lambda\pm i0$ where $\lambda\in (-1,1)$ so that $\lambda=\cos\theta$, $\theta\in (0,\pi)$. Now relations \e{eq:wwm} and \e{eq:ww3}  remain true with the numbers $\lambda_{n}$ given by
\begin{equation}
\lambda_{n} : =   \frac{\lambda-b_{n}}{2a_{n}}. 
 \label{eq:Joz}\end{equation}
 

   The following   result is a direct consequence of  Theorem~\ref{EIK}. It plays the role of Theorem~\ref{As+} stated for the case of increasing coefficients.
 
 \begin{theorem}\label{As}
 Let assumptions~\e{eq:comp} and \e{eq:LR} be satisfied.  For $ \lambda\in (-1,1)$, define the phases $\varphi_{n}(\lambda)$  
   by formula \e{eq:ww3} where $\lambda_{n}$ are numbers \e{eq:Joz}.
 Then
   \begin{equation}
f_{n}(\lambda\pm i0)=   e^{\mp i \varphi_{n}(\lambda)}  (1+ o(1)), \q n\to\infty.
   \label{eq:Jost1}\end{equation}
 \end{theorem}
 
Note that
 \[
 \varphi_{n}(\lambda)=n\arccos\lambda + o(n),
 \]
 but under additional assumptions the error term can be made more explicit.  In particular, we see that
 asymptotics \e{eq:Jost1} of  $ f_{n}(\lambda\pm i0)$    is oscillating as $n\to\infty$.  


  To find asymptotic behavior of the polynomials $P_{n} (\lambda)$ for $\lambda\in (-1,1)$,  that is, on the continuous spectrum of the Jacobi  operator $J$,
we have to consider two complex conjugate Jost solutions $f(\lambda\pm i0)= \big(f_{n} (\lambda\pm i0)\big)_{n=-1}^\infty$ for
$ \lambda\in (-1,1)$. 
Similarly to Lemma~\ref{Asc}, we  have 

\begin{lemma}\label{AsC}  
The Wronskian \e{eq:Wr} of $f (\lambda+ i0)$  and $f  (\lambda- i0)$  
equals
   \begin{equation}
\{f (\lambda+i0),  f (\lambda-i0)\}=   i   \sqrt{1- \lambda^2 }\neq 0, \q \lambda\in  (-1,1),
      \label{eq:wwc}\end{equation} 
and hence these solutions  are linearly independent.
\end{lemma}

  Note that compared to \e{eq:ww}, 
  $\sqrt{1-\beta^2_{\infty}}$ is   replaced by $\sqrt{1- \lambda^2 }$ because   the  role  of relation \e{eq:aabb2} is now played by \e{eq:Stabi2}. 
Besides,   the coefficient $2$ has  disappeared  in  \e{eq:wwc} because \e{eq:Jost1} does not contain an amplitude factor.

    Lemma~\ref{AsC} implies the following two  results (cf. Lemma~\ref{HH+} and Theorem~\ref{HX+}).

\begin{lemma}\label{HH}
 For $ \lambda\in (-1,1)$,  the representation 
  \begin{equation}
P_{n} (\lambda)=\frac{  \Omega (\lambda-i0) f_{n} (\lambda+i0) -  \Omega (\lambda+i0)  f _{n} (\lambda-i0)  }{i   \sqrt{1-\lambda^2}},    \q n=0,1,2, \ldots, 
\label{eq:HH4l}\end{equation}
 holds true.
 \end{lemma}

 \begin{theorem}\label{HX}
 The Wronskians $ \Omega (\lambda + i0)$ and $ \Omega (\lambda - i0)=
 \ov{\ \Omega (\lambda + i0) } $  are continuous functions of $\lambda\in (-1,1)$ and
 \begin{equation}
 \Omega (\lambda\pm i0) \neq 0 ,\q \lambda\in (-1,1)  .
\label{eq:HH5}\end{equation}
 \end{theorem}

Let us define the functions $\kappa (\lambda) $ and $\eta (\lambda)$ by equalities \e{eq:AP+}.
In the theory of short-range perturbations of the Schr\"odinger operator, these functions  are known as the limit amplitude and the limit phase, respectively; the function   $\eta (\lambda)$ is also called the scattering  phase or the   phase shift.  


Combined together,  relations \e{eq:Jost1} and \e{eq:HH4l} yield the   asymptotics of Bernstein-Szeg\H{o} type for the  polynomials  $P_{n} (\lambda)$
(cf. Theorem~\ref{Sz+}).

 \begin{theorem}\label{Sz}
 Let assumptions~\e{eq:comp} and \e{eq:LR} be satisfied,   let   $\lambda \in (-1,1)$ and let the phase $\varphi_{n}  (\lambda)$ be defined by formulas \e{eq:ww3} and 
 \e{eq:Joz}. Then the  polynomials $P_{n} (\lambda)$ have asymptotics 
   \begin{equation}
 P_{n} (\lambda)= \frac{2 \kappa ( \lambda)}   {\sqrt{1-\lambda^2}}\sin (\varphi_{n}  (\lambda) + \eta(\lambda) ) + o(1) 
\label{eq:Szold}\end{equation}
as $n\to\infty$. Relation \e{eq:Szold} is uniform in $\lambda$ on compact subintervals of $(-1,1)$.
 \end{theorem}

Note that Theorem~\ref{Sz} does not follow from Theorem~\ref{GE1} because asymptotics \e{eq:GEGEy} is not uniform as $z$ approaches the cut along $(-1,1)$.

Asymptotic formulas \e{eq:GEGEy} and \e{eq:Szold} are the  classical results of the Bernstein-Szeg\H{o} theory. They  are stated as Theorems~12.1.2 and 12.1.4  in the book  \cite{Sz} where the assumptions are imposed   on the spectral measure $d\rho(\lambda)$; in particular, it is assumed that
$\supp\rho\subset [-1,1]$. Under  assumptions \e{eq:comp} and \e{eq:LR} on recurrent coefficients $a_{n}$, $b_{n}$ asymptotic formulas for orthonormal polynomials were established in the paper   \cite{Mate}.  We followed here the presentation of  \cite{JLR}.

 \subsection{Spectral results}
 
 The construction of
the resolvent   of the Jacobi operator $J$ is quite similar to Theorem~\ref{AC}.  Note that for $z\not\in [-1,1]$, we have
 $|\z (z)| <1$  according to
  \e{eq:A22} and \e{eq:asq}. Therefore $ f_{n}(z)= O(\d^n)$ with some $\d<1$ as $n\to \infty$ whence $f(z)\in \ell^2 ({\Bbb Z}_{+})$. Therefore the  next statement directly follows from  Proposition~\ref{res}.

 \begin{theorem}\label{ACst}
  Let assumptions \e{eq:comp} and \e{eq:LR}  hold. Then
  \begin{enumerate}[{\rm(i)}]
 \item
The  resolvent $R(z)=(J-z I)^{-1}$ of the Jacobi operator $J$ is an integral operator with matrix elements  \e{eq:Rrpm}.
 For all $n,m\in{\Bbb Z}_{+}$,  it    is an analytic function of $z\in {\Bbb C}\setminus [-1,1]$ with simple poles at eigenvalues of the operator $J$. A point
$z\in {\Bbb C}\setminus [-1,1]$ is an eigenvalue of $J$ if and only if $\Omega(z)=0$.

 \item
 For all $n,m\in {\Bbb Z}_{+}$, the functions $\la R(z) e_{n}, e_{m} \ra$ are  continuous in $z$ up to the cut along $[-1,1]$ except, possibly, the points $\pm 1$.
 
  \item
Estimates
\[
|\la R (z)e_{n}, e_{m}\ra |\leq C | \Omega(z)|^{-1} \big|q_{m}(z)/q_{n}(z)\big|\leq C_{1}<\infty, \q n\leq m,
\]
are true  with some positive constants that do not depend on $n$, $m$ and on $z$ in compact subsets of the closure of the set $ {\Bbb C}\setminus [-1,1]$  as long as they are away from the points $\pm 1$.
 \end{enumerate}
 \end{theorem}
 
  The statement (ii)   is known as the limiting absorption principle. It implies that
the spectrum of the operator $J$ on the interval $(-1,1)$ is absolutely continuous.  Using   the Cauchy-Stieltjes-Privalov formula \e{eq:Priv}, we also see that
matrix elements  $\la E(\lambda) e_{n}, e_{m}\ra$ of the spectral projector $E(\lambda)$ of the operator $J$  are  continuously differentiable in $\lambda\in (-1,1)$. 

 Note that the points $1$ and $-1$ may be eigenvalues of $J$; see Example~4.15 in \cite{Y/LD}.
 
 A calculation of the spectral family  $d E(\lambda)$ of the   operator $J $ is quite similar to the case of increasing coefficients  $a_{n}$ so that we again have representation
 \e{eq:RE1}. 
 Combining this representation with  formula \e{eq:HH4l} for $P_m (\lambda) $, we obtain    the following result.  We recall that the spectral measure $d\rho(\lambda)$ of a Jacobi operator $J$ is defined by relation \e{eq:UD1}; if this measure is absolutely continuous, we define the weight $\tau (\lambda)$ by equality \e{eq:Jac}.

 \begin{theorem}\label{SFst}
  Let assumptions \e{eq:comp} and \e{eq:LR}  hold. Then, for all $n,m\in{\Bbb Z}_{+}$ and $\lambda\in(-1,1)$, we have the representation
 \begin{equation}
 \frac{d\la E (\lambda)e_n, e_m\ra} {d\lambda}= (2\pi)^{-1}\sqrt{1-\lambda^2}  | \Omega (\lambda\pm i0) |^{-2} P_{n} (\lambda) P_m (\lambda)   .
\label{eq:EE}\end{equation}
 In particular,
the spectral measure of the operator $J$  is absolutely continuous on the interval $(-1,1)$ and the corresponding weight equals
 \begin{equation}
\tau (\lambda)=  (2\pi)^{-1}\sqrt{1-\lambda^2}  \, | \Omega ( \lambda\pm i0) |^{-2}     
\label{eq:SF1S}\end{equation}
$($the right-hand sides here do not depend on the sign$)$.
 \end{theorem} 
 

 Putting together Theorem~\ref{HX} and formula \e{eq:SF1S}, we arrive at the next result.
 
  \begin{theorem}\label{SFrS}
  Under assumptions \e{eq:comp} and \e{eq:LR} the weight $\tau (\lambda)$ is a continuous strictly positive function of $\lambda\in (-1,1)$.
 \end{theorem}

 Note that this result was earlier obtained in \cite{M-N} by specific methods of the orthogonal polynomials theory.
 
 Theorems~\ref{SFst} and \ref{SFrS} are of course quite similar to  Theorems~\ref{SF} and \ref{SFr} for the case $a_{n}\to\infty$.  The difference is that now  the factor $\sqrt{1-\beta_{\infty}^{2}}$ is  replaced by $\sqrt{1-\lambda^{2}}$ and we have the restriction $\lambda\in (-1,1)$.

According to \e{eq:WRfr} for the operator $J_{0}$, we have
 \[
 \Omega_{0}(\lambda\pm i0)=-2^{-1} (\lambda\pm i\sqrt{1-\lambda^2}),
 \]
 and hence expression  \e{eq:SF1S}  reduces to \e{eq:fr}.


In view of   \e{eq:SF1S} the amplitude factor in  \e{eq:Szold}   equals
     \begin{equation}
 \kappa (\lambda) =  (2 \pi)^{-1/2}  (1-\lambda^2)^{1/4} \tau (\lambda)^{-1/2}  . 
\label{eq:Sz1S}\end{equation}
Substituting this expression into \e{eq:Szold}, we can reformulate Theorem~\ref{Sz}
in a form more common for the orthogonal polynomials literature.

 \begin{theorem}\label{SW}
 Under the assumptions of Theorem~\ref{Sz}
  the  polynomials $P_{n} (\lambda)$ have asymptotics 
   \begin{equation}
 P_{n} (\lambda)=   (2/\pi)^{1/2}  (1-\lambda^2)^{-1/4} \tau (\lambda)^{-1/2}  \sin (\varphi_{n}  (\lambda) + \eta(\lambda) ) + o(1)
\label{eq:SWS}\end{equation}
as $n\to\infty$. Relation \e{eq:SWS} is uniform in $\lambda$ on compact subintervals of $(-1,1)$.
 \end{theorem}
 
 

  \subsection{Discussion}
  
  As was already mentioned in  Sect.~1.3,  under assumptions \e{eq:comp},  \e{eq:LR} asymptotic formulas \e{eq:GEGEy} and  \e{eq:Szold} for the orthonormal polynomials were first  obtained in paper \cite{Mate}.  However, expressions for the coefficients in the right-hand sides were not, at least in the author's opinion,  very efficient.   It was conjectured in \cite{Mate} that the asymptotic coefficient in  \e{eq:Szold} can be obtained from that in \e{eq:GEGEy} as the limit on $(-1,1)$ from complex values of $z$. This conjecture was later justified in \cite{Va-As}. In our approach this problem does not even arise since both coefficients are expressed in terms of the Wronskian $\{ P(z), f(z)\}$ of the polynomial and Jost solutions of the Jacobi equation  \e{eq:Jy}.

 As far as spectral results are concerned, we   note that a large part of
  Theorem~\ref{ACst} can also be obtained by the Mourre method \cite{Mo1}. It was applied to Jacobi operators in \cite{BdM}; to be precise, the problem in the space $\ell^2 ({\Bbb Z})$ was considered in \cite{BdM}, but this is of no importance. However,  the Mourre method does not exclude eigenvalues of $J$ embedded in its continuous spectrum. It only shows that  these eigenvalues  do not have other points of accumulation
  except $1$ and $-1$.  The Mourre method applies also to some Jacobi operators with increasing coefficients; see \cite{Sah}. Note   that very general conditions of the absolute continuity of spectrum were obtained in \cite{Stolz} by the subordinacy method of \cite{GD}.
  
  We also note papers  \cite{Gon, Nik} where Bernstein-Szeg\H{o} results (see Sect.~1.3) were extended to measures with a finite number of point masses away from the interval $[ -1,1]$.

  \subsection{Hilbert-Schmidt perturbations}

 In addition to   \e{eq:LR}, assume now that the condition 
 \begin{equation}
 \sum_{n=0}^\infty  (v_{n}^2 + b_{n}^2 )<\infty, \q v_{n}=a_{n}-1/2,
\label{eq:HS}\end{equation}
 is satisfied,
 that is, $V=J-J_{0}$ is a Hilbert-Schmidt operator. Then  the asymptotic formulas  of Theorems~\ref{GE1} and \ref{Sz} can be made more explicit. We proceed from the following elementary assertion.
 
 \begin{lemma}[\cite{JLR}, Lemma~4.8]\label{HSq}
 Let   $z\neq \pm 1 $. 
 Under assumption \e{eq:HS} there exists a finite limit
  \[
\lim_{n\to\infty}  \Big(\z(z)^{-n} \exp \big(-\frac{1}{\sqrt{z^2-1}}\sum_{m=0}^{n-1}(2 z v_{m} + b_m)  \big)q_{n}(z)\Big)\neq 0.
\]
\end{lemma}

 Thus, the next statement is a direct consequence of Theorem~\ref{GE1}.

 \begin{theorem}\label{SzHSq}
  Let assumptions     \e{eq:LR} and  \e{eq:HS} be satisfied, and  let $z\in  {\Bbb C}\setminus [-1,1] $.  Then there exist  finite limits
  \begin{equation}
\lim_{n\to\infty}  \Big(\z(z)^{n} \exp \big(\frac{1}{\sqrt{z^2-1}}\sum_{m=0}^{n-1}(2 z v_{m} + b_m)  \big)P_{n}(z)\Big)\neq 0
\label{eq:HSq3}\end{equation}
if $z$ is not an eigenvalue of the  operator $J$ and
 \begin{equation}
\lim_{n\to\infty}  \Big(\z(z)^{-n} \exp \big(-\frac{1}{\sqrt{z^2-1}}\sum_{m=0}^{n-1}(2 zv_{m} + b_m)  \big)P_{n}(z)\Big)\neq 0
\label{eq:HSp3}\end{equation}
if $z$ is   an eigenvalue of   $J$.
 \end{theorem}
 
  \begin{corollary}\label{SzHSq1}
  Suppose additionally that the conditions
    \begin{equation}
\sum_{n=0}^\infty v_{n}< \infty \q \mbox{and} \q\sum_{n=0}^\infty b_{n}< \infty 
\label{eq:HSq4}\end{equation} 
$($these series should be convergent but perhaps not absolutely$)$ are satisfied. Then the exponential factors in  \e{eq:HSq3} and  \e{eq:HSp3} may be omitted.
\end{corollary}

 

  In some cases the exponential factors in  \e{eq:HSq3} and  \e{eq:HSp3} can be  simplified.
    
\begin{example}\label{HSe}
Let  the conditions
\begin{equation}
a_{n} = 1/2+ v n^{-r_{1}}+ \tilde{v}_{n} , \q b_{n} = b n^{-r_{2}} +  \tilde{b}_{n}
\label{eq:LR1}\end{equation}
be satisfied with some   $r_{1}, r_{2}\in (1/2,1)$. 
Then 
  \begin{multline}
\sum_{m=0}^n (2z\, v_{m} + b_m)=2 z v(1-r_{1})^{-1}n^{1- r_{1}}+ b
(1-r_{2})^{-1}n^{1- r_{2}}
\\
+ 2 z  v \boldsymbol{\gamma}_{r_{1}}+
b \boldsymbol{\gamma}_{r_{2}} + \sum_{m=0}^\infty (2 z  \tilde{v}_m+ \tilde{b}_m)+  o(1)
\label{eq:HSp}\end{multline}
where $\boldsymbol{\gamma}_r-(1-r)^{-1}$  is the Euler-Mascheroni constant.
With a natural modification, expression \e{eq:HSp} remains true if $r_{j} =1$ for one or both $j$. In this case
$(1-r_{j})^{-1}n^{1-r_{j}}$ should be replaced by $\ln n$ and  $\boldsymbol{\gamma}_1$ is    the Euler-Mascheroni constant.
  \end{example}

Asymptotic formula \e{eq:Szold} can also be simplified. Similarly to Lemma~\ref{HSq}, we have

\begin{lemma}[\cite{JLR}, Lemma~4.12]\label{HS}
 Under assumption \e{eq:HS} there exists a finite limit
  \[
\lim_{n\to\infty}  \big(\varphi_n (\lambda) - n\theta-(\sin\theta)^{-1}\sum_{m=0}^{n-1} (2\cos \theta\, v_{m} + b_m)\big)=: \gamma (\lambda) 
\]
where $  \lambda=\cos\theta\in (-1,1)$.
\end{lemma}

Thus,  the following statement is a direct consequence of Theorem~\ref{Sz}.

 \begin{theorem}\label{SzHS}
  Let assumptions     \e{eq:LR} and  \e{eq:HS} be satisfied. Then
for  $\lambda \in (-1,1)$, the asymptotic formula 
   \begin{equation}
 P_{n} (\lambda)=  \frac{2 \kappa ( \lambda)} {\sqrt{1-\lambda^2}}  \sin \big(n\theta + (\sin\theta)^{-1}\sum_{m=0}^{n-1} (2\cos \theta\, v_{m} + b_m) +\gamma (\lambda)+ \eta(\lambda)\big) + o(1)  
\label{eq:SzHS}\end{equation}
holds as $n\to\infty$. Relation \e{eq:SzHS} is uniform in $\lambda$ on compact subintervals of $(-1,1)$.
 \end{theorem}
 
   Under assumption \e{eq:LR1} the phase in \e{eq:SzHS} can be simplified if one takes relation \e{eq:HSp} where $z = \cos\theta$ into account.
  
  Of course   formulas \e{eq:HSq3} and \e{eq:SzHS} are consistent with asymptotic  formulas   for Pollaczek
  polynomials in the Appendix in the book \cite{Sz}. 
  
    \subsection{Related research}

Without any additional assumptions, Hilbert-Schmidt perturbations $V$ of the operator $J_{0}$
were investigated in the deep papers \cite{KS} and \cite{DS}.  In   \cite{KS}, necessary and sufficient conditions in terms of the spectral measure $d\rho(\lambda)$ of the operator $J=J_{0} + V$ were found for $V$ to be in the   Hilbert-Schmidt class.  Asymptotic behavior of the corresponding polynomials $P_{n} (z)$  was studied in \cite{DS}. It was proved in   Theorem~5.1 of this paper  that
  the   limit of $ \z(z)^{n} P_{n}  (z)$ as $n\to\infty$ exists 
if and only if  conditions \e{eq:HS} and \e{eq:HSq4} are satisfied.
As shown in   Theorem~8.1 of \cite{DS1},    assumptions \e{eq:HS},  \e{eq:HSq4} are sufficient also for the validity of  formula  \e{eq:SzHS} but only in some {\it averaged}  sense. Such a regularization seems to be necessary since under  these  assumptions    
 the structure of the essential spectrum of the operator $J$ can be quite wild.  

 Condition \e{eq:LR} accepted in this paper is different in nature from \e{eq:HS},  \e{eq:HSq4}. On the one hand, it excludes too strong oscillations of the coefficients $v_{n}$, $b_{n}$ but, on the other hand, it permits their arbitrary slow decay as $n\to\infty$.

  Weights  $\tau (\lambda)$ with  singularities or zeros inside $(-1,1)$
   were investigated in the papers \cite{Nev, Van, M-F}. In the first of them even weights behaving like $\kappa |\lambda|^\nu$ where $\nu  >-1$, $\kappa>0$ as $\lambda\to 0$ were considered. Such weights are either singular at the point $\lambda=0$ if $\nu <0$ or, on the  contrary,  $\tau (0)=0$ if $\nu >0$. It was shown in \cite{Nev} that the corresponding Jacobi coefficients $a_{n}$ satisfy the asymptotic relation
   \[
   a_{n} =1/2 + (-1)^n \nu  / (4n) + o(1/n)
   \]
   (the coefficients $b_{n}=0$ if the weight $\tau (\lambda)$ is even).
   Since $|a_{n}'|\sim |\nu |/ (2n)$, the condition \e{eq:LR} is now violated.   For such weights, the asymptotic behavior of the polynomials $P_{n}  (\lambda)$ in a neighborhood of the point $\lambda=0$ differs from \e{eq:Szold}.
   More general results of this type were obtained in \cite{Van} where weights had several exceptional points.
   The results of \cite{M-F}  for weights with a jump singularity are morally similar. 
   
   Thus, condition  \e{eq:LR} 
   is practically necessary even for our   results on the weight $\tau (\lambda)$. 
  We also note that   asymptotics \e{eq:Szold} obtained under assumption \e{eq:LR}  is, in some sense, more regular than  the asymptotics of $P_{n}  (\lambda)$ in \cite{Nev, Van, M-F}.

        \section{Short-range perturbations}
        
                \subsection{Asymptotic formulas}
        
        Under the short-range assumption  \e{eq:Tr} the constructions of the previous section can be made more explicit. Instead of \e{eq:A22} the Jost solutions are now distinguished by the asymptotics
               \begin{equation}
 f_{n} (z)= \z (z)^n \big(1+ o(1)\big) , \q  z\in {\Bbb C}\setminus [-1,1],
\label{eq:jost}\end{equation}
as $n\to\infty$. Vectors $f(z)=(f_{n}(z))$ are defined as solutions
 of  the discrete Volterra integral equation
    \begin{equation}
  f_{n}(z)=\z(z)^n -\frac{1}{\sqrt{ z^2-1}}\sum_{m=n+1}^\infty (\z(z)^{n-m}- \z(z)^{m-n}) (Vf(z))_{m}.
\label{eq:V}\end{equation}
The perturbation  $V=J-J_{0}$ acts on vectors 
     $u= (u_{0}, u_{1}, \ldots)^\top=: (u_{n})$ by the formula 
\[
(V u) _{0} =  b_{0} u_{0}+ v_{0} u_{1}, \q  
( Vu) _{n} = v_{n-1} u_{n-1}+b_{n} u_{n}+ v_{n} u_{n+1} \q \mbox{for}\q n\geq 1  , 
 \]
 where $ v_{n}=a_{n}-1/2$. As usual, solutions of equation \e{eq:V} are constructed by iterations.
  Note that in contrast to the general case there is now the canonical choice of the Ansatz ${\cal A}_{n} (z)=\z (z)^n$ which works for all $z\in \clos \Pi_{0}\setminus\{-1,1\}$ so that there is no need to use the local arguments of Sect.~8.2. 

Let us state a  relevant  particular case of Theorem~\ref{EIK}.

\begin{theorem}\label{Jost}
Let  assumption   \e{eq:Tr} be satisfied, and let $z\in\clos{\Pi}_{0}$, $z\neq \pm 1$. Then   equation \e{eq:Jy} has a   solution satisfying  condition
\e{eq:jost}
as $n\to\infty$. Every function $f_{n}(z)$, $n=-1,0,1,\ldots$, depends analytically on $z\in\Pi_{0}$, and it is continuous in $z$ up to the cut along $[-1,1]$ except, possibly, the points $\pm 1$.
\end{theorem}

\begin{corollary}\label{JoSR}
If $\lambda\in (-1,1)$, then, as $n\to\infty$,  
  \[
 f_{n} (\lambda\pm i0)= e^{\mp i n\theta} \big(1+ o(1)\big) , \q \theta=\arccos\lambda\in (0,\pi)  .
\]
\end{corollary}

With this definition of the Jost solution, the Jost function $\Omega (z)$ is, as before,  defined as the Wronskian  $\Omega (z)=\{ P(z),f(z)\}$.  
We note a representation for $\Omega (z)$  in terms of the orthonormal polynomials:
   \begin{equation}
 {\Omega} (z) = -\frac{1}{2\z (z)} +  \sum_{n=0}^{\infty} \z(z)^{n} (VP(z))_{n} ,\q z \in\clos\Pi, \q z\neq \pm 1,
 \label{eq:OM}\end{equation}
see Proposition~3.3 in \cite{Y/LD} for the proof. A similar representation in the continuous case is quite standard; see, e.g., Lemma~1.12 in Chapter~4 of \cite{YA}.
The Jost function $\Omega_{0} (z)$ for the free operator $J_{0}$ is given by equality \e{eq:WRfr}.

Relation \e{eq:HH4l} remains true, and a   short-range version of Theorem~\ref{Sz} can be stated in the following way.

\begin{theorem}\label{SzSR}
Let the amplitude $\kappa (\lambda)$ and the phase $\eta (\lambda)$ be defined by formulas \e{eq:AP+}.
  Under assumption \e{eq:Tr}  the  polynomials $P_{n} (\lambda)$ have asymptotics 
   \begin{equation}
 P_{n} (\lambda)= \frac{2 \kappa (\lambda) }{\sqrt{1-\lambda^{2}}}  \sin (n\arccos \lambda + \eta(\lambda) ) + o (1),\q \lambda \in (-1,1),
\label{eq:Szsr}\end{equation}
as $n\to\infty$. Relation \e{eq:Szsr} is uniform in $\lambda$ on compact subintervals of $(-1,1)$.
 \end{theorem} 
 
 In the case $z\in{\Bbb C}\setminus [-1,1]$  asymptotic  formulas  \e{eq:GEGEy} and \e{eq:GEGEx} for $P_n(z)$ are true with $q_{n} (z)$ replaced by $\z(z)^n$.
        
   Representations \e{eq:Rrpm} for the resolvent and      \e{eq:EE} for the spectral measure remain unchanged. We note also an equation for the orthonormal polynomials
     \begin{equation}
 P_{n}(z)=P_{n}^{(0)}(z) +\frac{1}{\sqrt{ z^2-1}}\sum_{m=0}^{ n-1} (\z(z)^{n-m}- \z(z)^{m-n}) (VP(z))_{m},\q n\geq 1,
\label{eq:OPR1}\end{equation}
where 
 \[
P_{n}^{(0)} (z)= \frac{1}{2\sqrt{z^2-1}} \big(\z(z)^{-n-1} -\z(z)^{n+1}\big)
\]
are the normalized Chebyshev polynomials of the second kind.



        
             \subsection{The perturbation determinant  and   the spectral shift function}
             
             First we recall abstract definitions of these notions for arbitrary bounded self-adjoint operators $J_{0}$ and $J$ with a trace class difference $V=J-J_{0}$. We refer to  the books \cite{GK, Sim,Ya} for a consistent presentation of this theory.  In view of our applications, we suppose that the spectrum of the  operator $J_{0}$ coincides with the interval $[-1,1]$.

If $V\in{\goth S}_{1}$ (the trace class), then 
 the perturbation  determinant
   \begin{equation}
D (z):=\det \big(I+V R_{0} (z)\big)
\label{eq:Tr1}\end{equation} 
for the pair $J_{0}$, $J$
 is well defined and is an analytic function of $z\in{\Bbb C}\setminus [-1,1]$. 
 Obviously, $D (\bar{z})=\overline{ D (z)}$ and
  \begin{equation}
D (z)\to 1\q {\rm as}\q |z|\to\infty .
\label{eq:PDhe}\end{equation}
Note also  the general formula
  \begin{equation}
\tr \big( R(z)-R_{0}(z)\big)=- \frac{   D' (z)  }{ D(z)   }.
\label{eq:PD}\end{equation}

  The  Kre\u{\i}n spectral shift function $\xi (\lambda)$ 
   is defined in terms of the perturbation determinant  \e{eq:Tr1}. According to \e{eq:PDhe} we can fix the branch of the function $\ln D (z) $ for $\Im z\neq 0$ by the condition
  \[
 \arg D (z)\to  0 \q {\rm as}\q |z|\to\infty.
\]
Then
   \begin{equation}
  \xi (\lambda):=\pi^{-1}\lim_{\varepsilon\to+0}\arg  D (\lambda+i\varepsilon).
\label{eq:SSF}\end{equation}
In the abstract setting,
this limit exists for almost every $ \lambda\in{\Bbb R}$,
\begin{equation}
 \int_{-\infty}^\infty |\xi(\lambda) |d\lambda \leq \| V\|_{1}
\label{eq:DS}\end{equation}
and the representation  
\begin{equation}
\ln D (z) =\int_{-\infty}^\infty \xi(\lambda) (\lambda-z)^{-1} d\lambda, \q \Im z\neq 0,
\label{eq:DD}\end{equation}  
holds.  The    function $\xi (\lambda)$ is constant on subintervals of $(-\infty,-1)$ and $(1,\infty)$
 not containing eigenvalues of $  J$. In particular,  $\xi(\lambda)=0$ for $\lambda$ below the smallest and above the largest eigenvalue of the operator $J$.  If $\mu$ is an isolated simple  eigenvalue   of the operator $J$, then 
\begin{equation}
\xi (\mu +0)- \xi (\mu -0)=-1.
\label{eq:ES}\end{equation}

Let us come back to Jacobi operators $J_{0}$ and $J$.

\begin{lemma}\label{Tr-SR}
The difference $V=J-J_{0}$ belongs to the trace class ${\goth S}_{1}$ if and only if assumption \e{eq:Tr} is satisfied.
 \end{lemma} 
 
 \begin{pf}
 Let $T$ be the shift in the space $\ell^{2} ({\Bbb Z}_+)$ defined by $(Te)_{n} =e_{n+1}$ and let
 $A=\diag\{a_{n}\}$,  $B=\diag\{b_{n}\}$.  Under assumption \e{eq:Tr} the diagonal operators $A-I/2$ and $B$ are trace class. 
 Therefore the same is true for the operator $V=A-I/2+ BT+ BT^{*}$. Conversely, if $V\in{\goth S}_{1}$, then (see, e.g., Theorem~11.2.3 in the book~\cite{BS})
 \[
 \sum_{n=0}^{\infty}| \la V e_{n}, e_{n}\ra |+  \sum_{n=0}^{\infty}| \la V e_{n}, e_{n+1}\ra |<\infty.
 \]
 Since $ \la V e_{n}, e_{n}\ra = a_{n}-1/2$ and $ \la V e_{n}, e_{n+1}\ra=b_{n}$, this proves  \e{eq:Tr}.
  \end{pf}
 
 Let us express the perturbation determinant  and  the spectral shift function
in terms of the Jost function $\Omega (z)$ and the phase $\eta(\lambda)$. It is  convenient to introduce the corresponding normalized   objects by the relations
 \[
\boldsymbol{\Omega} (z) : =\frac{\Omega (z)}{\Omega_{0} (z)}=-2 \z(z)\Omega (z), \q z\in \Pi_{0},
\]
and
\[
\boldsymbol{\eta} (\lambda): =\eta (\lambda)-\arccos\lambda, \q \lambda\in  (-1,1).
\]

Using    representation \e{eq:Rrpm} for the resolvent, one can prove  (see Theorem~5.4 in \cite{Y/LD}) a relation 
 \[
\tr \big( R(z)-R_{0}(z)\big)=- \frac{ \boldsymbol{\Omega}' (z)  }{  \boldsymbol{\Omega} (z)   }.
\]
Comparing it with formula  \e{eq:PD}, we see that
   \begin{equation}
D (z)=   A \, \boldsymbol{\Omega}  (z), \q z\in \Pi_{0} ,
\label{eq:HH}\end{equation}
for some constant $A\in{\Bbb C}$.   According to Theorem~5.6  in \cite{Y/LD} the normalized    Jost function has an asymptotics 
 \begin{equation}
\boldsymbol{\Omega}(z)= A^{-1} + O (|z|^{-1}),\q z\to\infty,
\label{eq:AA-}\end{equation}
where  
 \begin{equation}
A=\prod_{k=0}^{\infty} (2a_{k}) .
\label{eq:AA}\end{equation}
Note that   under assumption \e{eq:Tr} the infinite product here  converges and $A>0$.  Putting together  asymptotic relations \e{eq:PDhe} and  \e{eq:AA-}, we  obtain the next result.

 \begin{theorem}\label{RRY}
  Under assumption \e{eq:Tr}  equality \e{eq:HH}
 is  true with constant \e{eq:AA}.
 \end{theorem}

It follows from Theorem~\ref{Jost} and formula 
\e{eq:HH}  that the perturbation determinant $ D (z)$ is a continuous function of $z\in\clos\Pi_{0}$ except, possibly, the points $z=\pm 1$.  Moreover, according to \e{eq:HH5}, we have $D (\lambda\pm i0)  \neq 0$.  Therefore the spectral shift function 
$\xi (\lambda)$ is    a continuous function of $\lambda\in (-1,1)$.  

 Comparing  definitions    \e{eq:AP+} and \e{eq:SSF}, we find
 a link between     $\xi (\lambda)$ and the scattering  phase $\eta(\theta)$.

\begin{theorem}\label{SzShift}
  Under assumption \e{eq:Tr},  
 the relation
  \begin{equation}
  \xi (\lambda) =\pi^{-1}\eta (\arccos\lambda) 
\label{eq:SSF3}\end{equation}
holds   for all $\lambda\in ( -1 , 1)$.
 \end{theorem}

Substituting \e{eq:SSF3} into \e{eq:Sz} and taking into account relation  \e{eq:Sz1S}, we can reformulate Theorem~\ref{Sz} in terms 
of the weight  $\tau (\lambda)$ and   the  spectral shift function $\xi(\lambda)$ (cf. \e{eq:SW}).  This yields  the asymptotic formula as $n\to\infty$:
 \[
 P_{n} (\lambda)= (2/\pi)^{1/2}  (1-\lambda^2)^{-1/4} \tau (\lambda)^{-1/2} \sin ( (n+1)\arccos\lambda +\pi \xi(\lambda) ) + o(1)  ,  \q \lambda\in (-1,1),
\]
We emphasize that $\eta(\theta)$  is a continuous function of $\theta\in (0,\pi)$, but Theorem~\ref{Jost} yields no information about its behavior as $\theta\to 0$ and $\theta\to \pi$. Comparing relations \e{eq:DS} and \e{eq:SSF3}, we however see that
 \[
\int_{0}^\pi |\eta (\theta)| \sin\theta d\theta \leq\pi  \| V\|_{1}.
\]
 
    \subsection{Threshold behavior}
    
     A study of  the perturbation determinant in the limits  $z\to 1$ and $z\to -1$ is similar to the same problem for the Schr\"odinger operators in the space $L^{2}({\Bbb R}_{+})$ at zero energy discussed in Sect.~4.3 of \cite{YA}.  This problem  was already considered in Sect.~4 of \cite{Y/LD}, but here our presentation is closer to  \cite{YA}. 
     
     Now
      an additional assumption on the coefficients $a_{n}$, $b_{n}$ is required. 
      
\begin{theorem}[\cite{Y/LD}, Theorem~4.1] \label{JOSF}
Suppose that 
 \begin{equation}
\sum_{n=0}^\infty n(|  a_{n}-1/2|+|b_{n}|)<\infty. 
\label{eq:Trfd}\end{equation} 
Then all functions $f_{n} (z)$  are continuous   as $z\to\pm 1$, the sequence  $f_{n}(\pm 1)$ satisfies  the equation  
  \begin{equation}
  a_{n-1} f_{n-1} (\pm 1) +b_{n} f_{n} (\pm 1) + a_{n} f_{n+1} (\pm 1)= \pm f_{n} (\pm 1), \q n\in{\Bbb Z}_{+},
 \label{eq:Jypm}\end{equation}
 and  $f_{n}(\pm 1)=(\pm 1)^n +o(1)$ as $n\to\infty$.
\end{theorem}

In particular,   the Jost function $\Omega(z)$  is continuous in $z$ up to the cut along $[-1,1]$, including the points $\pm 1$.
   Representation \e{eq:Rrpm} for the resolvent $R(z)$ implies that    its matrix elements $\la R(z) e_{n}, e_{m}\ra$ are also continuous up to the point $\pm1$ provided 
    $\Omega (\pm 1) \neq 0$.
    
    Using   Theorem~\ref{GE}, we  introduce a solution  $g_{n} (\pm 1)$ of  equation
\e{eq:Jypm} by the formula
 \[ 
g_{n }(\pm 1)= f_{n}(\pm 1)\sum_{m=N_{0} + 1}^n (a_{m-1} f_{m-1}(\pm 1) f_{m}(\pm 1))^{-1},\q n \geq N_{0} + 1,
\]
 The following result is a  direct consequence of Theorem~ \ref{JOSF}.


 \begin{lemma}\label{GEpm}
Let  assumption \e{eq:Trfd} hold.  Then the sequence $g_{n} (\pm 1)$ satisfies equation \e{eq:Jypm}, 
\[
g_{n} (\pm 1)= 2(\pm 1)^{n +1}n (1+o(1))
\]
as $n\to\infty$ and $\{f(\pm 1),
g(\pm 1)\}=1$. 
\end{lemma}

Since neither of linearly independent solutions  $f_{n} (\pm 1)$,  $g_{n} (\pm 1)$ nor their linear combinations   tend to zero as $n\to \infty$, we obtain

\begin{theorem}\label{JOSF2}
Under  assumption \e{eq:Trfd} equations \e{eq:Jypm} do not have solutions tending to zero as $n\to \infty$. In particular, the operator $J$ cannot have eigenvalues $1$ and $-1$.
\end{theorem}

Next, we discuss an  asymptotic behavior of the orthonormal polynomials $P_{n} (z)$ at the critical points $z=\pm 1$. In view of  relation \e{eq:Pf}   where $z=\pm 1$, the next result is a  direct consequence of Theorem~\ref{JOSF} and Lemma~\ref{GEpm}.

\begin{theorem}\label{AScr}
Under  assumption \e{eq:Trfd} we have
\begin{equation}
P_{n} (\pm 1)= -2 \Omega (\pm 1) (\pm 1)^{n +1 } n +o(n).
\label{eq:P+}\end{equation}
\end{theorem}

Passing in  \e{eq:OPR1} to the limit $z\to \pm 1$
and taking into account that $P_{n}^{(0)}(\pm 1)= (n+1) (\pm 1)^{n}$, we obtain an equation  
 \begin{equation}
(\pm 1)^n   P_{n}(\pm 1)= n+1  -2 \sum_{m=0}^{n-1} (\pm 1)^{m+1}( n-m)  (V P(\pm 1))_{m}.
\label{eq:Vpmr}\end{equation}
Similarly, in the limit $z\to \pm 1$ relation \e{eq:OM} yields
   \begin{equation}
\Omega(\pm 1)=\mp  2^{-1} + \sum_{n=0}^{\infty} (\pm 1)^{n} (VP(\pm 1))_{n}  .
\label{eq:OM+}\end{equation}

    Our next goal is to  consider the exceptional case $\Omega (\pm 1) = 0$. Let us use the same terminology as for the Schr\"odinger operators.    

   \begin{definition}\label{resin}
Let  assumption \e{eq:Trfd} hold. If $\Omega (\pm 1) = 0$, we say that the Jacobi operator $J$ has a threshold resonance at the point $z=\pm 1$.
\end{definition}

Clearly, the condition  $\Omega (\pm 1) = 0$ is equivalent to the linear dependence of the solutions $P_{n} (\pm 1)$ and
$f_{n} (\pm 1)$ of equation \e{eq:Jypm} whence 
 \begin{equation}
P_{n} (\pm 1)=  f_{n} (\pm 1)/f_0 (\pm 1), \q n\in{\Bbb Z}_{+}.
\label{eq:Vxa}\end{equation}
Note that $f_0 (\pm 1)\neq 0$ if $\Omega (\pm 1) = 0$.    In this case Theorem~\ref{JOSF}  allows us to make  asymptotic formula  \e{eq:P+} more precise:
 \begin{equation}
P_{n}(\pm 1)=  (\pm 1)^{n} / f_0 (\pm 1)+ o (1), \q n\to\infty.
\label{eq:Vx1}\end{equation}
Conversely, if the sequence $P_{n}(\pm 1)$ is bounded, then it follows from \e{eq:Vpmr} that
 \begin{align}
 (\pm 1)^n   P_{n}(\pm 1)= &\Big(1- 2 \sum_{m=0}^\infty (\pm 1)^{m+1}   (V P(\pm 1))_{m}\Big) n
\nonumber \\ 
&+1+ 2\sum_{m=0}^\infty (\pm 1)^{m+1} m  (V P(\pm 1))_{m} + o(1).
\label{eq:Vpmrb}\end{align}
The coefficient at $n$ here is necessarily zero so that comparing \e{eq:Vx1} and \e{eq:Vpmrb}, we obtain a relation  
 \begin{equation}
 f_0 (\pm 1)^{-1}= 1+ 2 \sum_{m=0}^\infty (\pm 1)^{m+1} m  (V P(\pm 1))_{m} .
\label{eq:Vx2}\end{equation}
Thus, the definition of a threshold resonance can be equivalently reformulated in the following way.

\begin{lemma}\label{th}
Let  assumption \e{eq:Trfd} hold. Then the Jacobi operator $J$ has a threshold resonance at the point $z=\pm 1$ if and only if the polynomial solution
$P_{n} (\pm 1)$ of of equation \e{eq:Jypm}  is bounded for $n\to\infty$. In this case asymptotic formula \e{eq:Vx1} holds and $ f_0 (\pm 1)$ is given by relation \e{eq:Vx2}.
\end{lemma}

Now we find an asymptotic behavior of the Jost function $\Omega (z)  $ as $z\to \pm 1$. This requires 
an  estimate on the rate of convergence of $P_{n} (z)$ to $P_{n} (\pm 1)$.  The following technical assertion  is quite similar to Lemma~3.6 in \cite{YA} and its proof will be omitted.

\begin{lemma}\label{th1}
Let  assumption \e{eq:Trfd} hold and $\Omega (\pm 1) =0$. Then 
\[
| P_{n}(z)-P_{n}(\pm 1)| \leq C n  \sqrt{|z^{2}-1|} \, |\z(z)|^{-n}.
\]
\end{lemma}

The  next result is a translation to the discrete case of Proposition~3.7 in Chapter~4 of \cite{YA} stated there in the continuous framework.
 
  
\begin{theorem}\label{resin1}
Under  assumption \e{eq:Trfd} suppose that $\Omega (\pm 1) = 0$. Then
\begin{equation}
\Omega(z)=-2^{-1}  f_{0} (\pm 1)^{-1}  \sqrt{z^{2}-1} + o (\sqrt{|z^{2}-1|})
\label{eq:th1}\end{equation}
as $z\to \pm 1$.
\end{theorem}

\begin{pf}
Let us proceed from  representation \e{eq:OM} of the   Jost function   in terms of the orthonormal polynomials.
First, we observe that
\[
 \sum_{n=0}^{\infty}| \z(z)| ^{n} |(VP(z))_{n}- (VP(\pm 1))_{n}| = o (\sqrt{|z^{2}-1|}).
\]
Indeed, every term in this sum is $O( |z^{2}-1|)$, and, by the dominated convergence, the passage to the limit in the sum can be justified by Lemma~\ref{th1}.
Next, we observe that
\[
\z(z)^{n}=(\pm 1)^{n }-n(\pm 1)^{n -1}\sqrt{z^{2}-1} + O (z^{2}-1), 
\]
so that
\begin{align}
\Omega(z)=\mp  2^{-1} -2^{-1} \sqrt{z^{2}-1} + & \sum_{n=0}^{\infty} (\pm 1)^{n } (VP(\pm 1))_{n}
\nonumber\\
-  \sum_{n=0}^{\infty} & n(\pm 1)^{n -1} (VP(\pm 1))_{n} \sqrt{z^{2}-1} + o (\sqrt{|z^{2}-1|}).
\label{eq:th5}\end{align}
In view of \e{eq:OM+} 
the constant term here is zero because $\Omega (\pm 1)=0$.  The coefficient at $\sqrt{z^{2}-1}$ is $ -2^{-1} f_{0} (\pm 1)^{-1}$ by virtue of representation \e{eq:Vx2}.  Therefore \e{eq:th5} implies relation \e{eq:th1}.
 \end{pf}

Using formula   \e{eq:SF1S}, we obtain the following consequence for the weight function.

\begin{corollary}\label{resnw}
Under the  assumptions of Theorem~\ref{resin1}, we have
 \begin{equation}
\tau (\lambda)=2\pi^{-1} f_{0} (\pm 1)^{2}(1-\lambda^2)^{-1/2} (1+ o(1)), \q \lambda\in (-1,1),
\label{eq:th6}\end{equation}
as $\lambda\to \pm 1$.
\end{corollary}

 In view of Proposition~\ref{res} and equality \e{eq:Vxa}, Theorem~\ref{resin1} implies also the following result.

\begin{corollary}\label{resn1}
Under the  assumptions of Theorem~\ref{resin1}  for all $n,m\in{\Bbb Z}_{+}$, the representation 
 \[
 \la R(z)e_{n}, e_{m}\ra = -2\frac{ f_{n} (\pm 1) f_{m} (\pm 1) + o(1)}
{\sqrt{z^{2}-1}  }
\]
as $z\to \pm 1$ is satisfied.
\end{corollary}

\begin{remark}\label{find}
Under  assumption \e{eq:Trfd} 
  Theorems~\ref{JOSF} and \ref{resin1} ensure that $\Omega (z)\neq 0$ for $z\in (1,1+\epsilon)$ and $z\in (-1-\epsilon, -1)$ if $\epsilon>0 $ is small enough.
 Therefore the discrete spectrum of the operator $J$ is finite.
\end{remark}

  \begin{example}\label{find1}
  Let the Jacobi polynomials ${\sf G}^{(\alpha,\beta)}(z)$ be defined by spectral measure  \e{eq:Jac1}. 
    Asymptotics of ${\sf G}^{(\alpha,\beta)}(\pm 1)$ are given by formulas \e{eq:edge}.
   It follows from relations  \e{eq:norm6} that condition  \e{eq:Trfd} is satisfied if and only if $|\alpha|=|\beta|=1/2$.  The point $z=1$ (and the point $z=-1$) is regular for the corresponding Jacobi operator ${\sf J}^{(\alpha,\beta)}(z)$ if $\alpha=1/2$ (resp., $\beta=1/2$) and there is a resonance at this point if $\alpha =-1/2$ (resp., $\beta=-1/2$). 
  In the regular case formulas \e{eq:edge} are consistent with a general relation  \e{eq:P+}, and in the resonant case they are consistent with   \e{eq:Vx1}.
  
  Formulas \e{eq:edge} for the cases $|\alpha| \neq 1/2 $ or $|\beta| \neq 1/2 $ show that the asymptotics of  the orthonormal polynomials at the edge points is significantly   changed   if   condition \e{eq:Trfd} is even slightly relaxed.  
   \end{example}
   
   Let us finally consider the case when $a_{n}-1/2$ and $b_{n}$ decay as $n^{-\varrho}$ with $\varrho<2$. An important example of such coefficients is given by the Pollaczek polynomials when $\varrho= 1$; see Appendix~C. For  the Jacobi polynomials, we have $\varrho=2$.  For simplicity, we now assume that $a_{n}=1/2$ for all $n $ and, up to sufficiently rapidly decaying terms,    $b_{n}=\kappa n^{-\varrho}$ where $\varrho<2$ and,  for definiteness, $\kappa >0$.  In this case the operator $J$ has infinite discrete spectrum accumulating at the point $\lambda=1$. 
   Relying on an analogy with the continuous case considered in \cite{YS} (see also \S 4.3 in the book \cite{YA}), we conjecture that
    \begin{equation}
\ln \tau(\lambda)=-\gamma_{-1}(1+ \lambda)^{- (2-\varrho)/(2\varrho)}(1+ o(1))\q \mbox{as}\q \lambda\to -1+0
\label{eq:sl-d}\end{equation}
 and  
 \begin{equation}
  \tau(\lambda)=\gamma_{1}(1-\lambda)^{- 1/4} (1+ o(1))\q \mbox{as}\q \lambda\to 1-0
\label{eq:sl-d1}\end{equation}
where $\gamma_{-1}$ and $\gamma_{1} $ are positive constants depending on $\varrho$ and $\kappa$.

  According to \e{eq:sl-d} the weight $\tau(\lambda)$ tends to zero exponentially as $\lambda\to -1+0$. This result can be interpreted as a virtual shift to the right of the essential spectrum of the operator $J$ (obviously, it coincides with the interval $[-1+\kappa, 1+\kappa]$ if $\varrho=0$). The spectral point $\lambda=-1$ becomes quasiregular in this case. For the Schr\"odinger operator, this phenomenon in discussed in \cite{YS, YA}.
  
 It follows from  \e{eq:sl-d}  that   the integral
   \begin{equation}
\int_{-1}^1 (1-\lambda^2)^{s}  \ln \tau (\lambda)   d\lambda> -\infty 
\label{eq:Szeg}\end{equation}
if 
   \[
   s> (2-3\varrho) / (2\varrho).
\]
In particular, for all $\varrho >1$, we can take $s=-1/2$. This is consistent with the following conjecture of Nevai \cite{Nev2} proved in \cite{KS} with earlier partial results obtained in
\cite{Nev3} and \cite{Va-As1}.

  \begin{theorem}[\cite{KS}, Theorem~2]\label{Nev-Sz}
  Under assumption  \e{eq:Tr} condition \e{eq:Szeg} is satisfied for $s=-1/2$.   
     \end{theorem}
     
    Condition \e{eq:Szeg} for $s=-1/2$ is known as the Szeg\H{o} condition.  It means that the weight $\tau(\lambda)$ does not tend to zero too rapidly as $\lambda\to\pm 1$.
    
    Under the Hilbert-Schmidt assumption \e{eq:HS} a weaker  quasi-Szeg\H{o} condition holds.
    
     \begin{theorem}[\cite{KS}, Theorem~1]\label{KS-qS}
  Under assumption  \e{eq:HS} condition \e{eq:Szeg} is satisfied for $s=1/2$.   
     \end{theorem}
     
     In the case $\varrho>1/2$ assumption  \e{eq:HS} is of course satisfied. In the intermediary case $\varrho =1$ formula \e{eq:sl-d} is consistent with  relation 
     \e{eq:PolCou} for the Pollaczek polynomials.  This example shows that assumption  \e{eq:Tr} in Theorem~\ref{Nev-Sz} is very precise.

Let us now discuss  formula \e{eq:sl-d1}.  Recall that according to \e{eq:Szeg}  in the regular case $\Omega (1)\neq 0$, the weight function $\tau(\lambda)$ has a finite positive limit as $\lambda\to 1-0$.  According to \e{eq:th6} it has a singularity $\tau_{1} (1-\lambda)^{-1/2}$, $\tau_{1}>0$,  in the resonant case $\Omega (1) = 0$. Formula \e{eq:sl-d1} shows that  for slowly decaying diagonal elements $b_{n}\sim \kappa n^{-\varrho}$ where $\kappa>0$ and $\varrho<2$ (the elements $a_{n}=1/2$), the behavior of 
$\tau(\lambda)$ is intermediary between the regular and resonant cases. This can be interpreted as the existence of a weak but stable resonance for slowly decaying recurrence coefficients at
the threshold energy $\lambda=1$.

 
Of course in the case $\kappa <0$,  the results  are exactly the same, but the roles of the thresholds $\lambda=1$ and  $\lambda=-1$ are interchanged.

    \subsection{Szeg\H{o} function}
    
     We define 
  the Szeg\H{o} function $S(\z)$   by the formula
       \begin{equation}
S(\z)=  \exp\Big(\frac{1 }{4\pi}\int_{-\pi}^\pi  \frac{  e^{i\theta}+\z }{  e^{i\theta}-\z } \ln \big( \tau (\cos\theta)|\sin\theta|\big)d\theta\Big),\q |\z|< 1,
\label{eq:SzY}\end{equation}
where $\tau (\lambda)$ is the weight function  \e{eq:SFx}. This is exactly formula (10.2.10) (see also Theorem~12.1.2) in \cite{Sz}, but in contrast to \cite{Sz} we do not suppose that the Jacobi operator $J$ has no eigenvalues. Of course definition \e{eq:SzY} requires that  the   condition
    \begin{equation}
\int_{-\pi}^\pi | \ln \big( \tau (\cos\theta)|\sin\theta|\big)|d\theta<\infty
\label{eq:SzY1}\end{equation}
(or, equivalently, \e{eq:Szeg} for $s=-1/2$)
 be satisfied.

     Recall the standard Jensen-Poisson representation of analytic functions $f(\z)$ from the Hardy class ${\bf H}^{1}$ on the unit disc $|\z| <1$:
       \[
f(\z)= i \Im f(0) + \frac{1 }{2\pi}\int_{-\pi}^\pi  \frac{  e^{i\theta}+\z }{  e^{i\theta}-\z }   \Re f(e^{i\theta})d\theta .
   \]
 In particular, we  have
   \begin{equation}
   \begin{split}
    \ln   (1+ \z) &=
 \frac{1 }{2\pi}\int_{-\pi}^\pi  \frac{  e^{i\theta}+\z }{  e^{i\theta}-\z } \ln \cos (\theta/2) d\theta, 
\\
    \ln   (1- \z) &=
 \frac{1 }{2\pi}\int_{-\pi}^\pi  \frac{  e^{i\theta}+\z }{  e^{i\theta}-\z } \ln| \sin (\theta/2) |d\theta. 
    \end{split}
\label{eq:wwz}\end{equation}
Here  the branches of the functions $\ln   (1\pm \z)$ are fixed by the condition $\ln 1=0$.  We also recall that  the Nevanlinna class ${\bf N} $ of   functions $f(\z)$
  analytic  on the unit disc  is distinguished by the condition 
   \[
\sup_{r<1} \int_{0}^{2\pi}\ln^{+} |f(re^{i\theta})| d\theta <\infty
\]
(as usual $\ln^{+}a= \max\{\ln a,0\}$).


   Let $D (z)$ be the perturbation determinant  \e{eq:Tr1} and $\Delta(\z)= D (z)$ if $\z=\z(z)$. It follows from Corollary~\ref{JOST} that the function $\Delta(\z)$ is analytic on the unit disc, and it is continuous up to the unit circle  with a possible exception of the points $\pm 1$.
Moreover, according to  \e{eq:HH5},  $\Delta(\z)\neq 0$ if $|\z| =1$ but $\z\neq \pm 1$. 
  We also note that $\Delta(0)=1$ according to \e{eq:PDhe}.


      Let $\lambda_{k}$ be eigenvalues (lying on $(-\infty, -1)\cup(1,\infty)$) of the operator $J$. We suppose that  $|\lambda_{1}|\geq |\lambda_{2}|  \geq \cdots >1$  not distinguishing positive and negative eigenvalues in notation.
          The  numbers $\mu_{k}:=\z (\lambda_{k})\in (-1,1)$  are zeros of the function $\Delta(\z)$. It was shown in \cite{Hu-S} that under assumption \e{eq:Tr} 
     \begin{equation}
\sum_{k=1}^\infty (1-|\mu_{k}|)<\infty.
\label{eq:Bla1}\end{equation}
Using this result the inclusion ${\Delta}\in {\bf N}$ was established in \cite{KS}.   
 
 Let us define an outer function
      \begin{equation}
G(\z)= \exp\Big(\frac{1 }{2\pi}\int_{-\pi}^\pi  \frac{  e^{i\theta}+\z }{  e^{i\theta}-\z } \ln | \Delta(e^{i\theta})|d\theta\Big)
   \label{eq:Bla2}\end{equation}
   where the function $\ln |  \Delta(e^{i\theta})|$ belongs to $L^1 (-\pi,\pi)$ 
   because $ \Delta\in {\bf N}$. This is equivalent to the Szeg\H{o} condition \e{eq:SzY1} 
   since, by relations \e{eq:SF1}  and \e{eq:HH},  
\[
|  {\Delta}(e^{i\theta})|^{2}= |D(\cos\theta)|^{2}= A^{2}   \frac{2}{\pi}\frac{|\sin\theta|}{\tau (\cos\theta)},
\]
whence
 \begin{equation}
2 \ln |  {\Delta}(e^{i\theta})|=-  \ln (   \tau(\cos\theta)|\sin\theta|)  +\ln (2\pi^{-1} A^2 \sin^2\theta).
 \label{eq:wj}\end{equation}

Let us now compare definitions    \e{eq:SzY} and \e{eq:Bla2}.  
Substituting   expression \e{eq:wj} into \e{eq:SzY} and taking    formulas \e{eq:wwz} into account, we see that
 \begin{equation}
G (\z)= A  \frac{1-\z^2}{\sqrt{ 2\pi} S(\z)}
\label{eq:SG}\end{equation}

Next, we introduce
 the   Blaschke product
      \begin{equation}
B(\z) =\prod_{k=1}^\infty\frac{|\mu_{k}|}{\mu_k}\frac{\mu_{k}-\z}{1-\bar{\mu}_k\z},\q |\z|<1.
\label{eq:Bla}\end{equation}
According to condition \e{eq:Bla1} it is well defined. The function $B(\z)$ is continuous on the closed disc $|\z|\leq 1$ except, possibly, the points $\pm 1$ and 
 $|B(\z)|=1$ for $|\z| =1$. 

 It is shown  in \cite{KS} that  the function $\Delta(\z)$  does not have a singular inner component. 
Therefore the classical factorization (see, e.g., Theorem~2.9 in \cite{Duren}) reads as
    \[
  \Delta(\z)=      
B(\z)    G(\z)  .
 \]
Comparing this relation with    \e{eq:SG}, we arrive at the following result.

 \begin{theorem}\label{PertSze}
Let assumption   \e{eq:Tr} be satisfied, and let $|\z |< 1$. Let  $D(z)$ be the perturbation determinant  \e{eq:Tr1} and $\Delta(\z)= D(z)$. Define the Blaschke product  $B(\z)$ by formula  \e{eq:Bla}, the Szeg\H{o} function $S(\z)$ --  by \e{eq:SzY}, and the product $A$  --  by \e{eq:AA}.
Then the factorization
  \begin{equation}
 \Delta (\z)= A B(\z)\frac{1-\z^2}{\sqrt{ 2\pi} S(\z)}
\label{eq:SzD3}\end{equation}
holds.
\end{theorem}

We emphasize that Theorem~\ref{PertSze} is a direct consequence of classical results on the factorization of functions in the Nevanlinna class combined with the analytical results of \cite{KS}.

According to \e{eq:SzD3} formulas \e{eq:DD} and \e{eq:SzY}    provide two different representations of an essentially the same object named the perturbation determinant $D(z)$   or  the Szeg\H{o} function $S (\z)$. In view of \e{eq:SSF} the first of them is given 
  in terms of $\arg  D(\lambda+i0)$ while according to \e{eq:SF1} and \e{eq:HH} the second representation is stated  
  in terms of  $\ln |D(\lambda+i0)|$. Obviously, these two functions are  harmonic conjugate.
  
  Below it will be convenient to make the change of variables $\lambda=\cos \theta$ in integral  \e{eq:SzY}. 
Recall that  the weight function  $\tau_{0}(\lambda)$ of the Jacobi operator $J_{0}$ is given by formula \e{eq:ww}. Taking  into account that the function $w(\cos\theta)|\sin\theta|$ is even and using relations  \e{eq:wwz}, we see that
    \begin{equation}
S (\z)= \frac{1-\z^2}{\sqrt{ 2\pi}  }\exp\Big(\frac{1-\z^2}{2\pi}\int_{-1}^1  \frac{\ln \big(\tau (\lambda)/\tau_{0}(\lambda)\big)}{1-2\z \lambda+ \z^2}\frac{d\lambda}{\sqrt{1-\lambda^2}}\Big).
\label{eq:SzD1}\end{equation}
In terms of the variable $z=2^{-1}(\z+\z^{-1})\in\Pi_{0}$, 
the integral in \e{eq:SzD1} becomes the Cauchy integral which yields the representation
  \[
S (\z (z))= \z (z) \sqrt{\frac{2(z^{2}-1)}{  \pi}  }\exp\Big(-\frac{\sqrt{z^2-1}}{2\pi}\int_{-1}^1  \frac{\ln \big(\tau (\lambda)/\tau _{0}(\lambda)\big)}{ \lambda-z}\frac{d\lambda}{\sqrt{1-\lambda^2}}\Big).
\]

  \subsection{Trace identities}
 
 We discuss two types of trace identities. Both of them are obtained by studying an asymptotic behavior of the perturbation determinant $D(z)$ (more precisely, of  $\ln D(z)$)
 as $z\to \infty$ or, equivalently, $\z\to 0$. We will get the first set of identities  by considering an asymptotic expansion in powers of $z^{-1}$ while  identities of the second type known as the Case sum rules (see \cite{Case}) are derived by expanding in powers of $\z$.  Note that for the  Schr\"odinger operators, trace identities  were obtained by Buslaev and Faddeev in \cite{BF}; see, e.g.,  \S 4.6 in the book \cite{YA} for details. In  this case the role of the variable \e{eq:ome} is played by $\z=\sqrt{z}$, and 
 there are also two types of identities: of integer and half-integer orders.   The identities of half-integer orders correspond to  the Case sum rules.

 
 To find 
   expressions of $\tr  (  J^{n}-  J_{0}^n) $ in terms of the   spectral shift function, we only have to compare asymptotic expansions as $|z|\to\infty$ of both sides of  representation   \e{eq:DD}.   
It follows from relation
   \e{eq:PD} that
     \[
\ln  D (z)=-\sum_{n=1}^{\infty} n^{-1} \tr  (  J^{n}-  J_{0}^n)\, z^{-n}.
\]
Since
\[
 \int_{-\infty}^\infty \xi(\lambda) (\lambda-z)^{-1} d\lambda= -\sum_{n=1}^{\infty}  \int_{-\infty}^\infty \xi(\lambda) \lambda^{n-1} d\lambda \,  z^{-n}, 
\]
equating the coefficients at $ z^{-n}$, we see that
     \begin{equation}
  \tr  (  J^{n}-  J_{0}^n)=n   \int_{-\infty}^\infty \xi(\lambda) \lambda^{n-1} d\lambda .
\label{eq:BF1}\end{equation}

On the discrete spectrum,  the   spectral shift function  can be explicitly calculated. Indeed,
let $\lambda_{1}^{(+)}>\lambda_{2}^{(+)} >\cdots> 1$ (and $\lambda_{1}^{(-)}<\lambda_{2}^{(-)} <\cdots < - 1$) be eigenvalues of the operator $J$ lying above the point $1$ (respectively, below the point $-1$). It follows from formula  \e{eq:ES} that
$\xi(\lambda)=n$ for $\lambda\in (\lambda_{n+1}^{(+)},\lambda_{n}^{(+)} )$  and $\xi(\lambda)=-n$ for $\lambda\in (\lambda_{n }^{(-)},\lambda_{n+1}^{(-)} )$ whence 
\[
 n\int_{1}^\infty \xi(\lambda) \lambda^{n-1} d\lambda =\sum_{k=1}^{\infty} k \big((\lambda_{k}^{(+)})^{n }-(\lambda_{k+1}^{(+)})^{n }\big)
 \]
 and, similarly, for the integral over $(-\infty,-1)$. The series here is convergent by virtue of  estimate  \e{eq:DS}. Putting together this relation with \e{eq:BF1}, we obtain the following result.
 
  \begin{theorem}\label{BF}
Let assumption   \e{eq:Tr} be satisfied. Then
\begin{equation}
  \tr  (  J^{n}-  J_{0}^n)=n   \int_{-1}^1 \xi(\lambda) \lambda^{n-1} d\lambda +\sum_{k=1}^{\infty}k \big((\lambda_{k}^{(+)})^{n }-(\lambda_{k+1}^{(+)})^{n }\big) +\sum_{k=1}^{\infty}k \big((\lambda_{k}^{(-)})^{n }-(\lambda_{k+1}^{(-)})^{n }\big).
\label{eq:BF2}\end{equation}
\end{theorem}

In view of \e{eq:SSF3}, the integral on the right can be expressed in terms of the phase function:
\[
\int_{-1}^1 \xi(\lambda) \lambda^{n-1} d\lambda=\frac{1}{\pi} \int_0^\pi \eta(\theta) \cos^{n-1}\theta \sin \theta d\lambda.
\]

    From a somewhat different point of view, formulas of type \e{eq:BF2} were studied in the book \cite{Teschl}, Chapter~6.

   The trace formula of zero order (the Levinson theorem) requires a special discussion. Now we assume a stronger condition \e{eq:Trfd} on the coefficients of the operator $J$. Then according to Theorem~\ref{JOSF}, the corresponding perturbation determinant  $D (z)$ is continuous as $z\to \pm 1$. One has to distinguish the cases  $ D(1)=0$ and/or  $D (-1)=0$ when the operator $J $ has threshold resonances at the points $\lambda=1$ and/or $\lambda=-1$. Note  (see Remark~\ref{find}) that under assumption \e{eq:Trfd} the operator $J$ has only a finite number $N$ of discrete eigenvalues.

    \begin{theorem}[\cite{Y/LD}, Theorem~5.10] \label{Levinson}
Let assumption   \e{eq:Trfd} be satisfied. Then the limits $\xi(1-0)$ and $\xi(-1+0)$ exist and
\begin{equation}
 \xi(1-0)-\xi(-1+0)=N+ p
\label{eq:Lev}\end{equation}
where $p=0$ if $D (\pm 1)\neq 0$ for both signs, $p=1/2$ if $D (\pm 1)=0$ for one of the signs and $p=1$ if $D  (\pm 1)=0$ for both   signs.
\end{theorem}

 Let us, finally, obtain the Case sum rules for the pair $J_{0}$, $J$.
  Putting together relations \e{eq:SzD3} and \e{eq:SzD1}, we see that
   \begin{equation}
\ln  \Delta (\z) -\ln B(\z) -\ln A= \frac{\z^2-1}{2\pi}\int_{-1}^1  \frac{\ln \big(\tau (\lambda)/\tau _{0}(\lambda)\big)}{1-2\z \lambda+ \z^2}\frac{d\lambda}{\sqrt{1-\lambda^2}} .
\label{eq:trace}\end{equation}
First we set here $\z=0$ and recall that $ \Delta (0) =1$.  Therefore in view of definitions  \e{eq:AA} and \e{eq:Bla}, relation \e{eq:trace} implies
 the identity
  \begin{equation}
 \sum_{k=1}^\infty \ln (2a_{k}) + \sum_{k=1}^\infty \ln |\mu_{k}| =  \frac{1}{2\pi}\int_{-1}^1   \ln \big(\tau (\lambda)/ \tau_{0}(\lambda)\big) \frac{d\lambda}{\sqrt{1-\lambda^2}} 
\label{eq:case}\end{equation}
known as the Case sum rule of zero order. It is of course quite different from the Levinson theorem \e{eq:Lev}.

More generally,  we consider the asymptotic expansions of both sides of \e{eq:trace} as $\z\to 0$ and compare the coefficients at the same powers of $\z$. According to Theorem~2.13 in \cite{KS} we have
\[
\ln  \Delta (\z) =-2 \sum_{n=1}^\infty n^{-1} \tr \big( T_{n} (J)- T_{n} (J_{0})\big) \z^n 
\]
where $T_{n}(\lambda)=\cos (n\arccos\lambda)$ are the Chebyshev polynomials of the first kind.
It directly follows from definition \e{eq:Bla} that
     \[
\ln B(\z) =\sum_{k=1}^\infty \ln |\mu_k| + \sum_{n=1}^{\infty}n^{-1} \sum_{k=1}^{\infty} (\mu_k^{n}-\mu_k^{-n}) \z^{n}
\]
where the series over $k $ are convergent due to the condition  \e{eq:Bla1}.
Finally, we use formula (10.11.29) in \cite{BE}: 
  \[
  \frac{1-\z^2} {1-2\z \lambda+ \z^2} =1+2 \sum_{n=1}^\infty T_{n} (\lambda) \z^n. 
\]
Thus the equality of the coefficients at $\z^{n}$ in the left-  and right-hand sides of \e{eq:trace} yields the identity
\begin{multline}
\tr \big( J_{n} (H)- J_{n} (H_{0})\big) =-\frac{1}{2} \sum_{k=1}^\infty  (\mu_k^{n}-\mu_k^{-n})
\\+ \frac{n}{2\pi}
\int_{-1}^1   \ln \big( \tau (\lambda)/ \tau_{0}(\lambda)\big) T_{n} (\lambda) \frac{d\lambda}{\sqrt{1-\lambda^2}},\q n=1,2,\ldots. 
\label{eq:case1}\end{multline}

The  trace identities \e{eq:case} and \e{eq:case1} known as the Case sum rules are not new. They    were obtained by him  in  \cite{Case} and rigorously proven in \cite{KS}.
We note that, in the paper \cite{KS},  the identities \e{eq:case} and \e{eq:case1}  were first checked for finite rank perturbations  $J-J_{0}$ and then  \e{eq:case} and \e{eq:case1} (for $n=1$) were used for
  the proof of the inclusion $ \Delta \in {\bf N}$.

  \appendix
 
 \section{Favard's theorem}
 
 Here we prove elementary assertions stated in Sect.~1.1 and 1.2.
 
 {\it Proof of Proposition~\ref{Fav}.}  Since $\lambda P_{n}(\lambda)$ is a polynomial of degree $n + 1$, we have
  \begin{equation}
 \lambda P_{n}(\lambda)= c_{n, n+1} P_{n+1}(\lambda)+ c_{n, n} P_{n}(\lambda)+\cdots +
 c_{n, 1} P_{1}(\lambda) +  c_{n, 0}P _{0}(\lambda).
\label{eq:FAV1}\end{equation} 
 Taking the scalar product  in the space $L^2 ({\Bbb R}; d\rho)$ of this expression with the polynomial $P_{l}$, we find that
 \begin{equation}
\int_{-\infty}^\infty  \lambda P_n (\lambda) P_l (\lambda) d\rho (\lambda)= \sum_{m=0}^{n+1}c_{n, m}\la P_{m}, P_l\ra.
\label{eq:FAV2}\end{equation}  
Let    $l=0, 1, \ldots, n- 2$.
Then the left-hand side here is zero
 because $P_{n}$ is orthogonal in the space $L^2 ({\Bbb R}; d\rho)$ to all polynomials of degree $\leq n-1$.  
 Since the right-hand side of \e{eq:FAV2} equals $c_{n,l}$, we see that $c_{n,l}=0$ for all $l=0, 1, \ldots, n- 2$.
  Therefore it follows from \e{eq:FAV1} that
  \begin{equation}
 \lambda P_{n}(\lambda)= c_{n, n+1} P_{n+1}(\lambda)+ c_{n, n}P_{n}(\lambda) +
 c_{n, n-1}P_{n-1}(\lambda)  .
\label{eq:FAV3}\end{equation} 
 Taking  the scalar product    of this expression with $P_{n}$, we see that $c_{n, n}=: b_{n}$ is given by  \e{eq:RP1}. Then we set $a_{n}:= c_{n, n+1}$ and take the scalar product    of \e{eq:FAV3} with $P_{n+1}$ whence $a_{n}$ is given by  \e{eq:RP}.  Let us now compare the coefficients at $\lambda^{n+1}$ in \e{eq:FAV3}.  According to \e{eq:RRx} we have
  $k_{n}= a_{n}k_{n+1}$ whence $a_{n}>0$ because $k_{n}>0$ for all $n$. Finally, taking the scalar product    of \e{eq:FAV3} with $P_{n-1}$, we see that $ c_{n, n-1}= a_{n-1}$.  Thus  \e{eq:FAV3} yields relation  \e{eq:RR}. $\Box$
 
 \medskip
  
       {\it Proof of Proposition~\ref{adj}.}
       If $f\in{\cal D}$ and $g\in \ell^2 ({\Bbb Z}_{+})$, then according to \e{eq:ZP+1} we have
         \begin{align}
\la J_{\rm min}f, g \ra  =  &\sum_{n=1}^\infty a_{n-1} f_{n-1}\ov{g_n}+ \sum_{n=0}^\infty b_{n} f_{n }\ov{g_n}
\nonumber\\
= & \sum_{n=1}^\infty a_{n-1} f_{n-1}\ov{g_n}+ \sum_{n=0}^\infty b_{n} f_{n }\ov{g_n}
+ \sum_{n=0}^\infty a_{n} f_{n +1}\ov{g_n}+ \sum_{n=0}^\infty a_{n} f_{n +1}\ov{g_n}
\nonumber\\
= & \sum_{n=0}^\infty f_{n}   \ov{({\cal J}g)_{n}} .
\label{eq:ad2}\end{align}
If $g\in{\cal D} (J_{\rm max})$, the right-hand side here equals $\la f, J_{\rm max} g \ra$
 whence $J_{\rm max} \subset J_{\rm min}^*$. 
 
 Conversely, suppose that $g\in{\cal D}(J_{\rm min}^*)$, that is,  $\la J_{\rm min} f,  g \ra =\la f,  g_{*}\ra $ for some $g_{*}\in \ell^2 ({\Bbb Z}_{+})$ and all $f\in {\cal D}$. 
 Comparing this equality with \e{eq:ad2}, we see that $g_{*}= {\cal J} g$ whence $g\in {\cal D} (J_{\rm max})$ and $J_{\rm min}^*\subset J_{\rm max}$. $\Box$

\medskip
 
 {\it Proof of Proposition~\ref{simple}.}
 Relation \e{eq:Phi3} is obvious for $n=0$ when $P_{0} (J_{\rm min})=I$.  Supposing that \e{eq:Phi3} is true for all $n\leq N$, we will verify it for $n=N+1$. By definition \e{eq:RR}, we have
  \begin{align*}
a_{N} P_{N+1} (J_{\rm min})e_{0}= &   (J_{\rm min}-b_{N})P_{N} (J_{\rm min}) e_{0}- a_{N-1}P_{N-1} (J_{\rm min}) e_{0}
\\
= & (J_{\rm min}-b_{N})e_{N} - a_{N-1}e_{N-1}  
\end{align*}
where we have used \e{eq:Phi3} for $n=N$ and $n=N-1$. By definition \e{eq:ZP+1} of the operator $J_{\rm min}$, the right-hand side here equals $a_{N}e_{N+1}$.  $\Box$

\medskip

 {\it Proof of Proposition~\ref{FAV}.} Using \e{eq:Phi3} and \e{eq:UD1}, we see that
  \[
 d\la E_{J}(\lambda)e_{n}, e_{m}\ra = d\la E_{J}(\lambda)P_{n} (J)e_{0}, P_{m} (J) e_{0}\ra = P_{n} (\lambda) P_{m} (\lambda)  d\rho_{J}(\lambda) .
\]
Integrating this equality, we obtain relation \e{eq:mo1}.
It now follows from definition \e{eq:Phi1} that
\[
\la \Phi e_{n},\Phi e_{m}\ra= \d_{n,m}
\]
and hence
 operator \e{eq:Phi} is isometric.
 
 Next, we check the intertwining property \e{eq:Phi2}. It suffices to consider $f=e_{n}$. By definition of $J$, in this case relation  \e{eq:Phi2} means that
 \[
\big( \Phi (a_{n-1} e_{n-1}+ b_{n} e_{n}+ a_{n} e_{n+1})\big) (\lambda)= \lambda (\Phi e_{n} )(\lambda).
\]
In view of definition \e{eq:Phi1} of $\Phi$ this equality is equivalent to  relation 
  \e{eq:RR} defining  $P_{n} (\lambda)$.

It remains to verify that the operator $\Phi$ is unitary, that is, its range $\ran \Phi= L^2 ({\Bbb R}; d\rho_{J})$. Supposing the contrary, we find a vector $g\neq 0$ in $L^2 ({\Bbb R}; d\rho_{J})$ such that $\la \Phi f, g\ra=0$ for all $f\in \ell^2 ({\Bbb Z}_{+})$. In particular, this is true for all elements $f= E_{J} (\Lambda)e_{0}$ where $\Lambda$ is an arbitrary Borelian subset of $\Bbb R$. Using the intertwining property \e{eq:Phi2} and equality $\Phi e_{0}=1$, we find that  $(\Phi E_{J} (\Lambda)e_{0})(\lambda)$ is the characteristic function of $\Lambda$.  Therefore
\[
\la \Phi E_{J} (\Lambda)e_{0}, g\ra=\int_{\Lambda}g(\lambda) d \rho_{J}(\lambda).
\]
Since this integral is zero for all $\Lambda\subset {\Bbb R}$, we see that $g(\lambda)=0$ on a set of full $d \rho_{J}$ measure. 

It follows from \e{eq:Phi3} and \e{eq:Phi2} that
\[
(\Phi e_{n})(\lambda)= (\Phi P_{n}(J) e_{0}) (\lambda)=P_{n} (\lambda).
\]
Since $  e_{0}, e_{1}, \ldots$ is the basis in $\ell^2 ({\Bbb Z}_{+})$ and   operator \e{eq:Phi} is unitary,
  linear combinations of the polynomials $P_{n} (\lambda)$ are dense in the space $L^2 ({\Bbb R}; d\rho_{J})$.
$\Box$

        \section{Classical polynomials}
        
        Here we discuss some basic properties of the Jacobi, Hermite and Laguerre
        polynomials. Note that for the corresponding recurrence coefficients $a_{n}$, the Carleman condition \e{eq:Carl}   is satisfied.
            
    It is convenient to define    {\bf Jacobi} polynomials    by their spectral  measures $d\rho (\lambda)= d\rho^{(\alpha,\beta)} (\lambda)$ where the parameters $\alpha ,\beta>-1$. These measures are supported on the interval $ [-1,1 ]$, $\rho (\{\pm 1 \})=0$ and  
 \begin{equation}
d\rho (\lambda)= k (1-\lambda)^{\alpha} (1+\lambda)^{\beta} d\lambda,   \q  \lambda\in (-1,1) , \q \alpha ,\beta>-1.
\label{eq:Jac1}\end{equation}
The dependence of various objects  on $\alpha $ and $\beta$ is often omitted in notation. 
 The constant 
$  k = k^{(\alpha ,\beta)}$
 is chosen in such a way that the measure \e{eq:Jac1} is normalized, i.e., 
$\rho ({\Bbb R})=\rho ((-1,1))=1$. 
The orthonormal polynomials ${\sf G}_{n}(z)= {\sf G}_{n}^{(\alpha,\beta)}(z)$ determined by  measure
\e{eq:Jac1} are known as the Jacobi polynomials.  According to Proposition~\ref{refl} we have
 \[
 {\sf G}_{n}^{( \beta,\alpha)}(z) =  (-1)^{n} {\sf G}_{n}^{(  \alpha,\beta)}(-z).
\]
  In some particular but important cases,  the polynomials ${\sf G}_{n}^{(\alpha,\beta)}(z)$ have special names. They are called Gegenbauer polynomials if $\alpha=\beta$. In particular,  ${\sf G}_{n}^{(0,0)}(z)$ are known as the Legendre polynomials; ${\sf G}_{n}^{(-1/2,- 1/2)}(z)$ and ${\sf G}_{n}^{(1/2,1/2)}(z)$ are the Chebyshev polynomials of the first and second kinds, respectively.

Let ${\sf J}={\sf J}^{(\alpha,\beta)}$ be the Jacobi operator  with the spectral measure $   d\rho^{(\alpha ,\beta)}(\lambda)$. Explicit expressions for
its matrix elements  $a_{n}, b_{n}$ can be found, for example, in the books \cite{BE, Sz}, but we do not need them. We here note only  asymptotic  formulas 
 \begin{equation}
 a_{n}  =1/2 + 2^{-4} (1 -2\alpha^2-2\beta^2) n^{-2} +O\big( n^{-3}\big),\q
b_{n}  = 2^{-2}  (\beta^2-\alpha^2) n^{-2} +O\big( n^{-3}\big) 
\label{eq:norm6}\end{equation}
for  the matrix elements. In the case $\alpha=\beta=1/2$ we have $a_{n}=1/2$ and $b_{n}= 0$ for all $n$. The corresponding Jacobi operator denoted ${\sf J}^{(1/2,1/2)}=:J_{0}$
 is known as the  ``free"  discrete Schr\"odinger operator.  
Eigenvectors  of $J_{0}$  are   normalized  Chebyshev polynomials     of the second kind, and the corresponding spectral measure $d\rho_{0} (\lambda)= d\la E_{0} (\lambda) e_{0}, e_{0}\ra $ is given by the formula
 \begin{equation}
d\rho_{0} (\lambda)= 2 \pi^{-1}  \sqrt{1-\lambda^2 }\, d\lambda,\q \lambda\in (-1,1). 
\label{eq:fr}\end{equation}  
In the case $\alpha=\beta=-1/2$ we have  $b_{n}= 0$ for all $n$, but $a_{n}=1/2$ for $n\geq 1$ only while  $a_{0}=1/\sqrt{2}$. The corresponding Jacobi operator denoted ${\sf J}^{(-1/2,-1/2)}$ is a two-rank perturbation of the
  operator $J_{0}$.  
Its eigenvectors   are   normalized  Chebyshev polynomials     of the first second kind.

The Jacobi polynomials satisfy asymptotic relations 
  \begin{multline}
 {\sf G}_{n} (\lambda)=  2^{1/2}( \pi k)^{-1/2}(1-\lambda)^{-(1+2\alpha)/4}(1+\lambda)^{-(1+2\beta)/4} 
 \\
 \times \cos \big((n+ \gamma) \arcsin\lambda - \pi(2n+\beta-\alpha) /4\big) + O (n^{-1}), \q \gamma = (\alpha+\beta+1)/2,
\label{eq:GG} \end{multline}
(see formula (8.21.10)   in the book~\cite{Sz}) if $\lambda\in (-1,1)$ and 
 \begin{align}
 {\sf G}_{n} (z)=   ( 2\pi & k)^{-1/2}  2^{- (\alpha+\beta)/2}
(z-1)^{-(1+2\alpha)/4}(z+1)^{-(1+2\beta)/4}
\nonumber\\
\times &\big(\sqrt{z+1} + \sqrt{z-1}\big)^{\alpha+\beta} \big( z+ \sqrt{z^2-1}\big)^{n+1/2}  \big( 1+ o(1)\big)
\label{eq:GGcomp} \end{align}
(see formula (8.21.9)   in  \cite{Sz}) if $z\in {\Bbb C}\setminus [-1,1]$.  Here $\sqrt{z\pm 1}>0$ if $z\pm 1 >0$.
  Estimates of the remainders in \e{eq:GG} and \e{eq:GGcomp}  are uniform in $\lambda$ and $z$ from compact subsets of $(-1,1)$ and of ${\Bbb C}\setminus [-1,1]$, respectively. At the edge points of the spectrum,  the asymptotics of the Jacobi polynomials as $n\to\infty$ are given by the formulas 
   \begin{equation}
      \begin{split}
 {\sf G}_{n} (1)=  & k^{-1/2} 2^{-\alpha-\beta} \Gamma(\alpha+1)^{-1} n^{\alpha+1/2} \big(1+O (n^{-1})\big),
\\
 {\sf G}_{n} (-1)= & (-1)^{n}  k^{-1/2} 2^{-\alpha-\beta} \Gamma(\beta+1)^{-1} n^{\beta+1/2} \big(1+O (n^{-1})\big).
    \end{split}
    \label{eq:edge}\end{equation}

  The {\bf Hermite} polynomials ${\sf H}_{n}  (z)$ are defined by relations \e{eq:RR}, \e{eq:RR1} with the recurrence coefficients 
   \begin{equation}
a_{n} = \sqrt{(n+1)/2}, \q b_{n}  =0.
\label{eq:H}\end{equation} 
      According to Theorem~8.22.7   in the book \cite{Sz}  asymptotics of ${\sf H}_{n}  (z)$ as $n\to\infty$ is given by the Plancherel-Rotach formula
      \begin{align}
{\sf H}_{n}  (z)=
   2^{1/2}
  \pi^{-1/4} e^{z^2/2}  (2n+1)^{-1/4} &\cos \big( \sqrt{2n+1} \, z-  \pi  n/2\big) 
\nonumber  \\
  + &O(e^{\sqrt{2n+1}|\Im z|}n^{-3/4}).  
\label{eq:H2}\end{align}
This  asymptotics   is uniform in $z $  from compact subsets of $ \Bbb C$.  Obviously, the right-hand side of \e{eq:H2} exponentially  grows as $n\to\infty$ if $\Im z\neq 0$, and it is an oscillating function if $ z\in  {\Bbb R}$.


  Let us  consider the Jacobi operator $ J $  with  coefficients \e{eq:H}.       The spectral measure of $J$ equals  
   \[
    d\rho  (\lambda)=   \pi^{-1/2} e^{-\lambda^2} d\lambda, \q  \lambda\in {\Bbb R}
    \]
 (see, e.g., formula (10.13.1) in \cite{BE}).
Thus, $d\rho  (\lambda)$ is absolutely continuous and its 
    support is the whole axis  $\Bbb R$.

   Suppose now that  the recurrence coefficients $a_{n}$, $b_{n}$ are given  by formulas 
    \[
    a_{n} =  a_{n}^{(p)}= \sqrt{(n+1)(n+1+p)} \q\mbox{and}\q     b_{n} = b_{n}^{(p)}=  2n+p+1, \q p>-1. 
\]
The corresponding Jacobi operator   $J=J^{(p)}$      has the absolutely continuous  spectrum
coinciding with $[0,\infty)$, and the spectral measure is given by the relation (see, e.g., formula (10.12.1) in \cite{BE}) 
         \[
d\rho^{(p)} (\lambda)=  \tau^{(p)}(\lambda)d\lambda  \q \mbox{where}\q \tau^{(p)} (\lambda)=\Gamma(p+1)^{-1}\lambda^p e^{-\lambda}, \q \lambda\in {\Bbb R}_{+}.
\]
The eigenfunctions of $J^{(p)}$  are
 orthonormal  {\bf  Laguerre } polynomials ${\sf L}_{n}^{(p)} (z)$ defined by     relations  \e{eq:RR} and \e{eq:RR1}.  
  Note that the normalized polynomials ${\sf L}_{n}^{(p)} (z)$ we consider here are related to the Laguerre polynomials ${\bf L}_{n}^{(p)} (z)$ defined in \S 10.12 of the book \cite{BE} or in \S 5.1 of the book \cite{Sz} by the equality
\[
{\sf L}_{n}^{(p)} (z)= (-1)^n   \sqrt{\frac{\Gamma (1+n) \Gamma (1+p)}
{\Gamma (1+n+p)  }} \:{\bf L}_{n}^{(p)} (z).
\]
According to asymptotic formula (10.15.1) in \cite{BE} for positive $\lambda$, we have
\begin{equation}
{\sf L}_{n}^{(p)} (\lambda)=
(-1)^n   \sqrt{\frac{\Gamma (1+p)  }
 {\pi  }}  \lambda^{-p/2-1/4}  e^{\lambda/2} n^{-1/4} \cos \Big( 2\sqrt{n\lambda}-\frac{2p+1} {4} \pi\Big) +O(n^{-3/4})  
\label{eq:Lag4}\end{equation}
as $n\to\infty$.  For $z\in{\Bbb C}\setminus  [0,\infty)$, one has (see Theorem~8.22.3   in   \cite{Sz})
\begin{equation}
{\sf L}_{n}^{(p)} (z)=
(-1)^n   \sqrt{\frac{\Gamma (1+p)  }
 {\pi  }}  (-z)^{-p/2-1/4}  e^{z /2} n^{-1/4} e^{ 2\sqrt{-n z}}\big(1 +O(n^{-1/2})  \big)
\label{eq:Lag4C}\end{equation}
where $\arg (-z)>0$ if $z<0$.
Asymptotics \e{eq:Lag4} and \e{eq:Lag4C} are uniform in $\lambda $ and $z$ from compact subsets of ${\Bbb R}_{+}$ and ${\Bbb C}\setminus [0,\infty)$,  respectively.
 
 \section{Pollaczek polynomials}
 
 The normalized Pollaczek polynomials are defined (see, e.g., Appendix in the book \cite{Sz}) by recurrent relations \e{eq:RR}, \e{eq:RR1} with
     \begin{equation}
 a_{n}=\frac{n+1}{\sqrt{(2n+2\boldsymbol{\alpha} +1)(2n+2 \boldsymbol{\alpha} +3)}},\q  b_{n}=-\frac{2 \boldsymbol{\beta}}{ 2n+2 \boldsymbol{\alpha} +1 }; 
\label{eq:Poll1}\end{equation}
here the parameters $\boldsymbol{\alpha},\boldsymbol{\beta}\in {\Bbb R}$ and $\boldsymbol{\alpha}> |\boldsymbol{\beta}|$.
It follows   that
    \begin{equation}
 a_{n}=2^{-1} -\boldsymbol{\alpha} (2n)^{-1}+ O (n^{-2}),\q  b_{n}=- \boldsymbol{\beta} n^{-1}  +  O (n^{-2}) \q \mbox{as} \q n\to\infty.
\label{eq:PolX}\end{equation}
The spectrum of the corresponding Jacobi operator coincides with the interval $[-1,1]$, it is absolutely continuous and the  weight function is given by the formula
 \[
\tau (\lambda)= (\alpha+ 1/2) e^{(2\theta-\pi) \Xi (\theta)}\big( \cosh (\pi \Xi(\theta) \big)^{-1}\q {\rm where}\q  \Xi (\theta) =(\boldsymbol{\alpha} \cos\theta +\boldsymbol{\beta}) (\sin\theta)^{-1}
\]
and as usual $\lambda=\cos\theta$.
It is easy to see that
  \begin{equation}
\ln \tau(\lambda)= - \pi (\boldsymbol{\alpha}+ \boldsymbol{\beta})\theta^{-1}+ O(1) 
\label{eq:PolCou}\end{equation}
as $\lambda\to1-0$ and a similar formula is true as $\lambda\to -1+ 0$.

It follows from \e{eq:PolX} that $V=J-J_{0}$ is Hilbert-Schmidt, but the series $\sum_{n} (a_{n}-1/2)$ and $\sum_{n} b_{n}$ are divergent; in particular, assumption \e{eq:Tr} is not satisfied.  
According to \e{eq:PolCou} the Szeg\H{o} condition \e{eq:Szeg} where $s=-1/2$ is violated for Pollaczek polynomials. This is consistent with the classical theorem of Szeg\H{o}, Shohat, Geronimus, Kre\u{\i}n  and Kolmogorov; see, e.g., Theorem~4 in \cite{KS}. On the other hand, relation \e{eq:PolCou}  implies  condition \e{eq:Szeg}
for $s=1/2$
which is consistent with  Theorem~1 in \cite{KS} stated here as Theorem~\ref{KS-qS}.

Pursuing an analogy with differential operators, we note that the Jacobi operators with coefficients \e{eq:Poll1} correspond to the Schr\"odinger operators with Coulomb    potentials (see \S\S 36 and 133 in the book \cite{LL}).  Condition $\boldsymbol{\alpha}> |\boldsymbol{\beta}|$ distinguishes repulsive potentials, and formula \e{eq:PolCou} corresponds to the exponential decay as $\lambda\to 0$ of the corresponding weight function at low energies. If $\boldsymbol{\alpha}< |\boldsymbol{\beta}|$, then the infinite discrete spectrum appears which is also quite similar to the Schr\"odinger operators with Coulomb  attractive  potentials; see Sect.~5.4 and 5.5 in the book \cite{Ism}.

 \section{Moment problems}
 
 For a sequence $s_{0}, s_{1}, \ldots, s_{n}, \ldots$ of positive numbers such that $s_{0}=1$, we set 
  \begin{equation}
{\cal S}_{n} = 
\begin{pmatrix}
 s_{0}&s_{1}& \ldots& s_{n}\\
 s_{1}&s_{2}& \ldots& s_{n+1} \\
   \vdots&  \vdots& \cdots&   \vdots
      \\  s_{n-1}&s_{n}& \ldots& s_{2n-1} 
   \\  s_{n}&s_{n+1}& \ldots& s_{2n} 
\end{pmatrix}
\q \mbox{and} \q
{\cal P}_{n} (z)= 
\begin{pmatrix}
 s_{0}&s_{1}& \ldots& s_{n}\\
 s_{1}&s_{2}& \ldots& s_{n+1} \\
   \vdots&  \vdots& \cdots&   \vdots
   \\  s_{n-1}&s_{n}& \ldots& s_{2n-1} 
   \\ 1& z& \ldots& z^n
\end{pmatrix}.
\label{eq:MOM}\end{equation}
One calls  $s_{0}, s_{1}, \ldots, s_{n}, \ldots$ a moment sequence if
 $ \det {\cal S}_{n}>0$ for all $n\geq 1$. By the Hamburger theorem for an arbitrary moment sequence,
  there exist  measures $d\rho(\lambda)$  with infinite supports satisfying relations \e{eq:mo} for all $n\in{\Bbb Z}_{+}$.  Such measures are in general not unique. However the polynomials $ P_{n}(z)$ satisfying conditions \e{eq:mo1}  are the same for all these measures. They can be constructed by a formula
 \begin{equation}
 P_{n}(z)= (d_{n-1} d_{n})^{-1/2} \det {\cal P}_{n} (z)
\label{eq:MOM1}\end{equation}
where $d_{n}= \det {\cal S}_{n}$.

According to   Proposition~\ref{Fav}
 the   polynomials $P_{n}(z)$ satisfy   recurrence relation 
 \e{eq:RR} with the coefficients $a_{n}$, $b_{n}$ given by
 \e{eq:RP}  and \e{eq:RP1}. Let ${\cal J}$ be the Jacobi operator \e{eq:ZP+1} with  these recurrence coefficients. Then the moments $s_{n}$ can be recovered by the formula 
 \begin{equation}
 s_{n}= \la {\cal J}^n e_{0}, e_{0}\ra .
\label{eq:MOM2}\end{equation}
Indeed, let $J$ be an arbitrary self-adjoint extension of the operator $J_{\rm min}$. 
The right-hand side of \e{eq:MOM2} equals
\[
\la J^n e_{0}, e_{0} \ra=\int_{-\infty}^\infty
\lambda^n d\rho_{J} (\lambda)= s_{n}
\]
by the spectral theorem and relation \e{eq:mo}.
Conversely, let recurrence coefficients $a_{n}$, $b_{n}$ be given. Define 
  the numbers $s_{n}$   by formula \e{eq:MOM2} and the polynomials $P_{n}(z)$  by equalities  \e{eq:MOM},  \e{eq:MOM1}. Then using Proposition~\ref{Fav} we can recover coefficients $a_{n}$, $b_{n}$. Thus the moments $s_{n}$ and the recurrence coefficients $a_{n}$, $b_{n}$ are in a one-to-one correspondence.
 
 Recall that the moment problem \e{eq:mo}  is called determinate if the measure satisfying these relations is unique. Otherwise, it is called indeterminate. It is known (see Theorem~2 in \cite{Simon}) that the determinacy is equivalent to the essential self-adjointness of the operator $J_{\rm min}$.



\begin{thebibliography}{99}
    
     \bibitem{AKH} N.~Akhiezer,
\emph{The classical moment problem and some related questions in analysis},  Oliver and  Boyd, Edinburgh and London, 1965.

 \bibitem{Apt} A.~I.~Aptekarev and J.~S.~Geronimo, Measures for orthogonal polynomials with unbounded recurrence coefficients, J. Approx. Theory {\bf 207}  (2016), 339-347.
 
  
     \bibitem{Atk} F.~V.~Atkinson,
\emph{Discrete and continuous boundary value problems},  Academic Press, New York, 1964.

\bibitem {Ber} Yu. M. Berezanskii,  {\it Expansion in eigenfunctions of selfadjoint operators}, Amer. Math. Soc., Providence, R.I., 1968.


 \bibitem {Bern} S. Bernstein,
   Sur les polyn\^omes orthogonaux relatifs  \`a un segment fini, Journal de Math\'e\-matiques  {\bf 9}  (1930), 127-177;   {\bf 10}  (1931), 219-286.
   
   
 \bibitem {Birkh} G. D. Birkhoff, General  theory of linear difference equations, Trans. Amer. Math. Soc.
  {\bf 12}  (1911), 243-284.
   
   \bibitem {BS}M. Sh. Birman and M. Z. Solomyak, {\em Spectral theory of selfadjoint operators in 
 Hilbert space}, Reidel, Doldrecht, 1987.
 
 \bibitem {BdM} A.  Boutet de Monvel and J. Sahbani,  
Anisotropic Jacobi matrices with absolutely continuous spectrum,
C. R. Acad. Sci. Paris S{\'e}r. I Math. {\bf 328}, No. 5  (1999), 443-448. 


 \bibitem {BF}V. S. Buslaev and L. D. Faddeev,
  {\em Formulas for traces for a singular Sturm-Liouville differential operator}, Soviet Math. Dokl. {\bf 1}  (1960), 451-454.
   
   
  
  
  
\bibitem{Carleman}   T.~Carleman, \emph{ Les fonctions quasi-analytiques}, Gauthier-Villars, 1926.

 




\bibitem {Case1} K.~M.~Case,    Orthogonal polynomials from the viewpoint of scattering theory, 
 Journal of Math. Phys., {\bf 15} (1974),  2166-2174.

\bibitem {Case} K.~M.~Case,    Orthogonal polynomials. II, 
 Journal of Math. Phys., {\bf 16} (1975),  1435-1440.
 
 \bibitem{Chihara} T.~S.~Chihara, \emph{An introduction to orthogonal polynomials}, Gordon and Breach
Science Publishers, New York-London-Paris, 1978.


 



\bibitem {DS1} D.~Damanik and B.~Simon,   Jost functions and Jost solutions for Jacobi matrices, I. A necessary and sufficient condition for Szeg\H{o} asymptotics, Invent. Math. {\bf 165}  (2006), 1-50.

\bibitem {DS} D.~Damanik and B.~Simon,  Jost function, and Jost solutions for Jacobi matrices, II. Decay and
 analyticity, Int. Math. Res. Notices, No.~5 (2006); art. ID 19396, 1-32.
 
   
 
       
       

 
   
  
  
 
 
 
 
 
  





 
\bibitem {Deift}  P.~Deift, {\it  Orthogonal polynomials and random matrices. A Riemann-Hilbert approach}, NYU lectures, AMS, 2000.


  \bibitem{D-Z}     P.~Deift and X.~Zhou,   A steepest descent method for oscillatory Riemann-Hilbert problem,  Ann.
Math., {\bf 137} (1993), 295-368.


    \bibitem{Duren} P.~L.~Duren,  {\em Theory of $H^p$ spaces}, Academic Press, New York and London, 1970.




\bibitem{El} S.~Elaydi, \emph{An introduction to difference equations}, Springer Science+Business Media, New York, USA, 2005.

\bibitem{BE} A.~Erd\'elyi,  W.~Magnus, F.~Oberhettinger, F.~G.~Tricomi, \emph{Higher transcendental functions}, Vol. 1, 2,   McGraw-Hill, New York-Toronto-London, 1953.

\bibitem {F-I-K} A. Fokas, A. Its, A. Kitaev,    The isomonodromy approach to matrix models in $2D$ quantum gravity, 
Comm.
Math. Phys., {\bf 147} (1992), 395-430. 



   \bibitem{M-F} A.~Foulqui\'e~Moreno, A.~Mart\'\i nez-Finkelshtein, and V.~L.~Sousa,   Asymptotics of orthogonal polynomials for a weight with a jump on $[-1,1]$, Constr. Approx. {\bf 33}, No. 2,  (2011),  219-263. 

 
   \bibitem{Freud}  G.~Freud,   On polynomial  approximation with the weight   $\exp(-x^{2k}/2)$, Acta Math. Acad. Sci. Hungar. {\bf 24},  (1973),  363-371. 
   
   \bibitem {Ge-Le} I. M. Gel'fand and B. M. Levitan, {\em  On the determination of a differential equation from its spectral function}, Izv. Akad. Nauk SSSR, Ser. Mat. {\bf 15}  (1951), 309-361; Amer. Math. Soc. Transl. (Ser. 2)  {\bf 1}  (1955), 253-304.
   
     \bibitem{Va-As1}   J.~S.~Geronimo and W.~Van Assche, Orthogonal polynomials with asymptotically periodic recurrence coefficients,  Journal   Appr. Theory {\bf 55} (1988), 220-231.
  
    \bibitem{Va-As}   J.~S.~Geronimo and W.~Van Assche, Asymptotics of the  
 orthogonal polynomials on and off the essential spectrum, Journal   Appr. Theory {\bf 55} (1988), 220-231.
 
  
 
  


  \bibitem {GD} D.~Gilbert and D. B.~Pearson, On subordinacy and analysis of the spectrum of one-dimensional
    Schr\"odinger  operators, J. Math. Anal. Appl. {\bf 128}, no. 2 (1987), 30-56.
    
   

   
   
\bibitem {GK} I. C. Gokhberg and M. G. Kre\u{\i}n, {\em Introduction to the theory of linear
nonselfadjoint operators in Hilbert space}, Amer. Math. Soc., Providence, Rhode Island, 1970.

  
\bibitem {Gon} A. A.~Gon\v{c}ar, On convergence of Pad\'e approximants for some classes of meromorphic functions,   Math. USSR Sb. {\bf 26}  (1975), 555-575.

    



     \bibitem{Hu-S}     D.~Hundertmark, B.~Simon, {\em Lieb-Thirring inequalities for Jacobi matrices},  J.
Approx. Theory, {\bf 118} (2002), 106-130.


     \bibitem{Ism} M. E. H.~Ismail, {\em Classical and quantum orthogonal polynomials in one variable}, Cambridge University Press, Cambridge, 2005.

 \bibitem{J-M}
J.~Janas and  M.~Moszy\'nski, Spectral properties of
Jacobi matrices by asymptotic analysis, J. Approx. Theory {\bf 120}  (2003), 309-336.

  \bibitem{Jan-Nab}
J.~Janas and  S.~Naboko, 
Jacobi matrices with power-like weights -- grouping in blocks approach, Journal of Funct. Analysis  
{\bf 166} (1999), 218-243.

 \bibitem{Jost} R. Jost, {\em \"Uber die falschen Nullstellen der Eigenwerte des $S$-matrix}, Helv. Phys. Acta  {\bf 20} (1947), 250-266.



 \bibitem {KhD} S.~Khan and D. B.~Pearson, Subordinacy and   spectral theory for
 infinite matrices, Helv. Phys. Acta {\bf 65},  (1992), 505-527.
 
   \bibitem{KS}     R.~Killip and B.~Simon,   Sum rules for Jacobi matrices and their applications to spectral theory,  Ann. of Math., {\bf 158} (2003), 253-321.

 \bibitem{Kost}
  A.~G.~Kostyuchenko and K.~A.~Mirzoev, Generalized Jacobi matrices and deficiency indices of differential
operators with polynomial coefficients, Funct. Anal. Appl. {\bf 33}, No. 1 (1999), 38-48.

  \bibitem{Kriech}
T.~Kriecherbauer and K.~T-R~McLaughlin, Strong asymptotics of polynomials orthogonal with respect
to Freud weights, Int. Math. Res. Notices, No. 6 (1999), 299-333.



 
  
 
    
 

 
 \bibitem {LL} L. D.~ Landau and E. M.~Lifshitz, {\it Quantum mechanics}, Pergamon Press, 1965.

   


 \bibitem {Lub} D.~S.~Lubinsky,  A survey of general orthogonal polynomials for weights on finite and infinite intervals,  Acta Appl. Math.  {\bf 10},   (1987),  237-296. 
 
 
   \bibitem{LMS} D.~Lubinsky,  H.~Mhaskar, and E.~Saff,  A proof of Freud's conjecture for exponential weights, Constr. Approx. {\bf 4}   (1988),  65-83. 
   
   

  \bibitem{Magnus} A.~P.~Magnus, On Freud's equations for exponential weights. Papers dedicated to the memory of G\'eza Freud, J. Approx. Theory {\bf 46}  (1986), 65-99.
  
   \bibitem {M-N} A.~M\'at\'e and  P.~Nevai,  Orthogonal polynomials and absolutely continuous measures, In: Approximation Theory IV (C.~K.~Chui, L.~L.~Schumaker, J.~D. Ward, eds.), New York: Academic Press, pp. 611-617.
   
    
 
 \bibitem {Mate} A.~M\'at\'e,  P.~Nevai, and V.~Totik,    Asymptotics for orthogonal polynomials defined by a recurrence relation,  
 Constr. Approx. {\bf 1}   (1985),  231-248. 
    
\bibitem {Mo1}E.~Mourre,  Absence of singular spectrum for certain self-adjoint operators, Comm.
Math. Phys. {\bf 78} (1981), 391-400. 


    \bibitem {Nab} S. N.~Naboko and S. I.~Yakovlev,   Discrete  Schr\"odinger   operator. 
    The point spectrum on the continuous one, Saint-Petersburg Math. Journal {\bf 4}, No. 3 (1993), 559-568.
    
  


 
  
      
    
      



      \bibitem {Nev} P. G.~Nevai,   {\em Orthogonal polynomials}, Memoirs of the AMS {\bf 18}, No. 213, Providence,  R. I., 1979.
  
     \bibitem {Nev3} P.~Nevai,   Orthogonal polynomials defined by a recurrence relation, Trans. Amer. Mah. Soc. {\bf 250} (1979),  369-384.
     
     \bibitem {Nev1} P. G.~Nevai,  Asymptotics for orthogonal polynomials  associated with $\exp(-x^4)$, Siam J. Math. Anal.  {\bf 15}   (1984)  No. 6,  1171-1187. 

   
    \bibitem {Nev2} P.~Nevai,   Orthogonal polynomials,  recurrences, Jacobi matrices, and measures, in   {\em Progress in Approximation Theory} (Tampa, FL, 1990), pp. 79-104,  {\em Springer Ser. Comput. Math.}. {\bf 19}, Springer, New York, 1992.
  
  \bibitem {Nevan} R.~Nevanlinna, Asymptotische Entwickelungen beschr\"ankter Funktionen und das Stieltjessche  Momentenproblem, Ann. Acad. Sci. Fenn. A {\bf 18}, No. 5  (1922), 52 pp.
  
    \bibitem {Nik} E. M.~Nikishin,   Discrete Sturm-Liouville operators and some problems of function theory, J. Sov. Math. {\bf 35}  (1986), 2679-2744.
  
     \bibitem {Olver}  F.  W.  J.~Olver, {\it Introduction to asymptotics and special functions}, Academic Press, 1974.
 
 
   \bibitem {Pos} J.~P\"oschel,   Examples of discrete  Schr\"odinger   operators with 
    pure point spectrum, Comm. Math. Phys. {\bf 88}, No. 3 (1983), 447-463.   

    
  \bibitem {Rah} E.~A. Rakhmanov,   On asymptotic properties of polynomials orthogonal on the real axis,  Math. USSR-Sb. {\bf 47} (1984), 155-193.
  
  
  
  
     \bibitem{Rakh}  E. A. Rakhmanov,   {\it Strong  asymptotics for orthogonal polynomials}, Lecture Notes Math. {\bf 1550}, Springer, Berlin,   1993, 71-97. 
     
  
      \bibitem {Sah}   J.~Sahbani,  Spectral theory of certain unbounded Jacobi matricrs, J. Math. Anal. Appl. {\bf  342}  (2008), 663-681.
  
   \bibitem{Schm} K.~Schm\"udgen, \emph{The moment problem}, Graduate Texts in Mathematics, Springer,   2017.

 
 
 
  
  \bibitem {Sim}B. Simon, {\em Trace ideal methods}, London Math. Soc. Lecture Notes, Cambridge Univ.
Press,  London and New York, 1979.

 
  \bibitem{Simon}
  B.~Simon, 
 The classical moment problem as a self-adjoint finite difference operator, Advances in  Math.  
{\bf 137} (1998), 82-203.
   

 
  

 
     \bibitem {Stolz} G.~Stolz,   Spectral theory for slowly oscillating potentials I. Jacobi matrices, Manuscripta Math.   {\bf 84}   (1994), 245-260.
    
 
  
   \bibitem {Sw-Tr} G.~\'{S}widerski and B.~Trojan, Asymptotics of orthogonal polynomials with slowly oscillating recurrence coefficients,
 J. Funct. Anal. {\bf 278}  (2020), 108326.
 
 
    \bibitem{Sz} G.~Szeg\H{o},  {\em Orthogonal polynomials}, Amer. Math. Soc., Providence, R. I., 1978.
 
 \bibitem {Teschl}  G.~Teschl, {\em Jacobi operators and completely integrable nonlinear lattices}, Amer. Math. Soc.,   Providence,  R. I., 2000.
 
 \bibitem {Tot}V.~Totik,  Orthogonal polynomials, Surveys in Appr. Theory {\bf 1} (2005), 70-125. 
  
     
        \bibitem{Asshe} W.~Van Assche, {\em Asymptotics for  
 orthogonal polynomials},  Lecture Notes in Math. {\bf 1265} Springer-Verlag, Berlin, 1988.
 
   

  


 


 


  
     
     
   \bibitem{Van} M. Vanlessen, Strong asymptotics of the recurrence coefficients
of orthogonal polynomials associated to the
generalized Jacobi weight, Journal   Appr. Theory {\bf 125} (2003), 198-237.

  \bibitem {W-L} G.~Wong and H.~ Li, Asymptotic expansions for second-order linear difference equations,
 J. Comp. Appl. Math. {\bf 41}  (1992), 65-94.


  
    \bibitem {YS}D. R. Yafaev,  The low energy scattering for slowly decreasing 
  potentials,  Comm.
Math. Phys. {\bf 85} (1982), 177-196. 

 \bibitem{Ya} D. R. Yafaev, {\em Mathematical scattering theory: General theory}, Amer. Math. Soc.,   Providence,
  R. I., 1992.

 
  
   \bibitem{YA} D. R. Yafaev, {\em Mathematical scattering theory: Analytic  theory}, Amer. Math. Soc.,   Providence,
  R. I., 2010.
  
  
 
  \bibitem{Y/LD} D. R. Yafaev,   Analytic scattering theory for Jacobi operators and Bernstein-Szeg\H{o}
  asymptotics of orthogonal polynomials,  Rev. Math. Phys. {\bf 30}, No. 8  (2018),  1840019.

   \bibitem{Y-LR} D. R. Yafaev,   A note on the Schr\"odinger operator with a long-range potential,   Letters Math. Phys. {\bf 109}, No. 12  (2019), 2625-2648.  
 
   
    \bibitem{JLR} D. R. Yafaev,  Semiclassical  asymptotic behavior of orthogonal polynomials,    Letters Math. Phys.   {\bf 110}, No. 11  (2020), 2857-2891.  

    
     \bibitem{nCarl} D. R. Yafaev,  Asymptotic behavior of orthogonal polynomials without the Carleman condition,  J. Funct. Anal. {\bf 279}, No. 7   (2020), 108648.




     \bibitem{univ} D. R. Yafaev,   Universal relations in asymptotic formulas for orthogonal polynomials, Funct. Anal. Appl.  {\bf 55}, N 2  (2021), 77-99;  arXiv: 2011.14987 (2020).

 \bibitem{Jacobi-LC} D. R. Yafaev, Self-adjoint Jacobi operators in the limit circle case, J. Oper. Theory (accepted); arXiv  2104.13609.

    


 
 

 
 

 

  
    

       
     
      

  
    
  


 


 
  
 
  

  
 
 
 



 

   

  

   
  

   
 

 
 

 
 




 
 
        
    
 
 

  

  
 

 
         


 

 


  
  
   
 
    
 

 
  
 
       
     


   
       \end{thebibliography}
  \end{document}